\documentstyle[euscript,epsf,amssymb,12pt]{amsart}
\input xy
\xyoption{all}

\newcommand{\la}{\langle}
\newcommand{\ra}{\rangle}

\newtheorem{theorem}{\bf Theorem}[section]
\newtheorem{lemma}[theorem]{\bf Lemma}

\newtheorem{remark}[theorem]{\bf Remark}

\newtheorem{corollary}[theorem]{\bf Corollary}

\renewcommand{\AA}{{\Bbb A}}

\newcommand{\CC}{{\Bbb C}}
\newcommand{\DD}{{\Bbb D}}
\newcommand{\II}{{\Bbb I}}

\newcommand{\NN}{{\Bbb N}}

\newcommand{\QQ}{{\Bbb Q}}
\newcommand{\RR}{{\Bbb R}}

\newcommand{\ZZ}{{\Bbb Z}}

\newcommand{\av}{\operatorname{av}}
\newcommand{\area}{\operatorname{area}}

\newcommand{\hor}{\operatorname{hor}}

\newcommand{\crit}{\operatorname{cr}}

\newcommand{\diam}{\operatorname{diam}}
\newcommand{\dist}{\operatorname{dist}}

\newcommand{\End}{\operatorname{End}}

\newcommand{\gen}{\operatorname{ge}}

\newcommand{\Hol}{\operatorname{Hol}}

\newcommand{\Id}{\operatorname{Id}}

\newcommand{\Isom}{\operatorname{Isom}}
\newcommand{\Ker}{\operatorname{Ker}}
\newcommand{\Lie}{\operatorname{Lie}}
\newcommand{\length}{\operatorname{length}}
\newcommand{\loc}{\operatorname{loc}}

\newcommand{\Res}{\operatorname{Res}}

\newcommand{\st}{\operatorname{st}}

\renewcommand{\vert}{\operatorname{vert}}
\newcommand{\Vol}{\operatorname{Vol}}
\newcommand{\vol}{\operatorname{vol}}
\newcommand{\weight}{\operatorname{weight}}

\renewcommand{\exp}{\operatorname{exp}}

\newcommand{\ov}{\overline}

\newcommand{\curly}{\cal}

\newcommand{\CCC}{{\curly C}}

\newcommand{\HHH}{{\curly H}}

\newcommand{\MMM}{{\cal M}}
\newcommand{\NNN}{{\cal N}}
\newcommand{\OOO}{{\curly O}}
\newcommand{\PPP}{{\curly P}}

\newcommand{\TTT}{{\curly T}}

\newcommand{\XXX}{{\curly X}}
\newcommand{\qu}{/\kern-.7ex/}
\newcommand{\exh}{\to\kern-1.8ex\to}

\newcommand{\VP}{{\curly V}\kern-0.9ex\PPP}

\newcommand{\tV}{\widetilde{V}}

\newcommand{\tphi}{\widetilde{\phi}}

\newcommand{\wt}{\widetilde}

\newcommand{\imag}{{\mathbf i}}

\newcommand{\bq}{{\mathbf{q}}}
\newcommand{\bx}{{\mathbf{x}}}
\newcommand{\by}{{\mathbf{y}}}
\newcommand{\bz}{{\mathbf{z}}}
\newcommand{\bC}{{\mathbf{C}}}
\newcommand{\bP}{{\mathbf{P}}}
\newcommand{\bY}{{\mathbf{Y}}}

\newcommand{\oM}{\overline{\MMM}}
\newcommand{\oN}{\overline{\NNN}}
\newcommand{\Met}{\MMM et}

\newcommand{\YMH}{{\cal Y}{\cal M}{\cal H}}

\newcommand{\ox}{\overline{x}}
\newcommand{\oU}{\overline{U}}
\newcommand{\otau}{\overline{\tau}}

\newcommand{\ZI}{Z_I}
\newcommand{\ZII}{Z_{I\kern-.3ex I}}
\newcommand{\ZIII}{Z_{I\kern-.3ex I\kern-.3ex I}}

\date{31--4--2004}

\author[I. Mundet i Riera]{I. Mundet i Riera}
\address{Departament d'\`Algebra i Geometria, Facultat de Matem\`atiques,
Universitat de Barcelona, Gran Via de les Corts Catalanes 585,
08007 Barcelona, Spain} \email{mundet@@mat.ub.es}

\author[G. Tian]{G. Tian}
\address{Department of Mathematics, Princeton University,
Fine Hall, Washington Road, Princeton NJ 08544-1000 USA}
\email{tian@@math.princeton.edu}

\title{A compactification of the moduli space
of twisted holomorphic maps}

\begin{document}

\maketitle

\begin{abstract}
We construct a compactification of the moduli space of twisted
holomorphic maps with varying complex structure and bounded
energy. For a given compact symplectic manifold $X$ with a
compatible complex structure and a Hamiltonian action of $S^1$
with moment map $\mu:X\to\imag\RR$, the moduli space which we
compactify consists of equivalence classes of tuples
$(C,P,A,\phi)$, where $C$ is a smooth compact complex curve of
fixed genus, $P$ is a principal $S^1$ bundle over $C$, $A$ is a
connection on $P$ and $\phi$ is a section of $P\times_{S^1}X$
satisfying
$$\ov{\partial}_A\phi=0,\qquad \iota_{v}F_A+\mu(\phi)=c,$$
where $F_A$ is the curvature of $A$, $v$ is the restriction on $C$
of a volume form on the universal curve over $\oM_g$ and $c$ is a
fixed constant. Two tuples $(C,P,A,\phi)$ and $(C',P',A',\phi')$
are equivalent if there is a morphism of bundles $\rho:P\to P'$
lifting a biholomorphism $C\to C'$ such that $\rho^*A'=A$ and
$\rho^*\phi'=\phi$. The energy of $(C,P,A,\phi)$ is
$\|F_A\|_{L^2}^2+\|d_A\phi\|_{L^2}^2 +\|\mu(\phi)-c\|_{L^2}^2$,
and the topology of the moduli space is the natural one. We also
incorporate marked points in the picture.

There are two sources of non compactness. First, bubbling off
phenomena, analogous to the one in Gromov--Witten theory. Second,
degeneration of $C$ to nodal curves. In this case, there appears a
phenomenon which is not present in Gromov--Witten: near the nodes,
the section $\phi$ may degenerate to a chain of gradient flow
lines of $-\imag\mu$.
\end{abstract}

 \tableofcontents

\section{Introduction}

\subsection{}
Let $C$ be a compact smooth complex curve with a volume form
$\nu$, and let $X$ be a compact symplectic manifold with a
compatible almost complex structure. Suppose that $X$ supports a
Hamiltonian action of $S^1$ with moment map $\mu:X\to\imag\RR$. A
twisted holomorphic map (see \cite{M}, where the same object was
called twisted holomorphic curve) consists of a principal $S^1$
bundle $P$ over $C$, a connection $A$ on $P$ and a section $\phi$
of the associated bundle $P\times_{S^1}X$, satisfying the
equations
\begin{equation}
\label{eq:the-equations-intro} \ov{\partial}_A\phi=0
\qquad\text{and}\qquad \iota_{v}F_A+\mu(\phi)=c,
\end{equation}
where $F_A$ denotes the curvature of $A$, $v$ is a volume form on
$C$, and $c\in\imag\RR$ is a constant (see Section
\ref{s:symplectic-fibrations} for the definition of
$\ov{\partial}_A\phi$). The second equation is called the vortex
equation. Two triples $(P,A,\phi)$ and $(P',A',\phi')$ are said to
be equivalent if there is an isomorphism of bundles $\rho:P\to P'$
which lifts an automorphism of the curve $C$ and such that
$\rho^*A'=A$ and $\rho^*\phi'=\phi$ (note that in \cite{M} we only
consider isomorphisms of bundles lifting the identity, i.e., gauge
transformations). The set of isomorphism classes of triples
$(P,A,\phi)$ carries a natural topology, and the resulting
topological space $\MMM$ is what we call the moduli space of
twisted holomorphic maps over $C$. There is a notion of energy for
triples of the form $(P,A,\phi)$, called the Yang--Mills--Higgs
functional:
\begin{equation}
\label{eq:YMH-intro}
\YMH_c(P,A,\phi):=\|F_A\|_{L^2}^2+\|d_A\phi\|_{L^2}^2+
\|\mu(\phi)-c\|_{L^2}^2.
\end{equation}
Given a number $K>0$, define $\MMM(K)$ to be the subset of $\MMM$
consisting of isomorphism classes of triples $(P,A,\phi)$ with
energy $\leq K$. The space $\MMM(K)$ is not compact, but the only
source of noncompactness is the bubbling off phenomenon well know
in Gromov--Witten theory. Thus, it is not a surprise that
essentially the same methods as in Gromov--Witten theory (combined
with some standard techniques in gauge theory) allow to construct
a compactification of $\MMM(K)$. This was done in \cite{M}, and
the resulting compactification was used to construct what we call
the Hamiltonian Gromov--Witten invariants of $X$, which are very
analogous to Gromov--Witten invariants but depend not only on a
symplectic structure but also on a Hamiltonian action of $S^1$
(both things up to deformation).

\subsection{}
Now, in the general theory of Gromov--Witten invariants one
considers maps from curves $C$ to $X$ representing a given
homology class, and only the genus of $C$ is fixed: its complex
structure varies along the moduli space of curves of the given
genus. Looking for a compactification of the resulting moduli
space leads to the notion of stable map introduced by Kontsevich
in \cite{Ko}, which consist of a compact curve $C$ with nodal
singularities and with some marked points $\bx=(x_1,\dots,x_n)$,
and a holomorphic map $f:C\to X$. The only condition on
$(C,\bx,f)$ is that any rational component of $C$ on which $f$
restricts to a trivial map has to have at least $3$ exceptional
(that is, marked or singular) points. A definition of
Gromov--Witten invariants using moduli spaces of stable maps
was given independently by several authors
\cite{FO,LT,S,R} for general compact symplectic manifolds,
extending an earlier definition of Gromov--Witten invariants
by Ruan and Tian \cite{RT1,RT2}
for semipositive symplectic manifolds. One of the most important
features of this theory comes from the fact that the moduli space
of stable curves fibres in some sense over the Deligne--Mumford
moduli space of stable marked curves $\oM_{g,n}$. This allows to
introduce into the play the rich structure in the cohomology of
$\oM_{g,n}$. A beautiful and simple example which shows how
powerful this point of view is is the proof of the associativity
of the small quantum product in terms of the geometry of
$\oM_{0,4}$. More generally, the gluing formula or composition
axiom encodes the structure in Gromov--Witten invariants inherited
by the properties of the cohomology of $\oM_{g,n}$.

In view of this, it seems desirable to generalize the construction
given in \cite{M} replacing the fixed curve $C$ by arbitrary
marked stable curves, constructing in this way a compact moduli
space which fibres over $\oM_{g,n}$ and obtaining richer
invariants. The purpose of this paper is to make a first step
towards this aim: we define stable twisted holomorphic maps
($c$-STHM for short), we give a notion of energy of a $c$-STHM
(the Yang--Mills--Higgs functional), we define a topology on the
set of isomorphism classes of $c$-STHM's, and we prove that the
set of isomorphism classes $c$-STHM's with bounded energy is
compact.

\subsection{}
The definition of $c$-STHM involves several new features which do
not appear in the definition of twisted holomorphic maps. The
first one is that on nodal curves one has to consider connections
whose pullback to the normalization is not necessarily smooth, but
only meromorphic. This means that the connection does not extend
to the preimages of the nodes, but its curvature extends
continuously. To motivate this, consider a sequence of smooth
curves $C_u$ converging to a curve $C$ with a nodal singularity,
and choose simple closed curves $\gamma_u\subset C_u$ representing
the vanishing cycle. If we take bundles $P_u$ and connections on
$A_u$, it may perfectly happen that the holonomy of $A_u$ along
$\gamma_u$ (which is an element of $S^1$) converges to an element
different from the identity as $u$ goes to $\infty$. If the
curvature of $A_u$ is uniformly bounded in compact sets disjoint
with $\gamma_u$, passing to a subsequence we obtain a limit bundle
$P$ and connection $A$ defined on the smooth locus of $C$, but $A$
does not extend smoothly in the normalization of $C$ because of
the nontrivial holonomy. (The reason why we consider connections
whose curvature extends continuously will be shortly clarified.)

An important notion related to meromorphic connections is that of
limit holonomy: given a meromorphic connection $A$ on the
punctured disk $\DD^*$ with a pole at the origin, one can define
the (limit) holonomy of $A$ around the origin to be the limit of
the holonomies around circles centered in the origin (and oriented
in a way compatible with the complex structure) as the radius
converges to $0$. This limit exists because the curvature of $A$
extends continuously to $\DD$ (in fact boundedness is enough).
The connection $A$ extends to a connection on the whole disk $\DD$
(of type $C^1$) if and only if the holonomy is trivial (see
Corollary \ref{cor:trivial-hol-extension}).

Let $C'$ be the normalization of $C$
and let $y,y'$ be the preimages of the singular point $z\in C$.
The holonomy of the connection $A$ around $y$
is equal to the inverse of the holonomy around $y'$, because both
holonomies are obtained as the limits of the holonomies of $A_u$
along $\gamma_u$, one using an orientation and the other using
the opposite one.

\subsection{}
Suppose now that we have sections $\phi_u$ of $P_u\times_{S^1}X$
which satisfy equations (\ref{eq:the-equations-intro}) (this
involves a choice of volume form in each curve $C_u$, but for the
moment we will ignore this issue --- see Section
\ref{ss:universal-metrics}). Assume that the norm of the covariant
derivatives $|d_{A_u}\phi_u|$ is uniformly bounded away from the
vanishing cycles. Then the same standard arguments as in \cite{M}
allow to obtain, passing to a subsequence and regauging, a limit
triple $(P,A,\phi)$ defined on the smooth locus of $C$.
Furthermore, one can prove easily that $\|d_A\phi\|_{L^2}$ is
bounded and that $(P,A,\phi)$ satisfies equations
(\ref{eq:the-equations-intro}).

An important question to understand about $(P,A,\phi)$ is the
asymptotic behavior of $\phi$ as we approach the singularity of
$C$. If the holonomy of $A$ around the singular point is trivial,
then $A$ extends the normalization $C'$ (because its curvature is
uniformly bounded), and then Gromov's theorem on removal of
singularities (as proved for continuous almost complex structures
in \cite{IS}) proves that $\phi$ extends also to $C'$. (At this
point elliptic bootstrapping using (\ref{eq:the-equations-intro})
proves that in fact both $A$ and $\phi$ extend smoothly to $C'$.)
To describe what happens when the holonomy is nontrivial, let us
consider the real blow up of $C$ at $y$ and $y'$. This is a real
surface $\wt{C}$ whose boundary consists of two circles $S_y$ and
$S_{y'}$, called exceptional divisors, which can be identified
with circles centered at $0$ in the tangent spaces $T_yC$ and
$T_{y'}C$. The following is a particular case of Corollary
\ref{cor:existence-limit-orbit}.

\begin{theorem}
\label{thm:intro-extensio} Any triple $(P,A,\phi)$ defined over
the smooth locus of $C$ which satisfies
$$\ov{\partial}_A\phi=0\qquad\text{and}\qquad
\|d_A\phi\|_{L^2}<\infty$$ extends to the real blow up $\wt{C}$.
The restriction of the extension of $\phi$ to the exceptional
divisors is covariantly constant, and takes values in the set of
points which are fixed by the action of the limit holonomy of $A$
around $y$.
\end{theorem}

Let $F\subset X$ be the set of fixed points. We say that a
meromorphic connection has critical holonomy $H$ around a pole if
the set of points of $X$ which are fixed by $H$ is bigger than
$F$. The set of critical holonomies forms a finite subset of
$S^1$. Theorem \ref{thm:intro-extensio} implies, in particular,
that if the holonomy of $A$ around $y'$ is not critical then
$\phi$ converges somewhere in the fixed point set as we approach
$z$.
When the holonomy of $A$ is a root of unity (for example, if the
holonomy is critical) Theorem \ref{thm:intro-extensio} follows
from Gromov's theorem on removal of singularities, considering
local coverings ramified at $y$ and $y'$ so that the pullback connection
has trivial holonomy. In fact, the holonomy of $A$ specifies in this
case a structure of orbifold near $y$ and $y'$, and the statement
of the theorem is equivalent to saying that the section $\phi$ extends
to the normalisation $C'$ as a section of and orbibundle.

Applying Theorem \ref{thm:intro-extensio} to the limit triple $(P,A,\phi)$ and
recalling that the moment map $\mu$ is equivariant we deduce that
$\mu(\phi)$ extends continuously to the normalization of $C$. Then
the vortex equation implies that $F_A$ also extends continuously,
so the limit connection $A$ is indeed meromorphic.

\subsection{}
In case the norms $|d_{A_u}\phi_u|$ are not bounded one can obtain
a meaningful limit object by adding bubbles (rational components) to
the curve $C$. When bubbles are attached away from the nodes
the picture looks exactly like in Gromov--Witten theory (see \cite{M}).
While the first equation in (\ref{eq:the-equations-intro}) is conformally
invariant, the second is not, and as we zoom in the curve (which is
what we do to construct the limit object in the bubbles) the connection
becomes more and more flat. So both $P$ and $A$
are trivial on bubbles away from nodes, and the only nontrivial object
on the bubble is
the section $\phi$, which now can be seen as a holomorphic
map from $S^2$ to the fibre of $P\times_{S^1}X$ over the point in $C$ at
which the bubble is attached.

However, if bubbling off occurs near the nodes, new features appear.
In the simplest case the curve $C$ is replaced by the union of
the normalisation $C'$ and a rational curve
$C_0$, meeting in two nodes. The limit pair $(P,A)$ need not be trivial
on $C_0$: although the connection $A$ will still be flat, it
might have poles (nontrivial holonomy) around the nodes. In this case,
the resulting bubble is not quite a rational curve.

To give an example we use the cylindrical model $\RR\times S^1$
with coordinates $t,\theta$ for the rational curve $S^2$ minus two
points. Consider the real function $H:=-\imag\mu:X\to\RR$. For any
real number $l\neq 0$, take a gradient flow line $\psi:\RR\to X$
at speed $l$ for the function $H$ (so that $\psi$ satisfies the
equation $\psi'=l \nabla H(\psi)).$ Let $\phi:\RR\times S^1\to X$
be the map $\phi(t,\theta):=\psi(t)$, and let $\alpha$ be the
$1$-form $\imag l d\theta$. Then $d_A:=d+\alpha$ gives a covariant
derivative in the trivial bundle $P$ over $\RR\times S^1$ and we
can look at $\phi$ also as a section of the associated bundle with
fibre $X$. By convention, we assume that the complex structure in
$\RR\times S^1$ satisfies $\partial/\partial
t=\imag\partial/\partial\theta$. With these definitions, $A$ is
flat and we have $\ov{\partial}_A\phi=0$. Indeed, the latter
equation takes the following form in terms of the cordinates
$t,\theta$:
\begin{equation}
\label{eq:Floer} \frac{\partial \phi}{\partial t}
=I(\phi)\frac{\partial \phi}{\partial \theta}+l\nabla H(\phi),
\end{equation}
which is clearly satisfied. Note also that equation
(\ref{eq:Floer}) is the equation satisfied by connecting orbits in
Floer's complex for the Hamiltonian $lH$ (see \cite{F}).

Taken as a connection over $S^2$ minus two points, $A$ is
meromorphic, with holonomy around the poles equal to
$\exp(\pm 2\pi\imag l)$. Let $x_\pm$ be the limit of $\psi(t)$ as $t$ goes
to $\pm\infty$. The contribuition of a bubble $C'$ to the energy
of the limit object is given by the squared $L^2$ norm of
$d_A\phi$ on $C'$. In the
example which we have constructed this is
\begin{equation}
\label{eq:energia-intro}
\|d_A\phi\|_{L^2}^2=2\pi l (H(x_+)-H(x_-)),
\end{equation}
so it is finite.

\subsection{}
The factor $H(x_+)-H(x_-)$ in formula (\ref{eq:energia-intro}) is
uniformly bounded by the constant $\sup_X H-\inf_X H$. Hence,
taking $l$ very small, we can obtain bubbles with as small energy
as we wish. This makes a big difference with Gromov--Witten
theory, where the
energy of a rational curve cannot be arbitrarily small. This plays
an important role in the compactness theorem in Gromov--Witten
theory, where it is used to bound the number of bubbles in a
stable map in terms of the energy.

To get a bound on the number of bubbles in our situation (which is
crucial if we want to have a reasonable compactness theorem for
twisted holomorphic maps
of bounded energy), we need to control the geometry of
nontrivial bubbles with little energy. This is done in Theorem
\ref{thm:poca-energia-tururut}. In fact, the example constructed
above gives a hint of what the general situation is:
nontrivial bubbles with little energy can be identified with
gradient lines of $H$, which go upward or downward depending on
the holonomy of the connection (specifically, the holonomy has to
be {\it near but not quite} critical, and the direction of the
gradient line depends on the side of the nearest critical holonomy
in which the holonomy of $A$ is). A particular case of this fact
was already proved by Floer in Theorem 5 of \cite{F} (recall that the equation
$\ov{\partial}_A\phi=0$ for bubbles is equivalent to equation
(\ref{eq:Floer}) for connecting orbits in Floer's complex),
namely, that in which the critical points of $H$ are isolated and
the indices of the limit points as $t\to\pm\infty$ differ by one.
The argument given by Floer relies on index computations and works
only for generic almost complex structure.
Our method of proof is different and is based on an analysis of how
much a bubble with low energy deviates from being a gradient line
(see Theorems \ref{thm:small-energy} and \ref{thm:psi-gradient}; a
similar statement is given in Theorem \ref{thm:small-energy-intro}
below in this introduction). On the other hand, in our result we
require the hamiltonian $H$ to generate an action of $S^1$, and from
this point of view Floer's result is more general than ours.

Using the results in Theorem \ref{thm:poca-energia-tururut} we can
associate to any chain of nontrivial bubbles with little energy a
chain of gradient lines, each one going from one critical point to
another, with as many components as bubbles. Furthermore, since
the bubbles are consecutive all the holonomies have to be the
same, so the chain of gradient lines is monotone (either always
upward or always downward). Finally, since $X$ is compact, the
number of components of a monotone chain has to be bounded.

\subsection{}
Let $\phi'$ be the extension of $\phi$ to the blow up $\wt{C}$
given by Theorem \ref{thm:intro-extensio}. Another important
question is whether, in the absence of blow up, the restrictions
of $\phi'$ to each exceptional divisor $S_y$ and $S_{y'}$ coincide
(up to gauge transformation).  The answer turns out to be no: the
orbits to which $\phi'(S_y)$ and $\phi'(S_{y'})$ belong need not
be the same. Instead, there is a monotone chain of gradient
segments going from one to the other. Here by gradient segment we
mean the image of a map $\xi:T\to X$ satisfying the equation
$\xi'=-\nabla H(\xi)$, where $T\subset\RR$ is any closed interval
(not necessarily infinite). A chain of gradient segments is a
collection of gradient segments $\xi_1(T_1),\dots,\xi_r(T_r)$ such
that for any $j$ the segments $\xi_j(T_j)$ and
$\xi_{j+1}(T_{j+1})$ meet at a unique critical point $f_j\in F$,
and the chain is said to be monotone roughly speaking if $H(f_j)$
decreases with $j$.

This makes another difference with Gromov--Witten theory, where
the analogous question is certainly true. The key point is that
the diameter of the image of a holomorphic map $f:C\to X$, where
$C$ is any long cylinder, can be bounded in terms of the energy
but independently of the length of $C$. This follows from the
exponential decay of the energy density $|df|$ as we go away from
the boundary of the cylinder. To understand why there is such
exponential decay, we write the equation $\ov{\partial}f=0$ as an
evolution equation $f_t(\tau)=L(\tau)f(\tau)$, where $L(\tau)$ is
an elliptic operator on $S^1$ which is skew symmetric up to compact
operators. Then the nonzero spectrum of
$L(\tau)$ stays at distance $\geq\sigma>0$ from $0$, so the
function $|f_t(\tau)|^2_{L^2(S^1)}$ decays as $e^{-\sigma \tau}$
(note that $|df|^2=2|f_t|^2$). This is something very general (see
for example Chapter 3 in \cite{D} for the case of instanton Floer
homology).

If we replace the equation $\ov{\partial}f=0$ by
$\ov{\partial}_A\phi=0$, then the spectra of the operators
$L(\tau)$ are shifted by an amount depending on the holonomy of
$A$. And when the holonomy aproaches a critical value, nonzero
eigenvalues of $L(\tau)$ can approach arbitrarily $0$, so one does
not get exponential decay for $|d_A\phi|^2$. The standard way to
study the geometry of $\phi$ in such a situation is to {\it break}
$\phi$ in two pieces $\psi$ and $\phi_0$, in such a way that
$\psi$ is spanned by the small eigenvalues of the operator and
$\phi_0$ is spanned by the big ones. Then $\phi_0$ does decay
exponentially and so $\psi$ controls the geometry of $\phi$ away
from the boundary of the cylinder. Let us explain in concrete
terms how this works in our context.

Consider for simplicity the case in which the holonomy of $A$ is
nearly trivial but not quite. Suppose that $C=[-N,N]\times S^1$
and that $P\to C$ has been trivialized in such a way that $d_A$ is
approximatedly $d+\lambda\theta$, where
$\lambda\in\imag\RR\setminus\{0\}$ is very small. Suppose that
$\phi:C\to X$ has everywhere very little energy
$|d_A\phi|<\epsilon$. We can then define $\psi:[-N,N]\to X$ by
setting $\psi(t)$ to be the center of mass of $\phi(t,S^1)$, and
$\phi_0:C\to TX$ by the condition
$\phi(t,\theta)=\exp_{\psi(t)}\phi_0(t,\theta)$. The following
theorem, which combines parts of Theorem \ref{thm:small-energy}
and Theorem \ref{thm:psi-gradient}, gives the main properties of
$\phi_0$ and $\psi$ when $\phi$ is a solution of the equations.

\begin{theorem}
\label{thm:small-energy-intro} There are constants $K>0$ and
$\sigma>0$ with the following property. Suppose that $v=f dt\wedge
d\theta$ is a volume form on $C$ such that $f(t)<\eta
e^{-(N-|t|)}$. Assume that the pair $(A,\phi)$ satisfies the
equations $$\ov{\partial}_A\phi=0 \qquad\text{and}\qquad
\iota_{v}F_A+\mu(\phi)=c.$$
 Then we have
$|\phi_0(t,\theta)|\leq K e^{-\sigma(N-|t|)}$ and
\begin{equation}
\label{eq:psi-intro}|\psi'(t)+\imag\lambda\nabla H(\psi(t))|\leq K
e^{-\sigma(N-|t|)}(\epsilon+|\lambda|+\eta)^{1/4}.
\end{equation}
\end{theorem}

The condition on the volume form is natural, since we will
consider cylinders like $C$ as conformal models of curves of the
form $\{xy=\delta\}\subset\CC^2$ near the origin, with the metric
induced by $\CC^2$. The exponential bound on $|\phi_0|$ implies that
away from the boundary of $N$ the function $\phi$ can be
approximated by $\psi$, and equation (\ref{eq:psi-intro}) tells us
that $\psi$ flows near a gradient line at speed $-\imag\lambda$.

\subsection{}
To study the behavior of the solutions of
(\ref{eq:the-equations-intro}) when $C$ degenerates to a nodal
curve, we take cylinders as conformal models, and both $\epsilon$
and $\eta$ go to $0$. Hence, if $\lambda\to 0$ as well,
(\ref{eq:psi-intro}) implies that the map $\psi$ looks more and
more like a gradient line. The following result, which is Theorem
\ref{thm:convergeix-cap-a-gradient}, makes this fact precise. The
actual meaning of convergence of maps to a chain of gradient
segments is given in Section
\ref{ss:convergence-chain-grad-lines}. A good approximation of
this notion is the convergence of subsets in the Hausdorff metric.

\begin{theorem}
\label{thm:convergeix-intro} Suppose that $\{\psi_u:T_u\to X, l_u,
G_u\}$ is a sequence of triples, where $\psi_u$ is a smooth map,
$T_u\subset\RR$ is a finite interval, each $l_u$ is a nonzero real
number and each $G_u>0$ is a real number. Suppose that for some
$\sigma>0$ and any $u$ and $t\in T_u$ we have
\begin{equation}
\label{eq:psi-aprox-intro} |\psi_u'(t)-l_u \nabla
H(\psi_u(t))|\leq G_ue^{-\sigma d(t,\partial T_u)}.
\end{equation}
Suppose also that $G_u\to 0$ and that $l_u\to 0$. Passing to a
subsequence, we can assume that $l_u |T_u|$ converges somewhere in
$\RR\cup\{\pm\infty\}$ (here $|T_u|$ denotes the length of $T_u$).
Then we have the following.
\begin{enumerate}
\item If $\lim l_u|T_u|=0$ then $\lim \diam \psi_u(T_u)=0$. \item
If $\lim l_u|T_u|\neq 0$,  define for big enough $u$ and for every
$t\in S_u$ the rescaled objects $S_u:=l_u T_u$ and
$f_u(t):=\psi_u(t/l_u)$. There is a subsequence of $\{f_u,S_u\}$
which converges to a chain of gradient segments $\TTT$ in $X$.
\end{enumerate}
\end{theorem}

This theorem is almost evident when the critical points of $H$
(that is, the fixed points $F$) are isolated. In this case it
suffices to pay attention to the complementary of a small
neighbourhood of $F$ (where, when $u$ is big enough, $l_u\nabla H$
dominates the error term and hence the vector field $\psi_u'$ is
almost equal to $l_u\nabla H$), and then shrink the neighbourhood
to $F$. The main difficulty appears when the points of $F$ are not
isolated and we want to prove that the limit object is connected:
this essentially amounts to proving that the portion of
$\psi_u(S_u)$ which lies near $F$ has small diameter, which is not
at all obvious (near $F$ the vector field is small, so there is no
hope to control pointwise the error term using $\nabla H$
--- note that we also expect the preimage of a portion of
$f_u(S_u)$ near $F$ to become longer and longer as $u$ goes to
infinity, so the path $f_u$ may wander slowly near $F$ but for a
long time). It is for proving an estimate on the diameter of the
intersection of $f_u(S_u)$ with small neighbourhoods of $F$ set
(see Lemma \ref{lemma:diametre-petit}) that we use an exponential
bound in the error term.

\subsection{}
In view of all the preceding observations, it is clear that the
objects which we should take as limits of solutions of
(\ref{eq:the-equations-intro}) over smooth curves degenerating to
nodal curves have to be of a mixed nature,
combining two dimensional objects (holomorphic sections) with one
dimensional objects (gradient flow lines). A very similar thing
happens in other related moduli problems. In \cite{CT} the second
author and J. Chen construct a compactification of the moduli
space of harmonic mappings from compact surfaces to compact
Riemannian manifolds, and the limiting objects are a combination
of $2$-dimensional harmonic maps and geodesics (their one
dimensional version). Another instance is the approach suggested
by Piunikhin, Salamon and Schwarz \cite{PSS} to proving the
equivalence of Floer and quantum cohomology rings using {\it spiked disks}.

The actual definition of $c$-STHM's (see Section
\ref{ss:def-STHM}) incorporates two additional features. The first
one is that we consider marked points in the curves. Following the
philosophy of Gromov--Witten theory, we treat marked points
and nodes on an equal basis, so in particular we allow
poles of the meromorphic connections on marked points. The other feature is
some information, for each node of the curve, on how to define the
bundle in a smoothening of node (this will be relevant for doing
{\it gluing}, which will be addressed in a future paper). We call
this information gluing data (see Section \ref{ss:gluing-data} for
the precise definition). We also give a notion of equivalence
between $c$-STHM's, we define a topology on the set of isomorphism
classes of $c$-STHM's (see Section \ref{ss:convergence-STHM}), and
define the Yang--Mills--Higgs functional $\YMH_c$ for $c$-STHM's
(which is essentially like (\ref{eq:YMH-intro})). The main theorem
of the paper is then the following (see Theorem
\ref{thm:main-compactness}).

\begin{theorem}
Let $g$ and $n$ be nonnegative integers satisfying $2g+n\geq 3$.
Let $K>0$ be any number, and let $c\in\imag\RR$. Let $\{\CCC_u\}$
be a sequence of $c$-stable twisted holomorphic maps of genus $g$
and with $n$ marked points, satisfying $\YMH_c(\CCC_u)\leq K$ for
each $u$. Then there is a subsequence $\{[\CCC_{u_j}]\}$
converging to the isomorphism class of another $c$-stable twisted
holomorphic map $\CCC$ . Furthermore, we have
$$\lim_{j\to\infty}\YMH_c(\CCC_{u_j})=\YMH_c(\CCC).$$
\end{theorem}

\subsection{Hamiltonian Gromov--Witten invariants}
The main application of the results in this paper is the
construction of Hamiltonian Gromov--Witten invariants coupled to
gravity for compact symplectic manifolds, extending the
construction given in \cite{M}. This will appear in \cite{MT}, and
will be based on the compactness result proved here and the
technique of virtual moduli cycles. The invariants defined in
\cite{M} where obtained by integrating cohomology classes in the
moduli space $\NNN$ of twisted holomorphic maps with a fixed curve
$C$ and representing a given class in $H_2^{S^1}(X)$. In \cite{MT}
we will construct the moduli space $\oN_{g,n}$ of $c$-STHM's of
genus $g$ and $n$ marked points and we will define the Hamiltonian
Gromov--Witten invariants coupled to gravity. The main new feature
compared to \cite{M} (beyond removing the semipositivity
conditions imposed on the manifold and the action) is that the
moduli space fibres over the Deligne--Mumford moduli space
$\oM_{g,n}$. This implies that the resulting invariants satisfy a
gluing axiom similar to the one appearing in Gromov--Witten
theory.

\subsection{Contents}
We now briefly summarize the remaining sections of the paper.
Sections \ref{s:symplectic-fibrations} to \ref{s:orbibundles} are
preparatory. In Section \ref{s:symplectic-fibrations} we recall
the basic objects in the geometry of Hamiltonian almost Kaehler
fibrations (covariant derivatives, $d$-bar operators and minimal
coupling form). In Section \ref{s:meromorphic-connections} we
define meromorphic connections on bundles over punctured Riemann
surfaces, and in Section \ref{s:orbibundles} we recall the
definition of orbibundle over orbisurfaces. These sections are
included in order to fix notations in a way suitable for our
purposes, and there is essentially nothing new in them.

Sections \ref{s:STHM} to \ref{s:main-thm} are devoted to defining
$c$-STHM's and stating and proving the main theorem on
compactness. In Section \ref{s:STHM} we define $c$-STHM's and
isomorphisms between them, and in Section \ref{s:topology-STHM} we
define a topology on the set of isomorphism classes of $c$-STHM's.
Section \ref{s:YMH} is devoted to the Yang--Mills--Higgs
functional, which plays the role of energy for $c$-STHM's. In
Section \ref{s:bounding-bubbles} we prove that a bound on the
energy of a $c$-STHM imposes a bound on the number of irreducible
components of the underlying curve. This plays an important role
in the compactness theorem, which is the main result of the paper
and is stated and proved (modulo technical details) in Section
\ref{s:main-thm}.

All technical details of the proof of the compactness theorem are
given in Sections \ref{s:exp-decay} to
\ref{s:limits-gradient-lines}. In Section \ref{s:exp-decay} we
prove an exponential decay for pairs $(A,\phi)$ satisfying
$\ov{\partial}_A\phi=0$ and having finite energy on punctured
disks. As a Corollary we prove Theorem \ref{thm:intro-extensio} in
this introduction. We also prove that the energy of pairs on long
cylinders decays exponentially provided the holonomy of the
connection stays away from the critical residues. In Section
\ref{s:long-cylinders} we study pairs on long cylinders with
nearly critical holonomy, and we prove that in these pairs the
section stays near a curve in $X$ which satisfies up to
normalization and with some error the equation for gradient lines
(this is essentially Theorem \ref{thm:small-energy-intro}).
Finally, in Section \ref{s:limits-gradient-lines} we prove Theorem
\ref{thm:convergeix-intro}.

\subsection{Notation}

We fix here some notation and conventions which will be used
everywhere in the paper. We denote by $X$ a compact symplectic
manifold with symplectic form $\omega$. The manifold $X$ supports
an effective Hamiltonian action of $S^1$ with moment map
$\mu:X\to(\Lie S^1)^*$. We also take a $S^1$-invariant almost
complex structure $I$ on $X$, compatible with the symplectic form
$\omega$ in the sense that $g:=\omega(\cdot,I\cdot)$ is a
Riemannian metric.

We identify $\Lie S^1\simeq \imag\RR$ with its dual $(\imag\RR)^*$
by using the pairing $\la a,b\ra:=-ab$ for any $a,b\in\imag\RR$.
Thus in the rest of the paper we will assume that $\mu$ takes
values in $\imag\RR$. Let $\XXX$ be the vector field on $X$
generated by the infinitesimal action of $\imag\in\imag\RR$. Thus
the infinitesimal action of any element $\lambda\in\imag\RR$ gives
rise to the vector field $-\imag\lambda\XXX$. We will denote
$H:=-\imag\mu$. Using our convention, we have $H=\la
\mu,\imag\ra$, so that $dH=\la d\mu,\imag\ra=\iota_{\XXX}\omega$
by the definition of the moment map. It follows that the gradient
$\nabla H$ is equal to $I\XXX$.

Define, for any pair of points $x,y\in X$ the pseudodistance
$$\dist_{S^1}(x,y):=\inf_{\theta\in S^1}d(x,\theta\cdot y),$$
where $d$ denotes the distance in $X$ defined by the Riemannian
metric. Also, if $M\subset X$ is any subset, define
$$\diam_{S^1}M:=\sup_{x,y\in M}\dist_{S^1}(x,y).$$
Finally, if $A,B\subset X$ are subsets, we define
$$\dist_{S^1}(A,B):=\inf_{x\in A,\ y\in B}\dist_{S^1}(x,y).$$

Whenever we talk about (marked) curves we will implicitly mean
that they are compact connected complex curves with at most nodal
singularities and eventually with some marked points. These will be
a list of smooth labelled points $x_1,\dots,x_n$, and we will often
denote it by a boldface $\bx$. When we write $x\in\bx$ we will
mean that $x$ is one of the points of the list.

Many of the constants appearing in the estimates will be denoted
by the same symbol $K$, and often the value of $K$ will
change from line to line.

\section{Almost Kaehler geometry of Hamiltonian $S^1$ fibrations}
\label{s:symplectic-fibrations}

Let $C$ be a (nonnecessarily compact) complex curve and let $P$ be
a principal $S^1$ bundle over $C$. Denote by $Y$ the twisted
product $P\times_{S^1}X$ and by $\pi:Y\to C$ the natural
projection. Let $T^{\vert}$ be the vertical tangent bundle $\Ker
d\pi\subset TY$. Since the almost complex structure $I$ is
$S^1$-invariant, it defines a complex structure of each fibre of
$T^{\vert}$, which we still denote by $I$. We also denote by $g$
the Euclidean form on $T^{\vert}$ defined by the $S^1$-invariant
metric $g$ on $X$. Let $I_{C}$ be the complex structure on $C$ (so
$I_C$ is an endomorphism of the real tangent bundle of $C$), and
let $g_C$ be a conformal metric on $C$. Any connection $A$ on $P$
induces a splitting $TY\simeq \pi^*TC\oplus T^{\vert}$, which we
use to define $I(A)$ (resp. $g(A)$) as the sum of $\pi^*I_C$ and
$I$ (resp., $\pi^*g_C$ and $g$).

Suppose that $\phi$ is a section of $Y$. We define the covariant
derivative $d_A\phi$ of $\phi$ with respect to the connection $A$
to be the composition of $d\phi$ with the projection
$\pi_A^v:TY\to T^{\vert}$ induced by $A$. Hence, $d_A\phi$ is a
one form on $C$ with values in the pullback $\phi^*T^{\vert}$.

\begin{lemma}
\label{lemma:cov-der-diam-S1} Let $S$ be the interval $[0,1]$, let
$\gamma:S\to C$ be a smooth map, and let $\phi$ be a section of
$Y$. Then $\dist_{S^1}(\phi\gamma(1), \phi\gamma(0))$ is at most
the integral over $S$ of
$|\la d_A\phi(\gamma),\gamma'\ra|.$
\end{lemma}
\begin{pf}
Take a trivialization $\gamma^*Y\simeq [0,1]\times X$, and
consider a gauge transformation $g:[0,1]\to S^1$ which sends
$\gamma^*A$ to the trivial connection. Then we have
$\dist_{S^1}(\phi\gamma(1),\phi\gamma(0))=
\dist_{S^1}(g\phi\gamma(1),g\phi\gamma(0))$, and the latter can
easily be estimated in terms of the integral of
$|\la d_A\phi(\gamma),\gamma'\ra|$ (note
that by covariance we have $\la d_A\phi(\gamma),\gamma'\ra=\la
d_{g^*\gamma^*A}(g\phi),\partial/\partial t\ra$).
\end{pf}

We can split the space of forms $\Omega^1(C,\phi^*T^{\vert})$ as
the sum of the space of holomorphic forms
$\Omega^{1,0}(C,\phi^*T^{\vert})$ plus the space of
antiholomorphic forms $\Omega^{0,1}(C,\phi^*T^{\vert})$. Let
$\ov{\partial}_{I,A}\phi$ be the projection of $d_A\phi$ to
$\Omega^{0,1}(C,\phi^*T^{\vert})$.
In concrete terms,
\begin{equation}
\label{eq:dbar} \ov{\partial}_{I,A}\phi:=\frac{1}{2}(d_A\phi+I
\circ d_A \circ I_C).
\end{equation}
(When the complex structure on $X$ will be clear from the context,
we will simply write $\ov{\partial}_A\phi$.) Let us denote by
$\Phi:C\to Y$ the map defined by the section $\phi$. It is
straightforward to check that
\begin{equation}
\label{eq:holomorf}
\ov{\partial}_{I,A}\phi=0\qquad\Longleftrightarrow\qquad
\ov{\partial}_{I(A)}\Phi=0.
\end{equation}
On the other hand, for any section $\Phi$ of $Y$ we have
\begin{equation}
\label{eq:norma-dPhi}
|d\Phi|^2_{g(A)}=|d\Id_C|^2+|d_A\phi|^2=1+|d_A\phi|^2.
\end{equation}

Given a connection $A$ on $P$, there is a canonical way to pick a
closed $2$-form $\omega(A)$ on $Y$ which restricts to $\omega$ on
each fiber. This is the {\bf minimal coupling form} (see
\cite{GLS}). To give a local description of it we can assume that
there is a trivialisation $P\simeq C\times S^1$. Let
$\alpha\in\Omega^1(C,\imag\RR)$ be the $1$-form corresponding to
$A$ with respect to this trivialisation. Let also $\pi_X:Y\to X$
the projection induced by the trivialisation. Then we have
\begin{equation}
\label{eq:def-omega-A}
\omega(A)=\pi^*_X\omega-d(\pi^*\alpha\wedge \pi_X^*\mu)
=\pi^*_X\omega-\pi^*\alpha\wedge\pi^*_X\iota_{\XXX}\omega
-\pi^*F_A\wedge\pi^*_X\mu,
\end{equation}
where $F_A$ is the curvature of $A$.
If $A$ is flat, then $\omega(A)$ coincides with
the $2$ form $g(A)(\cdot,-I(A)\cdot)$ when restricted
to vertical tangent vectors.

\begin{remark}
\label{rmk:coh-class-omega} The cohomology class represented by
$\omega(A)$ is independent of $A$, and can in fact be identified
to the pullback of the class $[\omega-\mu t]$ in equivariant
cohomology (this denotes the class represented by the element
$\omega-\mu t$ of the Cartan--Weil complex
$\RR[t]\otimes\Omega(X)^{S^1}$, see for example \cite{GS}). More
precisely, if $c:C\to BS^1$ is the classifying map for $P$, then
we can identify $Y\simeq c^* X_{S^1}$ and $[\omega(A)]$ is equal
to $c^*[\omega-\mu t]$.
\end{remark}

\section{Meromorphic connections on marked nodal surfaces}
\label{s:meromorphic-connections}

\subsection{Meromorphic connections on marked smooth curves}

Let $(C,\bx)$ be a marked smooth complex curve, and let $P$ be a
principal $S^1$ bundle over $C\setminus\bx$. For any $x\in\bx$ we
denote by $T(P,x)$ the set of trivialisations up to homotopy of
the restriction of $P$ to a small loop $\gamma$ around $x$ (since
for two different loops $\gamma$ and $\gamma'$ the sets of
trivialisations of the corresponding restrictions can be
canonically identified, $T(P,x)$ is independent of the chosen
loop). Orient $\gamma$ counterclockwise with respect to the
natural orientation of $C$ as a Riemannn surface, $T(P,x)$ gets a
natural structure of $\ZZ$-torsor\footnote{Recall that if $\Gamma$
is a group, a $\Gamma$-torsor is a set $T$ with a free left action
of $\Gamma$.}, given by the action of gauge transformations
defined over $\gamma$: if $t:P|_{\gamma}\to \gamma\times S^1$ is a
trivialisation and $\tau:=[t]$ denotes its class in $T(P,x)$, then
for any $n\in\ZZ$ we define $n\cdot \tau$ to be the class of the
trivialisation $t\circ g$, where $g:P|_{\gamma}\to P|_{\gamma}$ is
the gauge transformation given by any map $\gamma\to S^1$ of index
$n$.

Any $\tau\in T(P,x)$ defines an extension of $P$ to $x$:
if $D\subset C$ is a disk containing $x$ and having $\gamma$
as boundary, then we glue $P$ to the trivial bundle over $D$
using any extension of a representative of $\tau$ to $D\setminus\{x\}$
as patching function. This extension is well defined up to homotopy.
In general, if we pick a collection of local trivialisations
$$\tau\in T(P,\bx):=\prod_{x\in\bx} T(P,x)$$
then we get an extension of $P$ to $C$, which we denote by $P^{\tau}$.
We define the {\bf degree} of the bundle $P$ over $(C,\bx)$ to be the map
$$\deg P: T(P,\bx)\to\ZZ$$
which sends any $\tau\in T(P,\bx)$ to $\deg P(\tau):=\deg P^{\tau}$.

We say that a connection $A$ on $P$ is {\bf meromorphic} with poles
in $\bx$ if its
curvature $F_A$ extends to the whole $C$ as a continuous $2$-form.

Suppose that $A$ is meromorphic. Chose any metric on $C$, and fix
some $x\in\bx$. Let $\gamma_{\epsilon}$ denote the boundary of a geodesic
disk centered on $x$ and of radius $\epsilon$, oriented
counterclockwise. Denote by $\Hol(A,\gamma_{\epsilon})\in S^1$
the holonomy of $A$ around $\gamma_{\epsilon}$.
We claim that the limit
\begin{equation}
\label{eq:def-Hol}
\Hol(A,x):=\lim_{\epsilon\to 0}\Hol(A,\gamma_{\epsilon})
\end{equation}
exists. Indeed, if $\epsilon>\epsilon'>0$ then the quotient
$\Hol(A,\gamma_{\epsilon})/\Hol(A,\gamma_{\epsilon'})$ is equal,
by Stokes' theorem, to the exponential of
the integral of the curvature $F_A$ over the region
enclosed between $\gamma_{\epsilon}$ and $\gamma_{\epsilon'}$.
But since $F_A$ extends continuously to $x$, this integral
tends to $0$ as $\epsilon\to 0$, so the quotient of the holonomies
tends to $1$. The same argument proves that the definition of
$\Hol(A,x)$ is independent of the chosen metric.

\begin{lemma}
\label{lemma:bones-coordenades}
Let $(r,\theta)$ be polar coordinates defined on a small disk
$D$ centered on $x\in\bx$. Any
$\tau\in T(P,x)$ has a representative which
extends to give a trivialisation of $P$ on $D\setminus\{x\}$
with respect to which the covariant derivative associated to
$A$ can be written as
\begin{equation}
\label{eq:bona-trivialitzacio}
d_A=d+\alpha+\lambda d\theta,
\end{equation}
where
$\alpha\in\Gamma(T^*D\otimes\imag\RR)$ is a $1$-form
of type $C^1$, and $\lambda\in\imag\RR$ satisfes
$\Hol(A,x)=e^{2\pi\lambda}.$
Furthermore, any two choices of trivialisations associated to
the same $\tau\in T(P,x)$ are related by a gauge
transformation which extends to $D$.
\end{lemma}
\begin{pf}
Take a trivialisation of $P$ on $D\setminus\{x\}$ extending
any representative of
$\tau$. With respect to this trivialisation we can write
$d_A=d+\alpha_0$, so that $F_A=d\alpha_0$. Let
$$\lambda:=\lim_{\epsilon\to 0}\frac{1}{2\pi}\int_{\theta\in S^1}
\alpha_0(\epsilon,\theta).$$
That this limit exists follows as in the definition of
$\Hol(A,x)$ from Stokes' theorem,
writing the difference of the integrals on the right hand side
for two choices of $\epsilon$ as the integral of $d\alpha$ on
the corresponding annulus, and then using the fact that
$d\alpha_0$ extends continuously to $D$.
By Poincar\'e's lemma, there exists $\alpha\in\Omega^1(D,\imag\RR)$
of type $C^1$ such that $d\alpha_0=d\alpha$. Hence
$\beta:=\alpha_0-\alpha-\lambda d\theta$ is closed. But it is also
exact, since for any $\epsilon>0$ we have
$$\int_{\theta\in S^1}\beta(\epsilon,\theta)
= \lim_{\epsilon\to 0}
\int_{\theta\in S^1}\beta(\epsilon,\theta)
= \lim_{\epsilon\to 0} \int_{\theta\in S^1}\alpha_0-\lambda d\theta
= 0 $$
(in the first equation we use that $\beta$ is closed, in the second
one that $\alpha_0$ extends to $D$ and in the third one the definition
of $\lambda$).
Consequently, we can write $\beta=dg$ for some $g:D\setminus\{x\}\to S^1$.
Now, the gauge transformation $G:=\exp g$ transforms
$d+\alpha_0$ into $d+\alpha+\lambda d\theta$.
If another gauge transformation $G'$ with winding
number $0$ transform $d+\alpha_0$ into $d+\alpha'+\lambda d\theta$
and $\alpha'$ extends to $D$, then we can write $G'=dg'$ (because
the winding number of $G'$ is $0$) and
$dg-dg'=\alpha-\alpha'.$
The right hand side is a closed $1$-form which extends to $x_j$,
hence $g-g'$ extends to $D$.
\end{pf}

\begin{corollary}
\label{cor:trivial-hol-extension} Let $P$ be a bundle over
$C\setminus\bx$ and let $A$ be a meromorphic connection on $P$.
Let $\bx'\subset\bx$ be the marked points around which the holonomy
of $A$ is nontrivial. Then there is a bundle $P'\to C\setminus\bx'$
with a smooth
connection $A'$ such that $\iota^*(P',A')\simeq (P,A)$, where
$\iota:C\setminus\bx\to C\setminus\bx'$ denotes the inclusion.
\end{corollary}

\begin{remark}
\label{rmk:radial-gauge}
In fact, it is even possible to trivialize $P$ around $x$ in such a
way that $d_A$ is in {\bf radial gauge}, that is,
it takes the form $d+\alpha+\lambda d\theta$ and
$\alpha=\alpha_{\theta} d\theta$ is of type $C^1$ (in particular,
$\alpha_{\theta}$ vanishes at $0$). Indeed, given a covariant derivative
$d+\alpha+\lambda d\theta$ as in Lemma \ref{lemma:bones-coordenades}
on the trivial bundle over the punctured
disk, we can define $g:D\to S^1$ to be the map
whose value at $z\in D$ gives
the parallel transport with respect to $d+\alpha$
from the fibre over $z$ to that over $0$. Then
$d+\alpha-g^{-1}dg$ is in radial gauge. Restricting to the punctured disk,
the gauge transformation $g$ sends $d+\alpha+\lambda d\theta$ to a
connection in radial gauge.
\end{remark}

We call $\lambda=:\Res(A,x,\tau):=\Res(A,x,\tau_x)$
the {\bf residue} of $A$ on $x$ with respect to $\tau$.
It is straightforward to check that for any $n\in\ZZ$ we have
\begin{equation}
\label{eq:variacio-residu}
\Res(A,x,n\cdot \tau)=\imag n +\Res(A,x,\tau).
\end{equation}

\begin{remark}
\label{rmk:curvatura-acotada}
One can check easily that if $A$ is a connection whose curvature
$F_A$ is uniformly bounded (but not necessarily extends continuously
to $C$) then Lemma \ref{lemma:bones-coordenades} is still true,
with the difference that the $1$-form $\alpha$ is in general only
continuous. Consequently, one can also define
in this situation the residue $\Res(A,x,\tau)$ and formula
(\ref{eq:variacio-residu}) still holds.
\end{remark}

\subsection{Meromorphic connections on marked nodal curves}
\label{ss:mero-connections}

Now suppose that $(C,\bx)$ is a marked curve with nodal
singularities. Denote by $\bz$ the set of nodes of $C$. Let
$\pi:C'\to C$ be the normalisation map. We denote the preimage by
$\pi$ of the marked points $\bx$ with the same symbol $\bx$.

If $P$ is a bundle over $C\setminus(\bx\cup\bz)$, by a {\bf
meromorphic connection} on $P$ we mean a meromorphic connection on
$\pi^*P$ such that for every node $z\in\bz$, denoting by
$y,y'$ the preimages of $z$ in the normalization, we have
$$\Hol(A,y)=\Hol(A,y')^{-1}.$$
Define\footnote{If $T$ and $T'$ are
two $\Gamma$-torsors and $\Gamma$ is abelian, we denote
$T\times_{\Gamma}T':=T\times T'/\sim$, where the equivalence
relation $\sim$ identifies, for any $(t,t')\in T\times T'$ and
$g\in\Gamma$, $(g\cdot t,t')\sim (t,g\cdot t').$}
for any such node $T(P,z):=T(\pi^*P,y)\times_{\ZZ}T(\pi^*P,y')$.
A {\bf marked $S^1$ principal bundle} over $(C,\bx)$ is a pair
$(P,\sigma)$, where $P$ is a principal $S^1$ bundle over
$C\setminus(\bx\cup\bz)$ and
$$\sigma\in T(P,\bz):=\prod_{z\in\bz}T(P,z).$$  Let $C^s$ be a
smoothening of $C$ (hence, $C^s$ is isomorphic to $C$ away from a
collection of neighbourhoods of the nodes of $C$, and near the
nodes $C^s$ takes the form $\{xy=\epsilon\}\subset\CC^2$ instead
of $\{xy=0\}\subset\CC^2$). The elements in $\sigma$ give a unique
way (up to homotopy) of extending $P$ to a bundle $P^s$ over
$C^s\setminus\bx$ (here we identify $\bx$ with points in $C^s$).
Then we define the {\bf degree} of $(P,\sigma)$ to be the map
$$\deg P^{\sigma}:=T(P,\bx)\to\ZZ$$
given by $\deg P^{\sigma}:=\deg P^s$.

A {\bf meromorphic connection} on $(P,\sigma)$ is a meromorphic
connection $A$ on $\pi^*P$ with poles in $\bx\cup\by$ which
satisfies the following compatibility condition on each node
$z$: if $y,y'$ are the preimages of $z$ in $C'$ and
we take a representative $(\xi,\xi')\in T(\pi^*P,y)\times T(\pi^*P,y')$
of $\sigma_z$ then
\begin{equation}
\label{eq:residues-match}
\Res(A,y,\xi)+\Res(A,y',\xi')=0.
\end{equation}
In other words, the sum of the residues in both preimages of the
node has to vanish. That this condition is well defined follows
from (\ref{eq:variacio-residu}).

It is straightforward to check that given a connection $A$ on $P$
there is a unique choice of $\sigma$ with respect to which $A$ is
meromorphic (that is, formula (\ref{eq:residues-match}) holds).
Using such $\sigma$, we define
$$\deg P^A:=\deg P^{\sigma}.$$

\subsection{Chern--Weil formula}
The following lemma gives a generalisation of the simplest
Chern--Weil formula to the case of meromorphic connections
on nodal curves.

\begin{lemma}
Denote by $F_A$ the curvature of a meromorphic connection $A$ on
$P$. For any $\tau\in T(P,\bx)$ we have
\begin{equation}
\deg P^A(\tau)=\frac{\imag}{2\pi}\int_{C}F_A+\imag\sum_{x\in\bx}
\Res(A,x,\tau). \label{eq:ChernWeil}
\end{equation}
\label{lemma:ChernWeil}
\end{lemma}
\begin{pf}
It is enough to consider the case of smooth $C$ (in the nodal
case, pulling back to the normalisation we are led to the smooth
case). Hence we assume that $C$ is smooth and that $A$ is a
meromorphic connection on a bundle $P$ over $C\setminus\bx$. Let
$\bx=(x_1,\dots,x_n)$. Pick some $\tau\in T(P,\bx)$ restriction of
$P$ to a small loop around each marked point and define
$\lambda_j:=\Res(A,x_j,\tau)$. Take disjoint disks $D_1,\dots,D_n$
centered at the marked points, small enough so that we can apply
Lemma \ref{lemma:bones-coordenades} to get trivialisations of $P$
near each $x_j$ for which the covariant derivative of $A$ takes
the form $d+\alpha_j+\lambda_j d\theta$ (here $(r,\theta)$ are
polar coordinates in $D_j$). For any $j$ chose a smooth 1-form
$\eta_j$ on $D_j$ which coincides with $d\theta$ away from a small
disk $\{r<\epsilon\}\subset D_j$. Now, modify $A$ near each $x_j$
by replacing $d+\alpha_j+\lambda_j d\theta$ by
$d+\alpha_j+\lambda_j\eta_j$. In this way we get a smooth
connection $A'$ on the bundle $P^{\tau}$ over $C$. Then using
standard Chern--Weil and Stokes we compute
\begin{align*}
\deg P(\tau)&=\deg P^{\tau}=\frac{\imag}{2\pi}\int_{C} F_{A'}
 = \frac{\imag}{2\pi}
\int_{C}F_A
+ \frac{\imag}{2\pi} \sum\lambda_j\int_{D_j}d\eta_j \\
& = \frac{\imag}{2\pi}\int_{C}F_A
+ \frac{\imag}{2\pi} \sum\lambda_j\int_{\partial D_j}\eta_j
 = \frac{\imag}{2\pi}\int_{C}F_A
+ \imag\sum\lambda_j.
\end{align*}
\end{pf}

\subsection{Gluing data}
\label{ss:gluing-data}

Suppose that $(C,\bx)$ is a smooth curve, $P$ is a principal
$S^1$ bundle over $C\setminus\bx$ and $A$ is a meromorphic
connection on $P$. For every $x\in\bx$, let $S_x$ be the quotient
of $T_xC^*/\RR$, where $T_xC^*\subset T_xC$ denotes the nonzero
elements and $\RR$ acts on $T_xC^*$ by multiplication. The
conformal structure on $C$ induces an orientation and a metric on
$S_x$ up to constant scalar. Imposing the volume of $C$ to be
$2\pi$, we get a well defined metric on $S_x$.
A useful point of view is to look at $S_x$ as the exceptional
divisor of the real blowup of $C$ at $x$.

The pair $(P,A)$ induces a {\bf limiting pair} $(P_x,A_x)$, where
$P_x$ is a bundle over $S_x$ and $A_x$ is a connection on $P_x$ in
the following way. Pick any conformal metric on $C$ and take a
small $\epsilon>0$. Let $D_{\epsilon}:=\exp_x\{v\mid v\in T_xC^*,\
|v|<\epsilon\}$. Let $D_{\epsilon}/\RR$ denote the quotient by the
equivalence relation which identifies $y,z$ if and only if
$y=\exp_x u$, $z=\exp_x v$ and $u \in\RR v$. We can lift this
equivalence to the restriction of $P$ on $D_{\epsilon}$ by using
parallel transport (with respect to $A$) along lines of the form
$\{\exp_x\lambda u\mid \lambda\in\RR\}$, and we denote the
quotient by $P_x:=P|_{D_{\epsilon}}/\RR$. We look at $P_x$ as a
bundle over $S_x$ using the obvious identification $S_x\simeq
D_{\epsilon}/\RR$. To define the limiting connection we proceed as
follows: for any $\delta>0$ smaller than $\epsilon$, let
$S_{\delta}$ be the exponential of the circle of radius $\epsilon$
in $T_xC$. The composition of inclusion and quotient:
$S_{\delta}\hookrightarrow D_{\epsilon}\to D_{\epsilon}/\RR\simeq
S_x$ is an isomorphism, and we denote by $g_{\delta}$ its inverse.
Using parallel transport $g_{\delta}$ lifts to an isomorphism
$P_x\simeq g_{\delta}^*P|_{S_{\delta}}$. One checks, using the
fact that $A$ is meromorphic, that the pullback connections
$g_{\epsilon}^*A$ converge to a limit connection, which we denote
by $A_x$. (For example, use a trivialisation around $x$
which puts $d_A$ in radial gauge and as $d_A=d+\alpha+\lambda d\theta$,
where $\alpha$ extends to $x$ --- see Remark \ref{rmk:radial-gauge}.)
The resulting pair $(P_x,A_x)$ is independent of the
chosen metric on $C$. Of course, the holonomy of $A_x$ around $S_x$ coincides
with the holonomy $\Hol(A,x)$ of $A$ around $x$ as defined by
(\ref{eq:def-Hol}).

In other words, the pullback of $(P,A)$ to the real blowup of
$C$ at $x$ extends to the exceptional divisor $S_x$, and the restriction
to $S_x$ of the extension is isomorphic to $(P_x,A_x)$.

Now suppose that $(C,\bx)$ is a nodal curve. Let $z\in C$ be a
node and let $y,y'$ be its preimages in the normalisation $C'$
of $C$. Define the set of {\bf gluing angles} at $z$ to be
$$\Gamma_z:=(T_{y}C'\otimes T_{y'}C')^*/\RR.$$
Recall that the nonzero elements of $T_{y}C'\otimes T_{y'}C'$
specify deformations of $C$ which smoothen the singularity at $z$
(see Section \ref{ss:convergencia-corbes}). On the other hand, the
set $\Gamma_z$ can be identified with the set of isometries
$\gamma:S_{y}\to S_{y'}$ which reverse the orientations.

Let $P$ be a principal $S^1$ bundle over $C\setminus
(\bx\cup\bz)$, where $\bz$ denotes the set of nodes. We define the
set of {\bf gluing data for $P$ at $z$} to be the set
$\Gamma(P,z)$ of pairs $(\gamma,\rho)$, where $\gamma\in\Gamma_z$
is a gluing angle at $z$ and $\rho:\gamma^*P_{y'}\to P_{y}$ is
an isomorphism of bundles satisfying $\rho^*A_{y'}=A_{y}$.
Clearly, since $\gamma$ reverses the orientation the condition for
$\rho$ to exist is that the holonomy of $A_{y'}$ is the inverse
of that of $A_{y'}$. Note that the projection
$\pi_z:\Gamma(P,z)\to\Gamma_z$ has a natural structure of
principal $S^1$ bundle. Finally, we define the set of {\bf gluing
data for $P$} to be the product
$$\Gamma(P):=\prod_{z\in\bz} \Gamma(P,z)$$ of sets of gluing data at each
of the nodes of $C$. The group of gauge transformations of $P$
acts in an obvious way on $\Gamma(P)$ preserving the gluing angles
(in other words, for every node $z$ the induced action on
$\Gamma(P,z)$ preserves the map $\pi_z$).

Any diagram of the form
$$\xymatrix{P\ar[r]^{g}\ar[d] & P'\ar[d] \\
C\ar[r]^{f} & C',}$$ where $g$ is an isomorphism of principal
$S^1$ bundles and  $f$ is a biholomorphism satisfying
$f(\bx)=\bx'$
induces a {\bf morphism of gluing data} $g^*:\Gamma(P')\to\Gamma(P)$
(note that the vertical arrows in the diagrams are not
surjections, since the bundles are defined over the set of smooth
points of $C$ and $C'$ which are not marked).

A choice of gluing data for $P$ should be understood as a device
to specify deformations of $P$ lifting certain smoothings of $C$.
More precisely, suppose that $G\in\Gamma(P)$ consists of  gluing
data $(\gamma_1,\rho_1),\dots,(\gamma_k,\rho_k)$. The set of
gluing angles $(\rho_1,\dots,\rho_k)$ give a real subspace $T$ in
the tangent space of $\oM_{g,n}$ at the point represented by
$(C,\bx)$ (suppose for simplicity that we only smoothen the
singularities and do not modify the complex structures of the
irreducible components of $C$). Then, for each deformation $C_t$
parametrized by some small segment $t\in [0,\epsilon)$ tangent to
$T$, $G$ gives a way to extend $(P,A)$ along $C_t$, in a {\it
infinitesimally unique} way.

This may seem a little unsatisfactory, since it would be
preferable to specify some data which tells how to deform $(P,A)$
along any deformation of $(C,\bx)$ in $\oM_{g,n}$. Such a thing is
certainly possible if $A$ has trivial holonomy around the
preimages of $z$. In this case, the bundle $\pi^*P$ (recall that
$\pi$ is the normalisation map) extends over $y$ and $y'$
(essentially by Lemma \ref{lemma:bones-coordenades}) to give a
bundle $P'$, and picking an identification of the fibres
$P'_{y}$ and $P'_{y'}$ (which should not be confused with
$P_{y}$ and $P_{y'}$) gives a way to extend $P$ to any
smoothing of $C$ at $z$. Furthermore, the set of possible
identifications between $P'_{y}$ and $P'_{y'}$ is a torsor over
$S^1$.

In view of this we could wonder whether it is possible to chose,
when $A$ has nontrivial holonomy $H$ around $y$ and $y'$, a
family parametrized by $S^1$ of sections of the bundle
$\Gamma(P,z)\to \Gamma_z$, varying continuously with $H$ (the
point is that such a section would give a way to deform $P$ for
each direction of smoothings of $C$). But this is unfortunately
not possible, as we now explain. Denote by $\HHH$, $\HHH'$ the moduli
spaces of connections on the trivial $S^1$ bundle over $S_{y}$, $S_{y'}$.
We then have $\HHH\simeq S^1\simeq\HHH'$ canonically, the
isomorphism being given by the holonomy. There are universal
Poincar\'e bundles $\PPP\to S_{y}\times\HHH$
and $\PPP'\to S_{y'}\times\HHH'$
(with a
universal connection $A$ whose restriction to the fibre
$S_{y}\times\{H\}$ has holonomy $H\in S^1$, and similarly there
is a connection $A'$ on $\PPP'$).
Consider the bundle
$\II:=\Isom(\PPP,\PPP')$ over $\Gamma_z\times S^1$ whose fibre
on $(\gamma,H)$ is the set of isomorphisms between
$\PPP|_{S_{y}\times\{H\}}$ and
$\gamma_I^*\PPP'|_{S_{y'}\times\{H^{-1}\}}$, where
$\gamma_I(\alpha,H):=(\gamma(\alpha),H^{-1})$, which preserve the
connections. What we are looking for is a family of sections
$\Sigma_H\subset\Gamma(\II|_{\Gamma_z\times\{H\}})$ for every $H$
and depending continuously on $H$, at least homeomorphic to $S^1$.
If such a thing existed, it would form a fibration $\Sigma$ over
$S^1$ with connected fibres, and hence it would admit sections.
Any section of $\Sigma$ would induce a section of $\II$, which is
impossible, because the bundle $\II$ is nontrivial (to see that
$\II$ is nontrivial observe that the connections $A$ and $A'$ induce a
connection $\AA$ on $\II$ whose holonomy around
$\Gamma_z\times\{H\}$ is $H$; this forces the curvature of $\AA$
to have integral over $\Gamma_z\times S^1$ equal to $-2\pi\imag$).
Of course, this is a manifestation of the impossibility of
extending the universal Jacobian to the whole Deligne--Mumford
moduli space as a fibration of smooth orbifolds.
See \cite{Fr} for a more general discussion.

\section{Orbibundles over orbisurfaces}
\label{s:orbibundles}

We will call an {\bf orbisurface} a smooth complex curve $C$
together with a list of points $\bx=x_1,\dots,x_n\in C$ and
corresponding positive integers $\bq=q_1,\dots,q_n$. Take an
orbisurface $(C,\bx,\bq)$. We now define (in a rather {\it ad hoc}
way) what is an orbibundle over $(C,\bx,\bq)$. (Of course,
orbibundles are well known objects and can be defined in much
greater generality than we do; we just recall the definition in
this particular case to fix notations.) Pick a neighborhood $U_j$
of each point $x_j$ and a holomorphic covering map
$\rho_j:\oU_j\to U_j$ of degree $q_j$, such that $\oU_j$ is
biholomorphic to the unit disk and $\rho_j$ has maximal
ramification at the unique point $\ox_j:=\rho_j^{-1}(x_j)$ and is
unramified everywhere else. Let $\Gamma_j:=\ZZ/q_j\ZZ$. The group
$\Gamma_j$ acts on $\oU_j$ leaving the map $\rho_j$ invariant. We
will assume that the sets $U_1,\dots,U_n$ are disjoint. A {\bf
principal $S^1$ orbibundle} over $(C,\bx,\bq)$ is a tuple
$$\bP=(P,P_1,\dots,P_n,\psi_1,\dots,\psi_n),$$ where $P\to
C\setminus\bx$ is a principal $S^1$ bundle and, for each $j$,
$P_j\to\oU_j$ is a $\Gamma_j$-equivariant principal $S^1$ bundle
and $\psi_j:\rho_j^*P|_{U_j\setminus\{x_j\}}
\stackrel{\sim}{\longrightarrow} P_j|_{\oU_j\setminus\{\ox_j\}}$
is an isomorphism. We leave to the reader the definition of
isomorphism of principal orbibundles.

To define the {\bf degree} of the orbibundle $\bP$, pick a
collection of local trivialisations $\tau=(\tau_1,\dots,\tau_n)\in
T(P,\bx)$. Each of them gives rise (pulling back through $\rho_j$
and using $\psi_j$) to a trivialisation $\otau_j$ of $P_j$
restricted to $\partial\oU_j$. Such trivalisation allows to define
a quotient of $P_j$ over $\oU_j/\partial\oU_j\simeq S^2$, and we
define $\deg(P_j,\otau_j)$ to be the degree of this quotient. Then
$$\deg\bP:=\deg P^\tau-\sum \frac{1}{q_j}\deg(P_j,\otau_j).$$
One checks, using (\ref{eq:variacio-residu}), that this expression
does not depend on the choice of $\tau$.

\subsection{Relation to meromorphic connections}
\label{ss:mer-conn-orbibundles} Suppose that $(C,\bx)$ is a marked
smooth curve and that $P\to C\setminus\bx$ is a principal $S^1$
bundle. Let $A$ be a meromorphic connection on $P$, all of whose
residues are of the form $\imag l$ for $l\in\QQ$. Then we obtain
in a natural way a structure of orbisurface on $C$ and an
orbibundle $\bP^A$ on it which restricts to $P$ on
$C\setminus\bx$, as follows. Pick a collection of trivialisations
$\tau\in T(P,\bx)$ and let the residue $\Res(A,x_j,\tau_j)$ be
$\imag p_j/q_j$, where $p_j$ and $q_j\geq 1$ are relatively prime
integers. The resulting orbisurface is then $\bC^A:=(C,\bx,\bq)$,
where $\bq=(q_1,\dots,q_n)$, and we define $P_j\to\oU_j$ to be the
trivial bundle with the trivial action of $\Gamma_j$. So it
remains to define the isomorphisms $\psi_j$ (which should be a
trivialisation of $\rho_j^*P$, since $P_j$ is the trivial bundle).
The residue of the connection $\rho_j^*A$ on the bundle
$\rho_j^*P$ at the point $\ox_j$ with respect to $\rho_j^*\tau_j$
is $p_j\in\ZZ$. Then we take $\psi_j$ to be the trivialisation of
$\rho_j^*P$ with respect to which the covariant derivative
$d_{\rho_j^*A}$ is equal to $d+\alpha$, where $\alpha$ extends
continuously to $\Omega^1(\oU_j,\imag\RR)$ (such trivialisation
exists by Lemma \ref{lemma:bones-coordenades}).

It is now a consequence of the Chern--Weil formula in
Lemma \ref{lemma:ChernWeil} that
$$\deg\bP^A=\frac{\imag}{2\pi}\int_{C\setminus\bx} F_A.$$

\subsection{Associated bundles}

Of course, one can consider in general orbibundles which are not
necessarily principal $S^1$ orbibundles, but which are, more
generally, {\bf locally trivial orbibundles}. Their definition is
completely analogous to that of $S^1$ principal orbibundle,
substituting in each case the words ``principal $S^1$ bundle'' by
``locally trivial bundle''. The notion of {\bf section} of an
orbibundle is also almost evident; so, for example, if
$\bY:=(Y,\{Y_j\},\{\psi^Y_j\})$ denotes a locally trivial
orbibundle over an orbisurface $(C,\bx,\bq)$, then a section of
$\bY$ is a collection of sections
$\Phi=(\phi,\phi_1,\dots,\phi_n)$, where
$\phi\in\Gamma(C\setminus\bx,Y)$ and
$\{\phi_j\in\Gamma(\oU_j,Y_j)\}$, satisfying the compatibility
condition $\psi_j\circ\rho_j^*\phi=\phi_j$.

This is the main example which we will encounter. Given a
principal $S^1$ orbibundle $\bP=(P,\{P_j\},\{\psi_j\})$ over
$(C,\bx,\bq)$, we can define a locally trivial orbibundle
$\bY:=\bP\times_{S^1}X$ on $(C,\bx,\bq)$ by setting
$\bY:=(Y,\{Y_j\},\{\psi^Y_j\})$, where $Y=P\times_{S^1}X$,
$Y_j=P_j\times_{S^1}X$ and $\psi^Y_j$ denotes the isomorphism
induced by $\psi_j$.

Recall that the de Rham complex of differential forms can be
defined for orbifolds, and that its cohomology is isomorphic to
the singular cohomology of the orbifold with real coefficients.
This means in particular that we can use Chern--Weil theory to
obtain representatives of Chern classes of orbibundles (allowing,
in particular, to compute their degree). Another consequence is
that one can define also the cohomology class $[\omega(A)]$
exactly as in Section \ref{s:symplectic-fibrations}.

\section{Stable twisted holomorphic maps}
\label{s:STHM}

\subsection{Critical residues}
Let $F\subset X$ be the fixed point set of the action of $S^1$.
For each connected component $F'\subset F$, there is an action of
$S^1$ on the corresponding normal bundle $N\to F'$. Then $N$
splits as a direct sum $N=\bigoplus_{\chi\in\ZZ} N_{\chi}$, where
$N_{\chi}\subset N$ is the subbundle on which $S^1$ acts with
weight $\chi$. Define the set of weights of $F'$ to be
$\weight(F'):=\{\chi\in\ZZ\mid N_{\chi}\neq 0\}$. Define also the
{\bf set of weights of $X$} to be
$$\weight(X):=\bigcup_{F'\subset F}\weight(F')\subset \ZZ,$$
where the union runs over the set of connected components of $F$.
Finally, define the set of {\bf critical residues} to be
$$\Lambda_{\crit}:=\{\lambda\in\imag\RR\mid
\text{there is some $w\in\weight(X)$ such that $w\lambda\in
\imag\ZZ$}\}.$$ For any $\lambda\in\imag\RR$ we will denote
$$X^{\lambda}:=\{x\in X\mid e^{2\pi\lambda}\cdot x=x\}.$$
Of course, for any $\lambda$ we have $F\subset X^{\lambda}$, and
the condition $\lambda\in\Lambda_{\crit}$ is equivalent to the
inclusion $F\subset X^{\lambda}$ being proper.

If $A$ is a meromorphic connection on a principal $S^1$ bundle
over a nodal curve $C$ and $z\in C$ is a node, we will say that
the {\bf holonomy of $A$ around $z$ is critical} if, denoting by
$y$ any preimage of $z$ in the normalisation of $C$, we have
$H:=\Hol(A,y)\in e^{2\pi\Lambda_{\crit}}$ (this is clearly
independent of the chosen preimage of the normalisation, and hence
well defined). In other words, the holonomy $H$ is critical if the
set of points in $X$ fixed by $H$ is bigger than $F$.

\subsection{Local behaviour of holomorphic sections near a marked point}

Let $(C,\bx)$ be a smooth marked curve and let $P$ be a principal
$S^1$ bundle over $C\setminus\bx$, endowed with a meromorphic connection $A$.
Let $Y:=P\times_{S^1}X$, let $x\in\bx$ be a marked point, and denote
by $D_{\epsilon}\subset C$ the disk of radius $\epsilon$ centered at $x$.
Following the same ideas as in Section \ref{ss:gluing-data}
we can define an equivalence relation on the restriction of
$Y$ on $D_{\epsilon}$ using parallel transport with respect
to $A$ along radial
directions. Taking the quotient by this equivalence relation
gives rise to a bundle $Y_x$ over $S_x$, which can be canonically
identified with $P_x\times_{S^1}X$.

Given a smooth section $\phi$ of $Y$, we say that $\phi$ {\bf extends
at $x$ to
give a section} $\phi_x$ of $Y_x$ if, for every $\theta\in S_y$ we have
$\lim_{\delta\to 0}[\phi(\exp_x\delta\theta)]=\phi_x(\theta)$, where
the brackets denote the equivalence class in $Y_x$.
In other words, this means that the pullback of the section $\phi$
to the real blowup of $C$ at $x$ extends to the exceptional divisor.

\begin{theorem}
\label{thm:extensio-pol} Suppose that a section $\phi$ of $Y$ satisfies
$\ov{\partial}_A\phi=0$  
and $\|d_A\phi\|_{L^2(C\setminus\bx)}<\infty.$ Then for every
marked point $x\in\bx$ the section $\phi$ extends at $x$ to give a
section $\phi_x$ of $Y_x$ and we have
$d_{A_x}\phi_x=0.$
Furthermore, the section $\phi_x$ takes values in $X^{\lambda}$.
In particular, if the holonomy of $A$ around $x$ is not critical,
then $\phi_x$ takes values in the fixed point set $F$, hence
$\phi_x$ is constant and the following limit exists
\begin{equation}
\label{eq:limit-point-global} \phi(x):=\lim_{z\to x}\phi(z)\in F.
\end{equation}
\end{theorem}

The result of Theorem \ref{thm:extensio-pol} is of
local nature, so it follows from the corresponding version
for the punctured disk, which is considered in Corollary
\ref{cor:existence-limit-orbit} (see Section \ref{s:exp-decay} for
the statement and the proof of Corollary
\ref{cor:existence-limit-orbit}).

\subsection{Chains of gradient segments}
\label{ss:gradient-chains} Recall that we denote $H:=-\imag\mu$.
Let $\xi_t:X\to X$ be the downward gradient flow at time $t$ of
$H$, so $\xi_0=\Id_X$ and
$$\frac{\partial\xi_s}{\partial t}
=-\xi_s^*\nabla H =-\xi_s^*I\XXX,$$ where $\XXX$ is the vector
field generated by the infinitesimal action of $\imag\in\Lie S^1$
on $X$. A {\bf pointed gradient segment} in $X$ is a pair $(x,T)$,
where $x$ is a point in $X$, not contained in the fixed point set,
and $T\subset\RR$ is a closed interval of positive measure. Define
$\xi_T(x):=\{\xi_t(x)\mid t\in T\}\subset X$. Two pointed gradient
segments $(x,T)$ and $(x',T')$ are said to be {\bf equivalent} if
$\xi_T(x)=\xi_{T'}(x')$. A {\bf gradient segment} in $X$ is an
equivalence class of pointed gradient segments. A {\bf chain of
gradient segments} is a finite sequence $\TTT$ of gradient
segments represented by a list of pointed gradient segments
$((x_1,T_1),\dots,(x_k,T_k))$ satisfying the following properties:
\begin{enumerate}
\item if $j>1$ then $\inf T_j=-\infty$, and if $j<k$ then $\sup
T_j=\infty$; \item if $1\leq j<k$ then
$\lim_{t\to\infty}\xi_t(x_j)=\lim_{l\to-\infty}\xi_l(x_{j+1}).$
\end{enumerate}
The {\bf beginnig} of $\TTT$ is the point $\lim_{l\to\inf
T_1}\xi_l(x_1),$ and the {\bf end} of $\TTT$ is $\lim_{t\to\sup
T_k}\xi_t(x_k).$ A {\bf degenerate chain of gradient segments} is
simply a point $x\in X$ (so this corresponds to the case $(x,T)$
where $T=\{0\}$).

Denote by $\TTT(X)$ the {\bf set of chains of gradient segments on
$X$}, including the degenerate ones. The group $S^1$ acts on
$\TTT(X)$ as follows: if $\theta\in S^1$ and $\TTT\in\TTT(X)$ is
represented by $((x_1,T_1),\dots,(x_k,T_k))$ then
$\theta\cdot\TTT$ is represented by $((\theta\cdot
x_1,T_1),\dots,(\theta\cdot x_k,T_k))$. The set $\TTT(X)$ carries
a natural topology induced by the Hausdorff distance between
subsets of $X$. With this topology, $\TTT(X)$ is clearly compact.

If $S$ denotes the circle and $P\to S$ is a principal $S^1$ bundle
provided with a connection $A$ with trivial holonomy, then we say
that a section $\TTT_S$ of the associated bundle
$P\times_{S^1}\TTT(X)$ is {\bf covariantly constant} (with respect
to $A$) if, given a trivialisation with respect to which $d_A=d$,
the map $S\to \TTT(X)$ given by the section $\TTT_S$ is constant.

\subsection{Metrics of fixed volume on stable curves}
\label{ss:universal-metrics}

Recall that the Deligne--Mumford moduli space $\oM_{g,n}$ of
isomorphism classes of stable curves $[C,\bx]$ admits a natural
structure of orbifold. The map $f:\oM_{g,n+1}\to\oM_{g,n}$ which
forgets the last point and stabilises gives $\oM_{g,n+1}$ the
structure of universal curve over $\oM_{g,n}$. Let $\Met_{g,n}$ be
the space of all smooth (in the orbifold sense) metrics in
$\oM_{g,n+1}$ whose restriction to the fibre of $f$ over $[C,\bx]$
gives a metric on $C$ of total volume $1$ and in the conformal
class defined by the complex structure. For any $\nu\in\Met_{g,n}$
and $[C,\bx]\in\oM_{g,n}$ we will denote by $\nu_{[C,\bx]}$ the
induced metric in $C$. If $\pi:C'\to C$ is the normalisation map,
the pullback $\pi^*\nu_{[C,\bx]}$ is a smooth metric in $C'$ (this
is true because $\nu$ is smooth in $\oM_{g,n+1}$). Also,
$\nu_{[C,\bx]}$ is invariant under the action of the automorphisms
of $(C,\bx)$. We give $\Met_{g,n}$ the obvious topology, which
makes it a contractible space.

\subsection{Definition of $c$-stable twisted holomorhic maps}
\label{ss:def-STHM} Let $g$ and $n$ be nonnegative integers
satisfying $2g+n\geq 3$. Take two natural numbers $n_{\crit}$ and
$n_{\gen}$ such that $n_{\crit}+n_{\gen}=n$.

Let $(C,\bx)$ be a nodal curve of genus $g$ and with $n$ marked
points. Repeatedly contracting the unstable components of $C$, we
obtain a stable curve $C^{\st}$ and a map $s:C\to C^{\st}$, called
the {\bf stabilization map}.
Let $C^b\subset C$ be the union of the irreducible components which
are contracted to a point by the stabilisation map, and let
$C^p\subset C$ be the union of the components not contained in $C^b$.
Then we have $C=C^p\cup C^b$. The components of
$C^b$ are called the {\bf bubble components} of $C$, and those of
$C^p$ the {\bf principal components} of $C$.
Finally, an {\bf exceptional point} of $C$ is a
point which is either a marked point or a node.

Pick a metric $\nu\in\Met_{g,n}$ and an element $c$ of $\imag\RR$.
A {\bf $c$-stable twisted holomorphic map} ($c$-STHM for short) of
genus $g$ and $n$ marked points is a tuple
$$\CCC=((C,\bx_{\crit},\bx_{\gen}),(P,A,G),\phi,\{\TTT_y\},\{\TTT_x\}),$$
where
\begin{enumerate}
\item $C$ is a connected compact nodal complex curve,
$\bx_{\crit}$ and $\bx_{\gen}$ are disjoint lists of smooth points
of $C$: $\bx_{\crit}$ is the list of {\bf critical marked points}
and $\bx_{\gen}$ is the list of {\bf generic marked points};
$\bx_{\crit}$ contains $n_{\crit}$ points and $\bx_{\gen}$
contains $n_{\gen}$ points, and we denote by $\bx$ the union
$\bx_{\crit}\cup\bx_{\gen}$. \item $P$ is a principal $S^1$ bundle
on the set of non exceptional points of $C$,
$$P\to C\setminus (\bx\cup \bz),$$ where $\bz\subset C$ is the set
of nodes. \item $A$ is a meromorphic connection on $P$, $G$ is a
choice of gluing data for $P$ and $\phi$ is a section of the
bundle $P\times_{S^1}X$. \item  For each preimage $y$ of a node in
$C$, $\TTT_y$ is a covariantly constant section of the bundle
$P_y\times_{S^1}\TTT(X^{\lambda})$, where $\lambda\in\imag\RR$ is
such that $\Hol(A,y)=e^{2\pi\lambda}$ (this gives a chain of
gradient lines for each tangent direction at $y$ in the
normalisation of the curve, varying in a $S^1$-equivariant way).
\item For each $x\in\bx_{\gen}$, $\TTT_x$ is a covariantly
constant section of the bundle $P_x\times_{S^1}\TTT(X^{\lambda})$,
where $\lambda\in\imag\RR$ is such that
$\Hol(A,x)=e^{2\pi\lambda}$.
\end{enumerate}

The tuple $\CCC$ must satisfy the following conditions:
\begin{enumerate}
\item {\bf The section is holomorphic.} The section $\phi$
satisfies the equation
\begin{equation}
\label{eq:phi-holomorfa} \ov{\partial}_A\phi=0.
\end{equation}

\item {\bf Vortex equation.} For any principal component
$C_j\subset C^p$, let $\nu_j$ (resp. $A_j$, $\phi_j$) be the
restriction of $s^*\nu_{[C,\bx]}$ (resp. $A$, $\phi$) to $C_j$,
where $s$ is the stabilisation map, and let $d\vol(\nu_j)$ be
the induced volume form; then
\begin{equation}
\label{eq:theequations}
\iota_{d\vol(\nu_j)}F_{A_j}+\mu(\phi_j)=c;
\end{equation}
we call this equation the vortex equation because in the case of
$X=\CC$ with the action of $S^1$ of weight $1$ this equation
coincides with the standard abelian vortex equation.

\item {\bf Flatness on bubbles.} The restriction of $A$ to each
bubble component is flat.

\item {\bf Finite energy.} The energy of $\phi$ as a section is
bounded:
\begin{equation}
\label{eq:energia-finita}
\|d_A\phi\|_{L^2}<\infty;
\end{equation}

\item {\bf Matching condition at the nodes.} Let $z\in C$ be a
node, and let $y,y'$ be its preimages in the normalization map. By
Theorem \ref{thm:extensio-pol}, (\ref{eq:phi-holomorfa}) and
(\ref{eq:energia-finita}) imply that $\phi$ extends to give
sections of $\phi_y\in\Gamma(P_{y}\times_{S^1}X)$ and
$\phi_{y'}\in\Gamma(P_{y'}\times_{S^1}X)$. Let $\rho:P_y\to
P_{y'}$ be the isomorphism given by the gluing data $G$. Then
$$\TTT_{y}=\rho^*\TTT_{y'}$$ and for
every $\theta\in S_y$, $\phi_y(\theta)$ has to be either the
beginning or the end of the chain $\TTT_y(\theta)$, and
$\rho^*\phi_{y'}(\theta)$ has to be the opposite extreme.
Furthermore, if the holonomy of $A$ around $y$ is not critical,
then the chain of gradient segments $\TTT_y$ has to be degenerate,
which implies that $\rho^*\phi_{y'}=\phi_y$.

\item{\bf Matching condition and the generic marked points.} Given
$x\in\bx_{\gen}$, let $\phi_x\in\Gamma(P_x\times_{S^1}X)$ be the
extension of $\phi$. For any $\theta\in S_x$, $\phi_x(\theta)$ is
either the beginning or the end of $\TTT_x(\theta)$, and the
opposite extreme of $\TTT_x(\theta)$ is a fixed point.
Furthermore, if the holonomy of $A$ around $x$ is not critical
then $\TTT_x$ has to be degenerate.

\item {\bf Stability condition for the bubbles.} If $C'\subset C$
is a bubble component with less than $3$ exceptional points, then
the restriction of $d_A\phi$ to $C'$ is not identically zero.
\end{enumerate}

Two $c$-STHM's $\CCC$ and $\CCC'$ are said to be {\bf isomorphic}
if there is a commuting diagram
$$\xymatrix{P\ar[r]^{g}\ar[d] & P'\ar[d] \\
C\ar[r]^{f} & C',}$$ where $g$ is an isomorphism of principal
$S^1$ bundles and $f$ is a biholomorphism satisfying
$f(\bx_{\crit})=\bx'_{\crit}$, $f(\bx_{\gen})=\bx'_{\gen}$ and
preserving the ordering of the marked points, such that
$$g^*A'=A,\qquad
g^*G'=G,\qquad g^*\phi'=\phi, \qquad g^*\TTT_y'=\TTT_y
\qquad\text{and}\qquad g^*\TTT_x'=\TTT_x.$$

\subsection{Remarks on the definition of
Hamiltonian Gromov--Witten invariants}

We now make a few comments with the hope of clarifying some
of the ingredients appearing in the definition of $c$-STHM.
We recall that the main application we have in mind of the
compactness result proved in this paper is
the definition of invariants of the manifold $X$ and the
Hamiltonian action of $S^1$, using the moduli space of $c$-STHM's.

First of all, in this paper we have imposed no restriction on
$c\in\imag\RR$. However, when considering the moduli space of
$c$-STHM's we will need to take $c$ away from a discrete set of
critical values. The invariants obtained from two different
choices of $c$ may vary if we cross critical values when passing
from one choice to the other (this is explained in \cite{M}, and
is similar to the well known phenomenon of wall crossing in gauge
theories). On the other hand, in order to define $c$-STHM's we
have made a choice of an element $\nu\in\Met_{g,n}$; the
invariants which we will construct do not depend on this choice,
thanks to a standard cobordism argument and the fact the set
$\Met_{g,n}$ is connected. This is analogous to the fact that
Gromov--Witten invariants are independent of the chosen compatible
almost complex structure. In our situation, the obtained
invariants will also be independent of the $S^1$-invariant and
compatible almost complex structure. Finally, the reason why we
distinguish two different kinds of marked points (critical and
generic) is the following: when constructing the moduli space, we
will allow the residue at generic marked points to vary, whereas
the residue at a critical marked point will always have to be
critical. The way we define the evaluation map at a marked point
will depend on whether the point is critical or generic.

\section{Topology on the set of $c$-STHM}
\label{s:topology-STHM}

\subsection{Convergence of lines to chains of gradient lines}
\label{ss:convergence-chain-grad-lines}
 Let $F\subset X$ be the
fixed point set. For any small $\delta>0$, denote by $F^{\delta}$
the $\delta$-neighbourhood of $F$, that is, the set of points of
$X$ at distance $\leq\delta$ from $F$. Let also
$X^{\delta}:=X\setminus F^{\delta}$. Finally, recall that we
denote by $H$ the function $-\imag\mu$.

Let $\{f_u:S_u\to X\}$ be a sequence of smooth maps, where each
$S_u\subset\RR$ is a closed interval. We call the pairs
$(f_u,S_u)$ {\bf lines in $X$}. Let $\TTT$ be a chain of gradient
flow lines (see Section \ref{ss:gradient-chains}) represented by a
list $((x_1,T_1),\dots,(x_k,T_k))$ of pointed gradient segments.
Let $F_1,\dots,F_l$ be the connected components of the fixed point
set which intersect the closure of $\bigcup \xi_{T_k}(x_k)$,
labelled in such a way that $H(F_1)>H(F_2)>\dots>H(F_l).$ We will
say that {\bf the sequence of lines $\{(f_u,S_u)\}$ converges to
$\TTT$} if for each small enough $\delta>0$ and any big enough $u$
we can write $S_u$ as a union of two sets
\begin{equation}
\label{eq:descomposicio-S} S_u=T_u^{\delta}\cup E_u^{\delta},
\end{equation}
in such a way that:
\begin{enumerate}
\item each $E_u^{\delta}$ is a union of $l$ closed intervals:
$E_u^{\delta}=E_{u,1}^{\delta}\cup\dots\cup E_{u,l}^{\delta},$ and
each $T_u^{\delta}$ is a union of $k$ closed intervals:
$T_u^{\delta}=T_{u,1}^{\delta}\cup\dots\cup T_{u,k}^{\delta},$
labelled in such a way that $\sup E_{u,j-1}^{\delta}< \inf
E_{u,j}^{\delta}$ and $\sup T_{u,j-1}^{\delta}\leq \inf
T_{u,j}^{\delta}$ for every $j$ for which the expression makes
sense; \item denoting by $t_{u,j}^{\delta}$ any point of
$T_{u,j}^{\delta}$, the maps $f_u:T_{u,j}^{\delta}\to X$
approximate gradient segments
\begin{equation}
\label{eq:f_u-quasi-gradient}\lim_{u\to\infty} \sup_{t\in
T_{u,j}^{\delta}}
d(\xi_{t-t_{u,j}^{\delta}}(f_u(t_{u,j}^{\delta})),f_u(t))=0;
\end{equation}
\item the images $f_u(T_{u,j}^{\delta})$ approximate in
$X^{\delta}$ the $j$-th gradient segment of $\TTT$:
$$\lim_{\delta\to 0}\lim_{u\to\infty}
D(f_u(T_{u,j}^{\delta},\xi_{T_j}(x_j)\cap X^{\delta})=0,$$ where
here $D$ denotes de Hausdorff distance between sets;
 \item the images of the
sets $E_{u,j}^{\delta}$ become smaller and smaller as $\delta\to
0$ and they accumulate near the fixed point set:
$$\lim_{\delta\to 0}\left(
\limsup_{u\to\infty}\diam(f_u(E_{u,j}^{\delta}))\right)
=\lim_{\delta\to 0}\left(\limsup_{u\to\infty}
d(f_u(E_{u,j}^{\delta}),F_j)\right)=0.$$
\end{enumerate}

This implies in particular that the sets $f_u(S_u)$ converge in
the Hausdorff metric to $\bigcup\xi_{T_j}(x_j)$.

\subsection{Convergence with gauge $\lambda$ of cylinders to chains
of gradient lines} Now suppose that $\{\phi_u:C_u\to X\}$ is a
sequence of smooth maps, where each $C_u$ is a cylinder
$C_u=S_u\times S^1$ and $S_u\subset\RR$ is a closed interval. We
call the pairs $(\phi_u,C_u)$ {\bf cylinders in $X$}. Let
$\lambda\in\Lambda_{\crit}$ be a critical residue. We say that
{\bf the sequence of cylinders $\{\phi_u,C_u\}$ converges with
gauge $\lambda$ to a chain of gradient lines} $\TTT$ in
$X^{\lambda}$ if there is a sequence of lines in $X$,
$\{(\psi_u,S_u)\}$, which converges to $\TTT$ and such that
$$\lim_{u\to\infty}\sup_{(t,\theta)\in
C_u}d(e^{\lambda\theta}\phi_u(t,\theta),\psi_u(t))=0.$$ Here
$e^{\lambda\theta}$ denotes {\it any} number of the form
$e^{\lambda\ov{\theta}}$, where $\ov{\theta}\in\RR$ is a lift of
$\theta\in\RR/2\pi\ZZ$. Since the chain $\TTT$ is contained in
$X^{\lambda}$, the resulting notion of convergence is independent
of the chosen lifts.

\subsection{Description of the topology of $\oM_{g,n}$}
\label{ss:convergencia-corbes}

We give here a description of the topology of the Deligne--Mumford
moduli space which is suitable for our purposes. For that it
suffices to specify what it means that a sequence of stable
(pointed) curves converges to a given curve.

Let $(C,\bx)$ be a stable curve, denote by $\bz$ the set of nodes,
let $\pi:C'\to C$ be the normalization and let $\by$ be the
preimages of $\bz$ by $\pi$. Pick a conformal metric on $C'$ and
let $\epsilon>0$ be a small number. For each $y\in\by$, take a
neighborhood $U_y\subset C'$, small enough so that it admits a
biholomorphism $\zeta_y:D(\epsilon)\to U_y$ with the disk
$D(\epsilon)\subset T_{y}C'$ centered at $0$ and of radius
$\epsilon$. Chose also a small neighborhood $B_x\subset C'$ of the
preimage of each marked point $x\in\bx$. Denote by $I\in\End TC'$
the complex structure of $C'$, which we now view as a compact real
surface.

Let now $I'\in \End TC'$ be another complex structure which
coincides with $I$ on each $U_y$ and each $B_x$ (call such complex
structure {\bf admissible}). Take also, for each node $z$ with
preimages $y,y'$ in $C'$, an element $\delta_z\in T_{y}C'\otimes
T_{y'}C'$ satisfying $|\delta|<\epsilon^2$. We call the collection
of numbers $\{\delta_z\}$ {\bf smoothing parameters}. Define a new
curve $C(I',\{\delta_z\})$ as follows: replace the complex
structure $I$ by $I'$, then remove for each for each pair $y,y'$
of preimages of a node $z$, the sets
$\zeta_y(D(|\delta_z|/\epsilon))$ and
$\zeta_{y'}(D(|\delta_z|/\epsilon))$ from $C'$, and finally
identify for each pair of elements $u \in D(\epsilon)\setminus
D(\delta_z/\epsilon)$ and $v \in D(\epsilon)\setminus
D(\delta_z/\epsilon)$ satisfying $u \otimes v=\delta_z$, the
images $\zeta_y(u)$ and $\zeta_{y'}(v)$ (such identifications
preserve the complex structure because $I'$ is admissible).

For later use, define for every $z$ the subset
$N_y(\delta_z)\subset C(I',\{\delta_z\})$ to be the image by
$\zeta_y$ of the annulus $D(\epsilon)\setminus
D(|\delta_z|/\epsilon)$ (this is equal to the image by
$\zeta_{y'}$ of the corresponding subset of $T_{y'}C'$). The set
$N_z$ is conformally equivalent to the cylinder: $[\ln
|\delta_z|-\ln\epsilon,\ln\epsilon]\times S^1$. We will say that
$y$ {\bf is in the side of} $\{\ln\epsilon\}\times S^1$, and that
$y'$ is in the side of $\{\ln|\delta_j|-\ln\epsilon\}\times S^1$
(if we consider $N_{y'}(\delta_z)$ instead, then the roles are
inverted).

Note that for any compact set $K\subset C\setminus (\bx\cup\bz)$
and small enough smoothing data $\{\delta_z\}$ there is a
canonical inclusion $K\to C(I',\{\delta_z\})$, which we will
denote by $\iota_K$. (For $\iota_K$ to exist it suffices to take
each $\delta_z$ so that $\zeta_y(D(\delta_z/\epsilon))$ and
$\zeta_{y'}(D(\delta_z/\epsilon))$ are disjoint from $K$.)

A sequence $\{(C_u,\bx_u)\}$ of stable curves {\bf converges} to
$(C,\bx)$ if for big enough $u$ there is an admissible complex
structure $I_u$, smoothing parameters $\{\delta_{u,z}\}$, and an
isomorphism of marked nodal curves
\begin{equation}
\label{eq:iso-corbes} \xi_u:(C(I_u,\{\delta_{u,z}\}),\bx)\to
(C_u,\bx_u),
\end{equation}
such that $I_u$ converges to $I$ in $C^{\infty}(\End TC')$ and for
each node $z$ we have $\delta_{u,z}\to 0$.

The topology on $\oM_{g,n}$ defined by this notion of convergence
coincides with the usual one (see for example Section 9 in
\cite{FO}, where the topology of $\oM_{g,h}$ is described in
similar terms).

\subsection{Connections in balanced temporal gauge}
\label{ss:balanced-temporal} Let $C=[p,q]\times S^1$ be a
cylinder, and denote by $(t,\theta)$ the usual coordinates. Let
$d_A=d+\alpha$ be a covariant derivative on the trivial
principal $S^1$ bundle over $C$. We will say that $\alpha$ is in
{\bf balanced temporal gauge} if it is in temporal gauge, so that
$\alpha=a d\theta$ for some function $a:C\to\imag\RR$, and
furthermore the restriction of $a$ to the middle circle
$\{(p+q)/2\}\times S^1$ is constantly equal to some
$\lambda\in\imag\RR$, which is called the {\bf residue} of $A$
(with respect to the trivialization). Any connection on the
trivial bundle over $C$ is gauge equivalent to a connection in
balanced temporal gauge. Furthermore, since
$d\alpha=\frac{\partial a}{\partial t}dt\wedge d\theta$, we have
the estimate:
\begin{equation}
\label{eq:bound-balanced}
|a(t,\theta)-\lambda|\leq\left|\int_{(p+q)/2}^t
|d\alpha(\tau,\theta)|d\tau\right|.
\end{equation}

\subsection{Convergence of $c$-STHM}
\label{ss:convergence-STHM}
 Our aim here is to define a topology
on the set of isomorphism classes of $c$-STHM's specifying as
before what it means for a sequence of $c$-STHM's to converge to a
given $c$-STHM's. Before defining the convergence of sequences, we
make the observation that the notion of convergence for stable
curves given in Section \ref{ss:convergencia-corbes} makes perfect
sense when considering nodal marked curves in general: if
$(C,\bx)$ is a nodal marked curve with $k$ nodes, we can define as
before the deformations $C(I_u,\delta_1,\dots,\delta_k)$. (Of
course, the topology induced by this notion on the set of
isomorphism classes of nodal marked curves is not Hausdorff.)

For simplicity, we will only define convergence of sequences of
$c$-STHM's with smooth underlying marked curve and with degenerate
chains of gradient segments at generic marked points. To pass from
this to the general case is routine.

So let $\{\CCC_u\}$ be a sequence of $c$-SHTC's. Suppose that
$(C_u,\bx_u)$ is the smooth marked curve underlying $\CCC_u$, and
that $\bx_u=\bx_{\crit,u}\cup\bx_{\gen,u}$, $(P_u,A_u)$ is the
bundle and connection on $C_u$ (since $C_u$ is smooth there is no
gluing data) and $\phi_u$ is the section of $P_u\times_{S^1}X$
(again, since $C_u$ is smooth there are no chains of gradient
segments $\TTT_y$); finally, for each $x\in\bx_{\gen}$ the chain
$\TTT_x$ is degenerate. Let now
$$\CCC=((C,\bx_{\crit},\bx_{\gen}),(P,A,G),\phi,\{\TTT_y\},\{\TTT_x\})$$
be another $c$-STHM, not necessarily with smooth underlying curve,
and let $\bz\subset C$ be the nodes of $C$.

We will say that the sequence of isomorphism classes
$\{[\CCC_u]\}$ converges to $[\CCC]$ if for any exhaustion
$K_1\subset\dots\subset K_l\subset\dots$ of $C\setminus
(\bx\cup\bz)$ by compact subsets the following holds.

\begin{enumerate}
\item {\bf Convergence of the underlying curves.} The curves
$(C_u,\bx_{\crit,u},\bx_{\gen,u})$ converge to
$(C,\bx_{\crit},\bx_{\gen})$. This implies that there are
isomorphisms
$$\xi_u:(C(I_u,\{\delta_{u,z}\}),\bx_{\crit},\bx_{\gen})\to
(C_u,\bx_{\crit,u},\bx_{\gen,u})$$ such that $I_u\to I$ and
$\delta_{u,z}\to 0$. Pulling back everything by $\xi_u$ we can
assume that the underlying curve of $\CCC_u$ is
$C(I_u,\{\delta_{u,z}\},\bx_{\crit},\bx_{\gen})$.

\item {\bf Convergence of gluing angles.} Let
$((\gamma_1,\rho_1),\dots,(\gamma_k,\rho_k))$ be the glueing data
at each of the nodes of $C$ given by $G$. Since we assume that
each $C_u$ is smooth, any gluing parameter $\delta_{u,j}$ is
nonzero and hence gives rise to a gluing angle $[\delta_{u,j}]\in
\Gamma_{z_j}$. Then, for any $j$ we must have
$[\delta_{u,j}]\to \gamma_j$.

\item {\bf Convergence of the connections and sections away from
the nodes.} For each $K_l$ and any big enough $u$ (so that
$\iota_{K_l}$ is defined) there must exist an isomorphism of
vector bundles
$$\rho_{u,l}:P|_{K_l}\to \iota_{K_l}^*P_u$$
such that $\rho_{u,l}^*\iota_K^*A_u$ converges to $A$ and
$\rho_{u,l}^*\iota_K^*\phi_u$ converges to $\phi$ on $K_l$ as $u$
goes to $\infty$ (here the convergence is assumed to be in
$C^{\infty}$). This implies that the holonomies of $A_u$ around
the marked points converge to those of $A$.

 \item {\bf Convergence of
gluing data.} The isomorphisms $\rho_{u,l}$ have to satisfy the
following additional condition. Suppose that $r>0$ is smaller than
the $\epsilon$ used in the definition of convergence of stable
curves. Take some node $z\in\bz$ with preimages $y,y'$. Let
$Y(r):=\zeta_y(S(r))$, where $S(r)\subset T_{y}C'$ is the circle
of radius $r$ centered at $0$, and define $Y'(r)$ similarly.
Suppose that $l$ is big enough so that both $Y(r)$ and $Y'(r)$ are
contained in $K_l$, and suppose that $u$ is big enough so that
$\iota_{K_l}$ exists. Define a map $\tau_{u,l}(r)$ by the
condition that the following diagram commutes:
$$\xymatrix{(\iota_{K_l}^*P_u)|_{Y(r)}\ar[r]^{\tau_{A_u}}
& (\iota_{K_l}^*P_u)|_{Y'(r)} \\
P|_{Y(r)}\ar[r]^{\tau_{u,l}(r)} \ar[u]^{\rho_{u,z}} &
P|_{Y'(r)}\ar[u]_{\rho_{u,l}},}$$ where $\tau_{A_u}$ denotes the
parallel transport along the images by $\zeta_y$ of the radial
directions in the annulus $D(\epsilon)\setminus
D(\delta_z/\epsilon)$. Finally, let $f':P|_{Y'(r)}\to P_{y'}$ and
$f:P|_{Y(r)}\to P_{y}$ be the natural projections. Then the
composition $f'\circ\tau_{u,l}(r)\circ f$ gives an isomorphism
between $P_{y}$ and $P_{y'}$ lifting the isometry between $S_{y}$
and $S_{y'}$ given by the gluing angle $[\delta_z]$, so it
specifies an element $g_{u,l,j}(r)\in\Gamma(P,z)$. The condition
is that for every $j$
$$\lim_{r\to 0}
\left(\limsup_{u\to \infty}
\limsup_{l\to\infty}d(g_{u,l,j}(r),(\gamma_j,\rho_j))\right)=0,$$
where $d$ denotes a fixed distance function defined on $\Gamma(P,z_{j})$.

\item {\bf Convergence near the nodes to chains of gradient
segments.} Fix some node $z\in C$ with preimages $y,y'$ in the
normalization. Suppose that the holonomy of $A$ around $y$ is
$e^{2\pi\lambda}$ for some $\lambda\in\imag\RR$. We will denote
$\delta_u:=\delta_{u,z}$. To specify the condition, we will use
the cylinders
$$N_u:=N_y(\delta_{u})\simeq [\ln
|\delta_u|-\ln\epsilon,\ln\epsilon]\times S^1$$ defined in Section
\ref{ss:convergencia-corbes}. Define, for any big $\Delta>0$ and
$u$ the cylinder
$$N_u(\Delta)=
[\ln |\delta_u|-\ln\epsilon+\Delta,\ln\epsilon-\Delta]\times
S^1\subset N_u.$$ We distinguish two situations.
\begin{enumerate}
\item If the chain of gradient segments $\TTT_{y}$ is degenerate
(for example, if the residue $\lambda$ is not critical), then we
must have
$$\lim_{\Delta\to 0}\left(
\limsup_{u\to\infty}\diam(\phi_u(N_u(\Delta))\right)=0.$$ \item
Otherwise, the following must happen; for each big enough $u$
there is a cylinder $M_u\subset N_u$ satisfying two properties.
First, for any $\Delta>0$ and big enough $u$, we must have
$M_u\subset N_u(\Delta)$. Then the complementary is the union of
two cylinders $M_u^{\Delta,-}$ and $M_u^{\Delta,+}$. The second
property is the following. Suppose that $P_u$ is trivialized on
$N(\delta_u)$ in such a way that $d_{A_u}=d+\alpha_u$ is in
balanced temporal gauge and with residue $\lambda_u$, and that
$\lambda_u\to\lambda$. Then the section $\phi_u$ restricts to give
a map $\phi_u:N_u\to X$, and we must have:
\begin{enumerate}
\item the images of the sets $M_u^{\delta,\pm}$ have smaller and
smaller diameter:
$$\lim_{\Delta\to\infty}\left(\limsup_{u\to\infty}
\diam_{S^1}(\phi_u(M_u^{\Delta,\pm}))\right)=0;$$ \item
$\lambda_u\neq\lambda$ for big enough $u$, so the number
$l_u:=-\imag(\lambda_u-\lambda)$ is nonzero; \item the sequence of
cylinders $\{\phi_u/l_u:l_uM_u\to X\}$ converges with gauge
$\lambda$ to the chain of gradient segments $\TTT_{y_j}$ (here, if
$M_u=[p,q]\times S^1$ then $l_uM_u=[l_up,l_uq]$ and
$\phi_u/l_u(t,\theta):=\phi_u(t/l_u,\theta)$).
\end{enumerate}
\end{enumerate}

\item {\bf Convergence near generic marked points to chains of
gradient segments.} This is similar to the previous condition. Let
$x\in\bx_{\gen}$. Let $N$ be a punctured neighborhood of $x$
biholomorphic to $[0,\infty)\times S^1$. Define also
$N({\Delta}):=[\Delta,\infty)$. We distinguish two cases.
\begin{enumerate}
\item If the chain $\TTT_x$ is degenerate then
$$\lim_{\Delta\to 0}\left(
\limsup_{u\to\infty}\diam(\phi_u(N(\Delta))\right)=0.$$

\item Suppose now that $\TTT_x$ is not degenerate. Let the
holonomy of $A$ around $x$ be $e^{2\pi\lambda}$. Then we must have
$\lambda\in\Lambda_{\crit}$. There must be, for each big enough
$u$, a cylinder $M_u=[\Delta_u,\infty)$ satisfying:
\begin{enumerate}
\item $$\lim_{\Delta\to\infty}\left(\limsup_{u\to 0}
\diam_{S^1}(\phi_u([\Delta,\Delta_u]\times S^1))\right)=0;$$ \item
take a trivialisation of $P_u$ with respect to which
$d_{A_u}=d+\alpha_u$ is in temporal gauge and the restriction of
$\alpha_u$ to $\{\Delta_u\}\times S^1$ is equal to $\lambda_u
d\theta$, where $\lambda_u\in\imag\RR$ is a constant, and
$\lambda_u\to\lambda$. We must have $\lambda_u\neq\lambda$ for big
enough $u$, so the number $l_u:=-\imag(\lambda_u-\lambda)$ is
nonzero; \item the sequence of cylinders $\{\phi_u/l_u:l_uM_u\to
X\}$ converges with gauge $\lambda$ to the chain of gradient
segments $\TTT_x$.
\end{enumerate}
\end{enumerate}
\end{enumerate}

\section{The Yang--Mills--Higgs functional $\YMH_c$}
\label{s:YMH}

Let
$\CCC=((C,\bx_{\crit},\bx_{\gen}),(P,A,G),\phi,\{\TTT_y\},\{\TTT_x\})$
be a $c$-STHM. Let $C^p$ (resp. $C^b$) be the union of the
principal components (resp. bubble components) of $C$. Let $s:C\to
C^{\st}$ be the stabilization map. Take on $C^p$ the conformal
metric $s^*\nu_{[C^{\st},\bx]}$ and extend it to a conformal
metric on $C$ by chosing an arbitrary conformal metric on the
bubbles. Define
$$\YMH_c(\CCC):=\|F_A\|^2_{L^2(C^p)}
+\|\mu(\phi)-c\|^2_{L^2(C^p)} +\|d_A\phi\|^2_{L^2(C)}.$$ We call
$\YMH_c$ the {\bf Yang--Mills--Higgs functional}. Its value is
independent of the conformal metric chosen in the bubbles, since
over the bubbles we only integrate the energy
$\|d_A\phi\|^2_{L^2(C)}$, which is conformaly invariant. In this
section we will compute $\YMH_c(\CCC)$ in terms of topological
data and the residues of the meromorphic connection $A$.

\subsection{Orbifold structure induced by the critical residues}

Let $(C,\bx)$ be a smooth marked curve, let $P\to C\setminus\bx$
be a principal $S^1$ bundle and let $A$ be a meromorphic conection
on $P$. We can divide the set of marked points as
$\bx=\bx_{\crit,A}\cup\bx_{\gen,A}$, where $\bx_{\crit,A}$ (resp.
$\bx_{\gen,A}$) denotes the set of points around which the holonomy
of $A$ is critical (resp. noncritical, which we call
generic). Since critical residues are of the form $\imag l$ with
$l\in\QQ$, we can use the construction
in Section \ref{ss:mer-conn-orbibundles} to obtain, for any
$\tau\in T(P,\bx_{\crit,A})$, an orbibundle
$\bP^{\tau,A}=(P^{\tau},\{P_j\},\{\psi_j\})$ on
$\bC^A:=(C,\bx_{\crit,A},\bq)$ (where $\bq$ denotes the collection
of denominators in the set of critical residues of $A$). The
following result gives a crucial property of holomorphic sections
of $P\times_{S^1}X$ with bounded energy.

\begin{theorem}
\label{thm:rem-singularities} Suppose that $\phi$ is a section of
$P\times_{S^1}X$ which satisfies $\ov{\partial}_A\phi=0$ and
$\|d_A\phi\|_{L^2(C\setminus\bx)}<\infty$. Then, for any $\tau\in
T(P,\bx_{\gen,A})$ the section $\phi$ extends to give a section
$\Phi^{\tau,A}$ of the orbibundle $\bP^{\tau,A}$.
\end{theorem}
\begin{pf}
The fact that $\phi$ extends to the points in $\bx_{\gen,A}$ follows
from Theorem \ref{thm:extensio-pol} (since in this case the
limiting orbit has to be a fixed point, there is indeed an
extension no matter what trivialisation $\tau$ we chose). The fact
that $\phi$ extends to $\bx_{\crit,A}$ is a consequence of Gromov
removal of singularities theorem. Indeed, pulling back $(A,\phi)$
to $\oU_j\setminus\{\ox_j\}$ and applying a suitable gauge
transformation $g$ the connection extends to the whole $\oU_j$,
hence defines a complex structure on $\oU_j\times X$. The
resulting map $g\cdot\rho_j^*\phi:\oU_j\setminus\{\ox_j\}\to
\oU_j\times X$ is complex and has bounded energy, so Gromov's
theorem applies (see for more detais the proof the case of
critical weight of Theorem \ref{thm:exp-decay} in Section
\ref{proof:thm:exp-decay}).
\end{pf}

\subsection{Computation of the Yang--Mills--Higgs functional}

\begin{lemma}
\label{lemma:integral-omega-A} Let $(C,\bx)$ be a smooth marked
curve, let $P\to C\setminus\bx$ be an principal $S^1$ bundle, let
$A$ be a meromorphic connection on $P$ and let $\phi$ be a section
of $P\times_{S^1}X$. Suppose that $\ov\partial_A\phi=0$ and that
$\|d_A\phi\|_{L^2}<\infty$. Then, for any $\tau\in
T(P,\bx_{\gen,A})$, we have
$$\int (\Phi^{\tau,A})^*[\omega(A)]=\frac{1}{2}\|d_A\phi\|_{L^2}^2
-\int F_A\mu(\phi)
-\sum_{x\in\bx_{\gen,A}}2\pi\Res(A,x,\tau_x)\mu(\phi(x)).$$
\end{lemma}

\begin{remark}
Note that, since $A$ might have nontrivial poles, the form
$\omega(A)$ is singular on $\bP^{\tau,A}\times_{S^1}X$. However,
one can define the cohomology class $[\omega(A)]$ simply by taking
any smooth connection $A'$ and setting
$[\omega(A)]:=[\omega(A')]$. By Remark \ref{rmk:coh-class-omega},
this is independent of $A'$.
\end{remark}

\begin{pf}
The proof is similar to that of Chern--Weil formula for
meromorphic connections in Lemma \ref{lemma:ChernWeil}. For any
small $\epsilon>0$ we modify $A$ near each point $x\in\bx_{\gen,A}$
as follows: if in a neighbourhood of $x$ we have
$d_A=d+\alpha+\lambda d\theta$ (where $\alpha$ extends
continuously to $x$ and $\lambda=\Res(A,x,\tau_x)$), then we
replace $\lambda d\theta$ by $\eta_{\epsilon}$, where
$\eta_{\epsilon}$ extends smoothly to $x$, it coincides with
$\lambda d\theta$ away from the disk $B(x,\epsilon)$ of radius
$\epsilon$ centered at $x$, and both the supremum of the norm
$|\eta_{\epsilon}|$ and the integral of $|d\eta_{\epsilon}|$ on
$B(x,\epsilon)$ is bounded above by a constant $C$ independent of
$\epsilon$ (for example, one can define $\eta_{\epsilon}$ to be
$r^2\lambda d\theta/\epsilon^2$ in $B(x,\epsilon-\epsilon^2)$).
Let $A_{\epsilon}$ be the resulting smooth connection on
$\bP:=\bP^{\tau,A}$ over the orbisurface $\bC=\bC^A$. Then we have
for any $\epsilon>0$
$$\int (\Phi^{\tau,A})^*[\omega(A)]
=\int (\Phi^{\tau,A})^*[\omega(A_{\epsilon})]
=\frac{1}{2}\|d_{A_\epsilon}\phi\|_{L^2}^2-\int
F_{A_\epsilon}\mu(\phi).$$ This follows from formula
(\ref{eq:def-omega-A}) and the fact that (taking local
trivialisations for (\ref{eq:def-omega-A}) to make sense) for any
point $x\in C$ and any choice of conformal metric on $C$ we have
$$\phi(x)^*(\pi_X^*\omega
-\pi^*\alpha_{\epsilon}\wedge\pi^*_X\iota_{\XXX}\omega)
=\frac{1}{2}(|\partial_{A_{\epsilon}}\phi|^2
-|\ov{\partial}_{A_{\epsilon}}\phi|^2)\ dvol(x).$$ To finish the
proof, we make two observations. First, that as $\epsilon\to 0$,
$\|d_{A_\epsilon}\phi\|_{L^2}^2$ converges to
$\|d_A\phi\|_{L^2}^2$. This is a consequence of the bound
$|\eta_{\epsilon}|<C$ and of Theorem \ref{thm:exp-decay}. Second,
we similarly have, as $\epsilon\to 0$,
$$\int F_{A_\epsilon}\mu(\phi)\to
\int F_A\mu(\phi)+\sum_{x\in\bx_{\gen,A}}2\pi\Res(A,x,\tau_x).$$
Indeed, setting $C_{\epsilon}:=C\setminus \bigcup_{x\in
\bx_{\gen,A}} B(x,\epsilon)$ we have
$$\int F_{A_\epsilon}\mu(\phi)=
\int_{C_\epsilon} F_A\mu(\phi) +\sum_{x\in\bx_{\gen,A}}
\int_{B(x,\epsilon)} F_{A_{\epsilon}}\mu(\phi).$$ The first
integral on the right hand side clearly converges to $\int_C
F_A\mu(\phi)$. To estimate the second term we compute for any
$x\in\bx_{\gen,A}$
$$\int_{B(x,\epsilon)} F_{A_{\epsilon}}\mu(\phi)
=\int_{z\in B(x,\epsilon)} F_{A_{\epsilon}}(z)\mu(\phi(x)) dz
+\int_{z\in B(x,\epsilon)}
F_{A_{\epsilon}}(z)(\mu(\phi(z))-\mu(\phi(x))) dz.$$ The last
integral converges to $0$ as $\epsilon\to 0$ because
$\int_{B(x,\epsilon)}|d\eta_{\epsilon}|<C$ and because
$\mu(z)\to\mu(x)$ as $z\to x$ (this follows from Corollary
\ref{cor:existence-limit-orbit}). Finally, we compute
\begin{align*}
\int_{z\in B(x,\epsilon)} F_{A_{\epsilon}}(z)\mu(\phi(x)) dz
&=\int_{z\in B(x,\epsilon)} d\eta_{\epsilon}\mu(\phi(x)) dz
+\int_{z\in B(x,\epsilon)} d\alpha\mu(\phi(x)) dz \\
&=\int_{\partial B(x,\epsilon)} \eta_{\epsilon}\mu(\phi(x)) dz
+\int_{z\in B(x,\epsilon)} d\alpha\mu(\phi(x)) dz \\
&=2\pi\Res(A,x,\tau_x)\mu(\phi(x)) +\int_{z\in B(x,\epsilon)}
d\alpha\mu(\phi(x)) dz,
\end{align*}
and the last integral converges to $0$ because $d\alpha$ is
integrable.
\end{pf}

\begin{theorem}
\label{thm:YMH} Let
$\CCC=((C,\bx_{\crit},\bx_{\gen}),(P,A,G),\phi,\{\TTT_y\},\{\TTT_x\})$
be a $c$-STHM.
We have $\bx_{\gen,A}\subset\bx_{\gen}$ and, for any
$\tau\in T(P,\bx_{\gen,A})$,
$$\YMH_c(\CCC)=
2\int (\Phi^{\tau,A})^*[\omega(A)] -4\pi\imag c\deg\bP^{\tau,A}
+4\pi\sum_{x\in\bx_{\gen}}\Res(A,x,\tau_x)(\mu(\phi(x))-c).$$
\end{theorem}
\begin{pf}
The inclusion $\bx_{\gen,A}\subset\bx_{\gen}$ follows from the
definition of $c$-STHM's. So we only have to prove the formula for
$\YMH_c(\CCC)$. For that we can compute separatedly the integral
on each principal component $C_i\subset C^p$ (the computation for
the bubbles is as in standard Gromov--Witten theory, except the
case in which the connection is flat but has poles in two points
of the bubble, which follows from Lemma
\ref{lemma:integral-omega-A}). Using Lemma
\ref{lemma:integral-omega-A} and the Chern--Weil formula in Lemma
\ref{lemma:ChernWeil} we have (here all integrals are over $C_i$,
$\by_{\gen}$ denotes the marked points and nodes of $C_i$ on which
$A$ has a pole with noncritical residue, $v$ denotes the volume
form, and all norms are $L^2$ over $C_i$):
\begin{align*}
0 &= \|\iota_{v}F_A+\mu(\phi)-c\|^2 \\
&= \|F_A\|^2 +\|\mu(\phi)-c\|^2
-2\int F_A c + 2\int F_A\mu(\phi) \\
&= \|F_A\|^2 +\|\mu(\phi)-c\|^2 +4\pi\imag c\left(\deg
\bP^{\tau,A}
-\imag\sum_{x\in\by_{\gen}}\Res(A,x,\tau_x)\right) \\
&+\|d_A\phi\|^2-4\pi\sum_{x\in\by_{\gen}}\Res(A,x,\tau_x)\mu(\phi(x))
-2\int(\Phi^{\tau,A})^*[\omega(A)].
\end{align*}
Rearranging the terms we obtain the desired formula.
\end{pf}

\section{Bounding the number of bubbles in terms of $\YMH_c$}
\label{s:bounding-bubbles}

\begin{theorem}
\label{thm:bounding-number-bubbles} For any $K>0$ and any $g,n$
satisfying $2g+n\geq 3$ there exists some $N$ with the following
property. Let $\CCC=((C,\bx),(P,A,G),\phi,\{\TTT_y\},\{\TTT_x\})$ be a
$c$-STHM of genus $g$ and $n$ marked points. Suppose that
$\YMH_c(\CCC)\leq K$. Then the number of bubbles in $C$ is less
that $N$.
\end{theorem}

The proof of Theorem \ref{thm:bounding-number-bubbles}
is given in Section \ref{ss:proof-thm:bounding-number-bubbles}.
Note that the corresponding result in Gromov--Witten theory is
an immediate consequence of the existence of a lower bound on
the energy of nontrivial bubbles. In our situation such a lower
bound does not exist, and this is why the theorem is not so obvious.
The idea of the proof is the observation that bubbles with very
little energy come in fact from gradient segments, and that a chain
of consecutive bubbles with little energy gives rise to a chain
of gradient segments, with at least as many components as the
chain of bubbles. The number of components in a chain of gradient
segments is at most equal to the number of connected components
of the fixed point set of $X$. Hence, there is an upper bound
for the length of a chain of consecutive bubbles with little energy.

\subsection{Connecting bubbles and tree bubbles}
\label{ss:connecting-tree-bubbles} Take some $c$-STHM $\CCC$, let
$(C,\bx)$ be the corresponding marked curve, let $P$ be the
principal bundle, and let $A$ be the meromorphic connection. Let
$\Gamma:=\Gamma(C,\bx)$ be the graph whose set of vertices is the
set of irreducible components of $C$ plus the set of marked
points, and whose set of edges is the following:
\begin{itemize}
\item for any pair of vertices $v',v''$ corresponding to components $C',C''$,
there are as many edges connecting $v'$ to $v''$ as there are nodes in
$C$ at which $C'$ and $C''$ meet
(in particular, $\Gamma$ has a loop for each node whose
two branches belong to the same component);
\item if $v_x$ is a vertex corresponding to a marked point $x\in\bx$,
then there is an edge connecting $v_x$ to the vertex corresponding to
the component of $C$ in which $x$ is contained.
\end{itemize}
We divide the set of vertices of $\Gamma$ in {\bf bubble vertices},
{\bf principal vertices} and {\bf marked point vertices}.
A bubble vertex $v$ of $\Gamma$ is said to be
{\bf exterior} if there is a unique edge having $v$ as one of its
extremes (such edge is not allowed to be a loop).
A bubble vertex $v$ of $\Gamma$ is called a {\bf tree vertex}
if there exists a saturated subgraph $T\subset\Gamma$ such that:
$v$ is a vertex of $T$, $T$ is a tree, all vertices of $T$ are
bubble vertices, and there is a unique edge
conneting a vertex of $T$ to a vertex of $\Gamma$.
For example, an exterior vertex is a tree vertex. A bubble vertex
which is not a tree vertex is called a {\bf connecting vertex}.
We define the {\bf depth} of a tree vertex $v$ to be the minimal $p$
for which there is a sequence of tree vertices $v=v_1,\dots,v_p$
such that $v_p$ is exterior and each $v_j$ is connected by an edge
to $v_{j+1}$.

A {\bf chain of connecting vertices} is a connected
saturated subgraph $R\subset\Gamma$ all of whose vertices are
connecting bubbles. The name is motivated by the fact that such a
graph is necessarily homeomorphic to a segment, as the reader
can easily check.

We say that a bubble $C'\subset C$ is a {\bf tree bubble} if the
corresponding vertex in $\Gamma$ is a tree vertex. Similarly, we
define the {\bf connecting bubbles}.
Note that a tree bubble cannot contain any marked point.

\begin{lemma}
\label{lemma:tree-bubbles-flat}
If $C'\subset C$ is a tree bubble, then the restriction of $A$
to $C'$ has trivial holonomy around any node contained in $C'$.
Hence $(P,A)$ extends smoothly to $C'$.
\end{lemma}
\begin{pf}
For the first assertion,
apply induction on the depth of tree vertices, using the fact that
$F_A$ restricts to $0$ on any bubble and the Chern--Weil formula
in Lemma \ref{lemma:ChernWeil}. The second statement follows from
Corollary \ref{cor:trivial-hol-extension}.
\end{pf}

\begin{lemma}
\label{lemma:chain-connecting-bubbles}
Let $C_1,\dots,C_u$ be a sequence of connecting bubbles corresponding
to a chain of connecting vertices $R\subset\Gamma$, labelled in such
a way that each $C_j$ shares a node with $C_{j+1}$.
There exists some $\lambda\in\imag\RR$ with the following property.
Let the exceptional points of $C_j$ be $y_+,y_-,z_1,\dots,z_s$,
where $y_{\pm}$ is the node which $C_j$ has in common
with $C_{j\pm 1}$, except that if $j=1$ then $y_-$ is either
a marked point or a node which $C_1$ shares with a principal
component or another connecting bubble,
or if $j=s$ then $y_+$ may be either a marked point or a node
(connecting $C_s$ to a principal component or a connecting bubble).
Then: each $z_t$ is a node which $C_j$ shares with a tree bubble,
the holonomy of $A$ around any $z_k$ is trivial, and the holonomy
of $A$ around $y_{\pm}$ is $e^{\pm 2\pi\lambda}$.
\end{lemma}
\begin{pf}
Let $R\subset\Gamma$ be a maximal chain of connecting vertices.
Any vertex of $\Gamma$ which is connected
by an edge to an interior vertex of $R$ is necessarily a tree bubble.
Furthermore, each of the two vertices at the extremes
of $R$ are connected either to a marked point vertex or
to a principal vertex, and all other vertices connected to them
and not contained in $R$ are tree vertices.
This explains why the exceptional points of $C_j$ can be labelled
as $y_+,y_-,z_1,\dots,z_s$ ($s$ depends on $j$)
and have the properties claimed in the statement
of the lemma.
By Lemma \ref{lemma:tree-bubbles-flat} the holonomy of $A$ around $z_t$
is trivial. So $A$ can have nontrivial holonomy only around $y_{\pm}$.
Finally, taking any $\lambda$ such that for $C_1$ the holonomy around
$y_-$ is $e^{-2\pi\lambda}$, we deduce that the holonomy around the
point $y_{\pm}\in C_j$ is $e^{\pm 2\pi\lambda}$ (use induction on $j$).
\end{pf}

\subsection{Twisted bubbles}

Let $\by\subset S^2$ be a finite subset. We call a {\bf
twisted bubble} over $(S^2,\by)$ a triple $(P,A,\phi)$ consisting
of a principal bundle $P\to S^2\setminus\{\by\}$, a flat
meromorphic connection $A$ on $P$ and a section $\phi$ of
$P\times_{S^1}X$ satisfying $\ov{\partial}_A\phi=0$. We say
that a twisted bubble $(P,A,\phi)$ is {\bf trivial} if the covariant
derivative $d_A\phi$ is identically zero.

Let $\epsilon>0$ be a small number. Let
$\Lambda^{\epsilon}_{\crit}\subset\imag\RR$ be the set of residues
at distance $<\epsilon$ from the set or critical residues
$\Lambda_{\crit}$. Suppose that $\epsilon$ is small enough so that
for any $\lambda\in\Lambda_{\crit}^{\epsilon}$ there is a unique
critical residue $\crit(\lambda)\in\Lambda_{\crit}$ lying in the
same connected component as $\lambda$. In this case, define
$$\Lambda_{\crit}^+=\{\lambda\in\Lambda_{\crit}^{\epsilon}
\mid -\imag\lambda>-\imag\crit(\lambda)\}, \qquad
\Lambda_{\crit}^-=\{\lambda\in\Lambda_{\crit}^{\epsilon} \mid
-\imag\lambda<-\imag\crit(\lambda)\}.$$ We then have a partition
$\Lambda_{\crit}^{\epsilon}=\Lambda_{\crit}\cup
\Lambda_{\crit}^+\cup \Lambda_{\crit}^-$. Define the set of
generic residues to be
$$\Lambda_{\gen}:=\imag\RR\setminus\Lambda_{\crit}^{\epsilon}.$$

\begin{theorem}
\label{thm:poca-energia-tururut} Let $\by=\{y_+,y_-\}$ consist of
two points. There exists some $\epsilon>0$ with the following
property.  Suppose that $(P,A,\phi)$ is a twisted bubble on
$(S^2,\by)$ satisfying $\|d_A\phi\|_{L^2}<\epsilon$. Take some
local trivialization of $P$ around $y_+$ and let $\lambda$ be the
corresponding residue of $A$ at $y_+$. Then, if $\epsilon$ is
small enough, we have
\begin{enumerate}
\item if $\lambda\in\Lambda_{\crit}$ then $(P,A,\phi)$ is trivial;
\item otherwise, by Theorem \ref{thm:extensio-pol} the following
limits exist
$$\phi(y_+):=\lim_{z\to y_+}\phi(z)\in F,\qquad
\phi(y_-):=\lim_{z\to y_-}\phi(z)\in F;$$ let $F_{\pm}$ be the
connected component of $F$ containing $\phi(y_{\pm})$. Then we
have
\begin{enumerate}
\item if $\lambda\in\Lambda_{\crit}^+$ then $H(F_+)\leq H(F_-),$ with
equality only if $(P,A,\phi)$ is trivial, \item
if $\lambda\in\Lambda_{\crit}^-$ then $H(F_+)\geq H(F_-),$ with
equality only if $(P,A,\phi)$ is trivial, \item if
$\lambda\in\Lambda_{\gen}$ then $(P,A,\phi)$ is trivial.
\end{enumerate}
\end{enumerate}
\end{theorem}
\begin{pf}
Suppose first that $\lambda\in\Lambda_{\crit}$, and write
$\lambda=\imag p/q$ for relatively prime integers $p,q$. Take a
covering $\pi:S^2\to S^2$ of degree $q$ ramified at $y_+$ and
$y_-$, in such a way that $\pi^*A$ has trivial holonomy around
$y_+$ and $y_-$. It follows from Corollary
\ref{cor:trivial-hol-extension} that both $P$ and $A$ extend to
give a bundle $P'\to S^2$ with a smooth connection $A'$. Then
$A'$ is flat, so trivializing any fibre of $P'$ we obtain a
trivialization of the whole bundle $P'\simeq S^2\times S^1$, with
respect to which the connection $A'$ is trivial. This induces in
particular a trivialization of $P$, so we can view the section
$\phi$ as a map $\phi:S^2\setminus\by\to X$, and as such it is
$I$-holomorphic and has bounded energy. It follows from Gromov's
removal of singularities theorem that $\phi$ extends to give a map
$\Phi:S^2\to X$. On the other hand, $\|d\Phi\|_{L^2}<\epsilon$,
and if $\epsilon$ is small enough this implies that $\Phi$ is
constant (this is a standard result in Gromov--Witten theory).
Hence $d_A\phi=0$, so $(P,A,\phi)$ is trivial.

For the remaining cases,
fix a conformal isomorphism between $S^2\setminus\by$ and the
cylinder $\RR\times S^1$, and take in the latter the standard coordinates
$(t,\theta)$, in such a way that as $t$ goes to $\pm\infty$ we approach
the point $y_{\pm}$.

Consider first the case $\lambda\in\Lambda_{\crit}^+$.
Suppose that $\epsilon$ is very small, and
let $(P,A,\phi)$ be a bubble satisfying:
$\|d_{A_u}\phi_u\|_{L^2}< \epsilon$,
$d_A=d+\lambda d\theta$ in some trivialisation of $P$,
and $-\imag(\lambda-\lambda_{\crit})\in (0,\epsilon)$.
Let $\alpha:=\lambda d\theta$.
In the rest of the proof we are going to use the notation introduced
in Section \ref{s:exp-decay}.
Using (2) in Theorem \ref{thm:ineq-perturb},
we can assume that both $\|d_{\alpha}\phi\|_{L^{\infty}}$
and $\|\alpha-\lambda_{\crit}d\theta\|_{L^{\infty}}$
are less than the $\epsilon$ in Theorem \ref{thm:small-energy}.
By Theorem \ref{thm:small-energy} there is some
$\psi:\RR\to X^{\lambda_{\crit}}$
and some $\phi_0:\RR\times S^1\to TX$ such that
$\phi=\exp_\psi\phi_0$.
Now, Theorem \ref{thm:psi-gradient} together with $d\alpha=0$ implies
that $|\psi'+\imag\lambda I\XXX(\psi)|<Ke^{\sigma(|t|-N)}$ for every
$N$. Making $N\to\infty$ we deduce that
$\psi'=-\imag\lambda I\XXX(\psi)$, so $\psi$ follows a downward
gradient line of $H$. Consequently, either
$H(F_+)<H(F_-)$ or $H(F_+)=H(F_-)$, and in the latter case $\psi$ is
constantly.
Now, using (\ref{eq:estimate-phi0-flat}) in Theorem \ref{thm:small-energy}
and taking $N$ bigger and bigger as above we deduce that $\phi_0=0$.
It follows that, if $H(F_+)=H(F_-)$, then $\phi$ is constant.

The case $\lambda\in\Lambda_{\crit}^-$ is proved in the same way.
Finally, if $\lambda\in \Lambda_{\gen}$ and $\epsilon$ is small enough,
then we can use (\ref{eq:diametre-energia}) in Theorem
\ref{thm:exp-decay-noncritical} to deduce that the image of $\phi$ is
contained in a small and $S^1$-invariant
ball $B$ centered at $\phi(y_+)\in F$. Similarly, for any
$\tau\in T(P,\by)$ we have $\Phi^{\tau,A}(S^2)\subset P\times_{S^1}B$.
The equivariant cohomology of $B$ comes entirely from
the classifying space. Picking $\tau$ such that $\deg P(\tau)=0$
we deduce that $(\Phi^{\tau,A})^*[\omega(A)]=0$.
This choice of $\tau$ ensures that $\Res(A,y_+,\tau)+\Res(A,y_-,\tau)=0$.
Also, $\phi(y_-)$ belongs to the same connected component of $F$ as
$\phi(y_+)$, so $\mu(\phi(y_+))=\mu(\phi(y_-))$.
Hence, applying Lemma \ref{lemma:integral-omega-A}
to $(P,A,\phi)$ we deduce that
$\|d_A\phi\|_{L^2}=0$, so the bubble is trivial.
\end{pf}

\subsection{Proof of Theorem \ref{thm:bounding-number-bubbles}}
\label{ss:proof-thm:bounding-number-bubbles} Let
$\CCC=((C,\bx),(P,A,G),\phi,\{\TTT_y\},\{\TTT_x\})$ be a $c$-STHM of genus
$g$ and $n$ marked points satisfying $\YMH_c(\CCC)\leq K$. A
bubble $C'\subset C$ is said to be {\bf unstable} if it contains
less than $3$ exceptional points. Let $\epsilon$ be as in Theorem
\ref{thm:poca-energia-tururut}. (Note that a particular
consequence is that if $\psi:S^2\to X$ is a $I$-holomorphic map
with energy less than $\epsilon$, then $\psi$ is trivial --- this
corresponds to the case of trivial holonomy around $y_{\pm}$ in
the lemma.)

Let $\Gamma$ be the graph associated to $(C,\bx)$. All subgraphs
of $\Gamma$ which we shall mention will be saturated. The subgraph
of $\Gamma$ consisting of connecting vertices can be written as
the disjoint union of the set of maximal chains of connecting
vertices $R_1,\dots,R_l$. Furthermore, $l$ is at most equal to the
number of marked points $\bx$ plus the number of nodes of $C$
(because the bubbles corresponding to vertices in $R_j$ can be
identified with the connecting bubbles which are contracted by the
stabilisation map either to a given marked point or to a node).
Hence, $l\leq 3g-3+n$.

Each tree vertex of $\Gamma$ belongs to a maximal tree $T\subset\Gamma$,
all of whose vertices are tree vertices.
Hence, the subgraph of $\Gamma$ consisting of
tree vertices is the union of trees $T_1,\dots,T_r$.
Each such tree $T_j$ has a distinguished vertex, which we call the {\bf root},
which is connected by an edge to a vertex of $\Gamma\setminus T_j$.
If $T$ is a tree, we say that a vertex of $T$ is stable (resp. unstable)
if its degree is $\geq 3$ (resp. $\leq 2$). Let $|T|$ denote the
number of vertices of $T$. Then the number of unstable vertices
of $T$ is at least $(|T|+2)/3$. This follows from estimating the number
of stable vertices by counting edges: $T$ has $|T|-1$ edges and
each stable vertex contributes at least $3/2$ to the total number of
edges. Now, an unstable vertex of a tree $T_j$ corresponds to an
unstable bubble of $C$, unless the vertex is the root of $T_j$ and
it has degree $2$ (hence its degree as a vertex of $\Gamma$ is $3$).
In this case, the root can not
be the unique unstable vertex of $T_j$, so there is at least one unstable
bubble for each tree $T_j$. Since each such bubble
is nontrivial, it contributes at least $\epsilon$ to the total energy
$\YMH_c(\CCC)$. Hence, the number of trees can be bounded
as $r\leq \epsilon^{-1}K$.
On the other hand, the previous arguments tell us that for each
$T_j$ there are at least $(|T_j|-1)/3$ unstable bubbles
(we are substracting here a unit in case the root does not correspond
to an unstable bubble). Hence we have
$\sum_{j=1}^r(|T_j|-1)/3\leq\epsilon^{-1}K$, which combined with
the bound on $r$ yields $\sum_{j=1}^r|T_j|\leq \epsilon^{-1}4K$.
Hence, the number of tree bubbles in $C$ is bounded.

It only remains to prove that the length of any chain of
connecting vertices of $\Gamma$ is bounded by a constant
independent of $\CCC$. Let $R\subset\Gamma$ be any such chain. A
vertex of $R$ will be said to be stable if its degree {\it as a
vertex of $\Gamma$} is $\geq 3$, and unstable otherwise. Each
interior vertex of $R$ which is stable is connected by an edge to
at least one of the trees $T_j$. Since each tree $T_j$ is
connected by an edge to a unique vertex in $\Gamma\setminus T_j$,
it follows that the number of interior stable vertices in $R$ is
at most $\epsilon^{-1}K$. So if we prove that each sequence of
consecutive unstable vertices of $R$ is bounded, it will follow
that the number of vertices of $R$ is bounded. Now, this follows
from combining Lemma \ref{lemma:chain-connecting-bubbles} with
Theorem \ref{thm:poca-energia-tururut}, and the fact that the
number of connected components of $X$ is finite.

\section{Main theorem: compactness}
\label{s:main-thm}

The following is the main result of the paper.

\begin{theorem}
\label{thm:main-compactness} Let $g$ and $n$ be nonnegative
integers satisfying $2g+n\geq 3$. Let $K>0$ be any number, and let
$c\in\imag\RR$. Let $\{\CCC_u\}$ be a sequence of $c$-stable
twisted holomorphic maps of genus $g$ and with $n$ marked points,
satisfying $\YMH_c(\CCC_u)\leq K$ for each $u$. Then there is a
subsequence $\{[\CCC_{u_j}]\}$ converging to the isomorphism class
of another $c$-stable twisted holomorphic map $\CCC$ .
Furthermore, we have
\begin{equation}
\label{eq:YMH-convergeix}
\lim_{j\to\infty}\YMH_c(\CCC_{u_j})=\YMH_c(\CCC).
\end{equation}
\end{theorem}

The rest of this section is devoted to the proof of the theorem.
Many ideas involved in the proof are the same that appear in the
compactness theorem for stable maps in Gromov--Witten theory,
which we will assume that the reader is familiar with (see for
example \cite{T,RT1,IS,FO}). First of all, note that if a
subsequence satisfies $\CCC_{u_j}\to\CCC$ then, by Theorem
\ref{thm:YMH}, (\ref{eq:YMH-convergeix}) holds automatically.

\subsection{Getting the first limiting curve}
\label{ss:limit-curve}
 Let $(C_u,\bx_u)$ be the nodal marked curve
underlying $\CCC_u$, so that $\bx_u$ is the union of the critical
points $\bx_{\crit,u}$ and the generic ones $\bx_{\gen,u}$. By
Theorem \ref{thm:bounding-number-bubbles} the number of bubbles in
each $C_u$ is uniformly bounded. It follows that we can assume
(passing to a subsequence) that all curves $C_u$ have the same
topological type. For any $u$, let $\bx_u^0\subset C_u$ be a list
of points such that, setting $\bx_u':=\bx_u\cup\bx_u^0$, the
marked curve $(C_u,\bx_u')$ is stable, and suppose that each
$\bx_u^0$ has as few elements as possible. Taking a subsequence,
we can assume that $(C_u,\bx_u')$ converges to a stable curve
$(C',\bx')$. The new set of marked points $\bx'$ contains the
limit $\bx$ of the sequences $\bx_u\subset\bx_u'$. We call the
points in $\bx$ {\bf original marked points}. It is necessary to
make this distinction because the connection may have poles in the
original marked points, whereas in the other marked points it will
always be smooth. We remark also, to avoid confusion, that in the
course of the proof the lists of marked points $\bx_u'$ and $\bx'$
will increase, and that the curve $C'$ will change from time to
time (the changes will be addition of rational components).

In the names of the following three sections (where we describe
how bubbles are to be added to the limit curve $C'$) we use the
terminology of Section \ref{ss:connecting-tree-bubbles}.

\subsection{Adding tree bubbles, first part}
\label{ss:simple-bubbling} We first consider bubbling off away
from nodes and marked points. Following the approach in \cite{FO},
we force the appearence of bubbles in the domains $(C_u,\bx_u')$
by adding new marked points near points where $|d_{A_u}\phi_u|$
blows up. (In order to have a choice of Riemannian metric on each
stable curve
--- which we need to define $|d_{A_u}\phi_u|$ --- we chose an
element in $\Met_{g,n'}$, where $n'$ is the number of points in
$\bx'$, see Section \ref{ss:universal-metrics}. Later on we will
use metrics on stable curves with more than $n'$ marked points,
and it will be implicitly assumed that we have chosen an element
of the corresponding space $\Met$.) Suppose that $K\subset C'$ is
a compact set which does not contain any original marked point nor
any node. Since $(C_u,\bx_u')$ converges to $(C',\bx')$, we can
assume, provided $u$ is big enough, that there is a canonical
inclusion of $\iota_K:K\to C_u$, whose image we call $K_u$ (see
Section \ref{ss:convergencia-corbes}). Assume that
$\sup_{K_u}|d_{A_u}\phi_u|$ is not bounded as $u$ goes to
infinity. Then we pick for each $u$ a point $x_u^1\in K_u$ where
$|d_{A_u}\phi_u|$ attains its supremum, which we denote by $s_u$,
and another point $x_u^2$ at distance $s_u^{-1}$ from $x_u^1$. We
add $x_u^1$ and $x_u^2$ to the list $x_u'$ and obtain a new marked
nodal curve still denoted by $(C_u,\bx_u')$. Passing to a
subsequence, we can assume that $\{(C_u,\bx_u')\}$ converges to
another limit curve $(C',\bx')$. Then we repeat the process: we
take a compact set $K\subset C'$ not containing original marked
points nor smooth points, and so on. The process stops when, for
any such $K$, the sequence $\sup_{K_u}|d_{A_u}\phi_u|$ is bounded.
That the process stops after adding a finite number of points is
proved exactly as in Gromov--Witten theory (see for example
Proposition 11.3 in \cite{FO}): namely, each time we add two
points, there appears a new bubble, and each bubble which appears
contributes more than a certain amount $\epsilon>0$ to the energy
$\YMH_c(\CCC_u)$, so there cannot be infinitely many. The key
point here is that each $A_u$ is smooth on $K_u$ and has bounded
curvature. Hence, when we zoom up the connections become more and
more flat (this is because we are in real dimension $2$), and in
the limit they become trivial. So the bubbles which appear in this
situation are actually holomorphic maps in the usual sense.

\subsection{Adding tree bubbles, second part}
\label{ss:the-necks} Suppose that the limit curve $(C',\bx')$ has
$k$ nodes, so that the curves of the form $(C_u,\bx_u')$ can be
identified with $C'(I',\delta_{u,1},\dots,\delta_{u,k},\bx')$. Now
take some $j$ such that $\delta_u:=\delta_{u,j}$ is nonzero for
big enough $u$. We say in this case that the node is {\bf new},
i.e., it does not appear in the curves $C_u$ if $u$ is big. Let
$N_u:=N_j(\delta_{u,j})$ be the neck defined in Section
\ref{ss:convergencia-corbes}. We fix an isomorphism
$$N_u\simeq [\ln |\delta_u|-\ln\epsilon, \ln\epsilon]\times
S^1=:[p_u,q_u]\times S^1$$ and we denote as usual by $t,\theta$
the standard cylindrical coordinates. Also, we take on $N_u$ the
standard cylindrical metric. For any $\Delta>0$, let
$N_u(\Delta):=[p_u+\Delta,q_u-\Delta]$. Suppose that
\begin{equation} \label{ss:bombolles-al-cilindre}
 \limsup_{\Delta\to\infty}
\left(\limsup_{u\to\infty}\sup_{N_u(\Delta)}|d_{\alpha_u}\phi_u|
\right)=\infty.
\end{equation}
In this case, we proceed as in \S\ref{ss:simple-bubbling}, picking
pairs of sequences of points $x_u^1$ and $x_u^2$ at larger and
larger distance from the boundary of $N_u$, such that
$s_u:=|d_{\alpha_u}\phi_u(x_u^1)|$ goes to infinity and $x_u^2$ is
at distance $s_u^{-1}$ from $x_u^1$. We then add $x_u^1$ and
$x_u^2$ to $\bx_u'$, go to a subsequence so that there is a new
limit $(C_u,\bx_u')\to (C,\bx')$, and repeat the process. We
repeat this, going again to the beginning of \S\ref{ss:the-necks},
as many times as possible, and at some point we must stop.

In case the $j$-th node was not new, that is, $\delta_u=0$ for big
enough $u$, we consider the normalization of each $C_u$ near the
node, and apply the same technique as in the previous case in a
neighborhood of each of the two preimages of the node. The only
difference in this case is that instead of having a finite
cylinder $N_u$ we will have a semiinfinite cylinder $N_u\simeq
[\ln\epsilon,\infty)\times S^1$, and we define
$N_u(\Delta):=[\ln\epsilon+\Delta,\infty)\times S^1$.

Finally, we apply the same technique around each original original
marked point, modelling again a punctured neighborhood of it with
a semiinfinite cylinder $[0,\infty)\times S^1$.

\subsection{Connecting bubbles appear}
\label{ss:creating-tree-bubbles} Take again some node in $C'$ and
assume that it is new. Define the cylinders $N_u$ and
$N_u(\Delta)$ and in \S\ref{ss:the-necks}. Since we repeated the
process in \S\ref{ss:the-necks} as many times as possible, we must
have now
$$ \limsup_{\Delta\to\infty}
\left(\limsup_{u\to\infty}\sup_{N_u(\Delta)}|d_{\alpha_u}\phi_u|
\right)<\infty.$$
This means that for some $\Delta>0$ we have
$\limsup_{u\to\infty}\sup_{N_u(\Delta)}|d_{\alpha_u}\phi_u|<\infty.$
We now replace $N_u$ by $N_u(\Delta)$ and denote the extremes of
the new cylinder $N_u$ again by $p_u$ and $q_u$.

The neck $N_u\subset C_u$ can either belong to a principal or a
bubble component. Suppose we are in the first case. Then the
vortex equation $\iota_{d\vol(\nu_u)}F_A+\mu(\phi)=c$ is satisfied
($\nu_u$ denotes the restriction of the metric
on the componenet $C_u$ to which $N_u$ belongs).
Let us write $d\vol(\nu_u)=f_u dt\wedge d\theta$.
We have exponential bounds on each derivative of $f$: for any
$l>0$,
\begin{equation}
\label{eq:volume-form-bounded} |\nabla^l f_u(z)|\leq K_l
e^{-d(z,\partial N_u)},
\end{equation}
where $K_l$ is independent of $u$. This follows from the fact that
$\nu_u$ is the restriction to $C_u$ of a smooth metric on the
universal curve over a moduli space of stable curves, which is a
compact orbifold; hence all its derivatives are bounded. On the
other hand, the vortex equation takes the following form:
\begin{equation}
\label{eq:vortex-cilindric} d\alpha_u=f_u(c-\mu(\phi_u)).
\end{equation}

Let $\epsilon>0$ be smaller than the $\epsilon$'s in Theorems
\ref{thm:exp-decay-noncritical} and  \ref{thm:small-energy}. For
any $u$, let $d_u:[p_u,q_u]\to\RR$ be the function defined as
$$d_u(t):=\sup_{\{t\}\times S^1}|d_{\alpha_u}\phi_u|.$$
We define an {\bf $\epsilon$-bubbling list} to be a list of
sequences $(\{b_u^1\},\dots,\{b_u^r\})$ satisfying:
\begin{enumerate}
\item each $b_u^j$ belongs to $N_u$, \item denoting $b_u^0:=p_u$
and $b_u^{r+1}:=q_u$ we have, for each $j$ between $0$ and $r$,
$$\lim_{u\to\infty} b_u^{j+1}-b_u^j=\infty.$$
\item for each $j$ we have $\liminf_{u\to\infty}
d_u(b_u^j)\geq\epsilon$.
\end{enumerate}

\begin{lemma}
The number of sequences in an $\epsilon$-bubbling list is bounded
in terms of $K:=\sup_u\YMH_c(\CCC_u)$.
\end{lemma}
\begin{pf}
For any $b\in [p_u+1/2,q_u-1/2]$, let $C(b)$ be the cylinder
$[b-1/2,b+1/2]\times S^1$. The lemma follows from this claim:
there exists some $\eta>0$ (independent of $\{\CCC_u\}$) such
that, whenever $d_u(b)\geq\epsilon$, we have
$\|d_{\alpha_u}\phi_u\|_{L^2(C(b))}\geq\eta$. Indeed, then the
length of a $\epsilon$-bubbling list is at most $\eta^{-1}K$. To
prove the claim, fix some $K_0>0$. If
$\|d_{\alpha_u}\phi_u\|_{L^2(C(b))}\leq K_0$, then putting
momentarily the connections in Coulomb gauge and using elliptic
bootstrapping alternatively with equation
(\ref{eq:vortex-cilindric}) and $\ov{\partial}_{\alpha_u}\phi_u=0$
(taking the bound (\ref{eq:volume-form-bounded}) into account) we
deduce a uniform bound on the $L^2_2$ bound of $d\alpha_u$. Since
$L^2_2\subset L^p_1$ for any $p$, we may apply Corollary
\ref{cor:cota-inferior} to deduce that
$\|d_{\alpha_u}\phi_u\|_{L^2(C(b))}\geq\delta$ for some $\delta>0$
independent of $\{\CCC_u\}$. Finally, we set $\eta$ to be the
minimum of $\delta$ and $K_0$.
\end{pf}

Take an $\epsilon$-bubbling list $(\{b_u^1\},\dots,\{b_u^r\})$ of
maximal length. If the list is empty, then we do nothing, and
begin again the process in \S\ref{ss:creating-tree-bubbles} with
another node. If there is no node near which we can construct a
nonempty $\epsilon$-bubbling list, then we pass to the next step
in \S\ref{ss:construction-limit}. If instead the list is nonempty,
then we define the points $x_u^{i,1}:=(b_u^i,1)\in N_u$ and
$x_u^{i,2}:=(b_u^i+1,1)\in N_u$. We add to $\bx_u'$ these $2r$ new
points and do as always: pass again to a subsequence so that there
is a new limiting curve $(C_u,\bx_u')\to (C',\bx')$ and begin
again with a node of the new curve $C'$. The process has to stop
at some moment because the energy is bounded.

As in \S\ref{ss:the-necks}, we do the same for the nodes of $C'$
which are not new: take the normalization and do what we did in
the cylinders $N_u$ in a neighborhood of each preimage of the
node, which is conformally equivalent to a semiinfinite cylinder.
Finally, we consider neighbourhoods of each original marked point
and do exactly the same.

\subsection{Vanishing old bubbles}

After all this process of adding marked points at the curves $C_u$
we end up with a limiting curve $(C',\bx')$. We denote from now
one $C:=C'$. Suppose that $C$ has $k$ nodes (this $k$ will most of
the times be bigger than the one in \S\ref{ss:the-necks}), so that
for big enough $u$ there is an isomorphism of marked curves
$$\xi_u:C(I',\delta_{u,1},\dots,\delta_{u,k},\bx)\to (C_u,\bx_u).$$

We say that a bubble $C_0\subset C$ is {\bf old} if for each node
$z_j$ contained in $C_0$ the smoothing parameters $\delta_{u,j}$
vanish. This means that the bubble $C_0$ already existed in the
curves $C_u$, so for any $u$ we have an inclusion $C_0\subset
C_u$. We say that an old bubble $C_0$ {\bf vanishes in the limit}
if for big enough $u$ we have $d_{A_u}\phi_u|_{C_0}\neq 0$ and
$$\lim_{u\to\infty}\|d_{A_u}\phi_u\|_{L^2(C_0)}=0.$$

\subsection{Constructing the limiting $c$-STHM}

\label{ss:construction-limit} Take a compact set $K\subset C$
disjoint from the nodes of $C$ and the original marked points.
Since we have the bound
$$\sup_u\sup_{K}|d_{\xi_u^*A_u}\xi_u^*\phi_u|<\infty,$$
standard arguments (for example, Lemma
\ref{lemma:easy-compactness} combined with a patching argument as
in \S4.4.2 in \cite{DK}) imply that there is a subsequence of the
sequence $\xi_u^*(P_u,A_u,\phi_u)$ which, after regauging,
converges to a limiting triple $(P_K,A_K,\phi_K)$. Taking an
exhaustion of the smooth locus of $C\setminus\bx$ by compact sets,
we obtain a limiting triple $(P,A,\phi)$ which satisfies the
vortex equations on principal components, and such that $\phi$ is
holomorphic with respect to $A$. Furthermore, $\|d_A\phi\|_{L^2}$
is finite.

The triple $(P,A,\phi)$ will be part of the limiting $c$-STHM
$\CCC$. Let $\bz$ be the set of nodes of $C$. By construction the
bundle $P$ is defined over $C\setminus(\bx\cup\bz)$, and we now
prove that $A$ is meromorphic. On each principal component the
pair $(A,\phi)$ satisfies the equations
\begin{equation}
\label{eq:equacions-del-limit}
\ov{\partial}_A\phi=0\qquad\text{and}\qquad
\iota_{d\vol(\nu_{[C^{\st},\bx]})}F_A+\mu(\phi)=c
\end{equation}
(here $C^{\st}$ is the stabilization of $C$). The second equation
implies that the curvature of $A$ is bounded, and the first one,
combined with the fact that the energy $\|d_A\phi\|_{L^2}$ is
finite, allows to apply Corollary \ref{cor:existence-limit-orbit}
and deduce that $\mu(\phi)$ extends continuously to $C$. Going
back now to the second equation again, we deduce that $F_A$
extends continuously to $C$, so $A$ is meromorphic. On the other
hand, the restriction of $F_A$ to each of the bubbles of $C$ is
zero. Finally, the stability condition in the bubbles is
satisfied.

Hence, it only remains to construct the limiting gluing data $G$
for $P$ and the collections of chains of gradient segments
$\{\TTT_y\}$, $\{\TTT_x\}$. This will be done in the next three sections.

\subsection{Chains of gradient segments and gluing data at the old
nodes} Let $y,y'$ be the preimages of an old node in the
normalization of $C$. Then $y$ and $y'$ also belong to each of
the normalizations of the curves $C_u$, hence we have chains of
segments of gradient lines $\TTT_{u,y}$ and $\TTT_{u,y'}$. Since
the space of gradients of segment lines is compact, we can assume
that there are limit chains $\TTT_{y}$ and $\TTT_{y'}$. The same
happens with gluing data.

On the other hand, each vanishing old bubble $C_0\subset C$ gives
rise to an infinite gradient segment $\TTT$ in $X$ by Theorem
\ref{thm:poca-energia-tururut}. If $C_0$ has only two exceptional
points (it can't have only one, because then it would be a tree
bubble), then we collapse $C_0$. This identifies two different
nodes in $C$. Suppose that $\TTT_a$ and $\TTT_b$ are the chains of
gradient segments in each of them (taken with the same
orientation). Then the chain of segments in the new node is the
concatenation of $\TTT_a$, $\TTT$ and $\TTT_b$. The gluing data is
the obvious one. We leave the details of this construction to the
reader.

If $C_0$ has more than two exceptional points,
$y_+,y_-,z_1,\dots,z_r$, and the poles of $A_u$ are in $y_{\pm}$,
then we get as before a limit gradient segment $\TTT$. But instead
of collapsing $C_0$ we substitute it by a chain of $r$ trivial
bubbles, each of them containing one of the points $z_1,\dots,z_r$
and two nodes. And in the node shared by the bubbles containing
$z_j$ and $z_{j+1}$ we take the portion of $\TTT$ which lies
between the limits $\lim_{u\to\infty}S^1\cdot \phi(z_j)$ and
$\lim_{u\to\infty}S^1\cdot \phi(z_{j+1})$. Again we leave the
details to the reader.

\subsection{Chains of gradient segments and gluing data at the new nodes}

Take some new node $z\in\bz$ with preimages $y,y'$ in the normalisation
of $C$. For each $u$ let, as always,
$N_u:=N_y(\delta_{u,z})\subset
C(I',\{\delta_{u,z}\},\bx')\simeq C_u$ be the neck
stretching to $z$.  Assume that $N_u=[p_u,q_u]\times S^1$ and
define $N_u(\Delta):=[p_u+\Delta,q_u-\Delta]\times S^1$. There is
some $\Delta_0>0$ such that for any $\Delta\geq\Delta_0$ we have
\begin{equation}
\label{eq:energia-molt-petita}
\limsup_{u\to\infty}\sup_{N_u(\Delta)}|d_{A_u}\phi_u|<\epsilon,
\end{equation}
where $\epsilon>0$ is less than the $\epsilon$'s in Theorems
\ref{thm:exp-decay-noncritical} and  \ref{thm:small-energy}.

Assume that $y$ is in the side of $\{p_u\}\times S^1$ and $y'$
in the side of $\{q_u\}\times S^1$. By Corollary
\ref{cor:existence-limit-orbit} there is a limiting triple
$(P_{y},A_{y},\phi_{y})$ defined over $S_{y}$, and similarly
for $y'$. Let $\OOO\subset X$ be the orbit on which $\phi_{y}$
takes values, and define $\OOO'$ similarly. It follows from our
construction that
\begin{equation}
\label{eq:les-orbites-limit} \lim_{\Delta\to\infty}
\left(\limsup_{u\to\infty}d_{S^1}(\phi(\{p_u+\Delta\}\times
S^1),\OOO)\right)=0,
\end{equation}
and similarly for $\OOO'$.

Pick a trivialization of $P_u$ on $N_u$ with respect to which
$d_{A_u}=d+\alpha_u$ is in balanced temporal gauge with residue
$\lambda_u$ (see Section \ref{ss:balanced-temporal}). We can
assume that $|\lambda_u|\leq 1$ and, passing to a subsequence,
that there is a limit $\lambda_u\to\lambda$. We distinguish two
possibilities.
\begin{enumerate}
\item Suppose first that $\lambda$ is not critical, so that both
$\OOO$ and $\OOO'$ lie in the fixed point set. In this case we
can take the chains $\TTT_{y}$ and $\TTT_{y'}$ to be degenerate,
and in fact we have $\OOO=\OOO'$. To see this it suffices to
prove that
$$\lim_{\Delta\to\infty}
\left(\limsup_{u\to\infty}\diam(\phi_u(N_u(\Delta)))\right)=0.$$
Now, this formula follows from applying Theorem
\ref{thm:exp-decay-noncritical}. Indeed, on the one hand we have
for any $z=(t,\theta)\in N_u(\Delta_0)$ a bound
$$|d_{\alpha}\phi(z)|\leq K e^{-\sigma d(z,\partial
N_u(\Delta_0))}\epsilon.$$ Now, if $\Delta\geq\Delta_0$ and $z\in
N_u(\Delta)$ we have, for the same reason,
$$|d_{\alpha}\phi(z)|\leq K e^{-\sigma(\Delta-\Delta_0)}
e^{-\sigma d(z,\partial N_u(\Delta))}\epsilon.$$ In particular,
$\|d_{\alpha}\phi\|_{L^{\infty}(N_u(\Delta))}$ goes to $0$ as
$\Delta\to\infty$. Combined with (\ref{eq:diametre-energia}), we
deduce that
$$\lim_{u\to\infty}\diam_{S^1}(\phi_u(N_u(\Delta)))=0$$ which,
using (\ref{eq:les-orbites-limit}) and the fact that $\OOO$ is a
fixed point, implies that the actual diameter goes to zero:
$$\lim_{u\to\infty}\diam (\phi_u(N_u(\Delta)))=0$$

\item Now suppose that $\lambda$ is critical. We have to construct
for big enough $u$ a cylinder $M_u:=N_u(\Delta_u)\times S^1$ such
that the conditions given in Section \ref{ss:convergence-STHM}
(subsection {\it Convergence near the nodes to chains of gradient
segments}) are satisfied. Observe first of all that for any
$\Delta>0$ there is some $u(\Delta)$ such that if $u\geq
u(\Delta)$ then
$$M_u^{\Delta}:=([p_u,p_u+\Delta] \cup [q_u-\Delta,q_u])\times S^1$$
is contained in the smooth locus of the limit curve $C$. In
particular, $\phi$ restricts to give a map
$\phi:M_u^{\Delta}\to X.$
On the other hand, for any $\epsilon>0$ there is some
$u(\Delta,\epsilon)\geq u(\Delta)$ such that if $u\geq
u(\Delta,\epsilon)$ then
$$\sup_{z\in M_u^{\Delta}}d_{S^1}(\phi(z),\phi_u(z))<\epsilon$$
(because $\phi_u$ converge modulo gauge to $\phi$). Now take
sequences $\Delta_r\to\infty$ and $\epsilon_r\to 0$, and define
for every $r$ the number $u_r:=u(\Delta_r,\epsilon_r)$. We can
assume without loss of generality that the sequence $\{u_r\}$ is
strictly increasing. Now, take any $u>0$ and define $r$ by the
condition $u_r\leq u<u_{r+1}$. Then define $\Delta_u:=\Delta_r$
and also
$M_u:=N_u(\Delta_u).$
Given $\Delta$ and a big enough $u$, let $M_u^{\Delta,\pm}$ be the
two connected components of $N_u(\Delta)\setminus M_u$. It follows
from the construction that
$$\lim_{\Delta\to\infty}
\left(\limsup_{u\to\infty}
\diam_{S^1}(\phi_u(M_u^{\Delta,\pm}))\right)=0.$$ So it remains to
prove that the sequence of cylinders $(\phi_u,M_u)$ converges with
gauge $\lambda$ to a certain chain $\TTT\in\TTT(X^{\lambda})$.

If $u$ is big enough and $\Delta\geq\Delta_0$, Theorem
\ref{thm:small-energy} gives us a map
$\psi_u(\Delta):T_u(\Delta):=[p_u+\Delta,q_u-\Delta]\to X$. Let
$l_u:=-\imag(\lambda_u-\lambda)$. Formula (\ref{eq:der-psi}) in
Theorem \ref{thm:psi-gradient} (applied first to $T_u(\Delta_0)$
and then to $T_u(\Delta)$) implies that
\begin{equation}
|\psi_u(\Delta)'(t)-l_u
I(\psi(t))\XXX(\psi_u(\Delta)(t))|<Ke^{-\sigma(\Delta-\Delta_0)}
e^{-\sigma d(t,\partial S_u(\Delta)} K
\end{equation}
for some constant $K$. Let $M_u=T_u\times S^1$. Since
$\Delta_u\to\infty$,  the previous inequality allows us to apply
Theorem \ref{thm:convergeix-cap-a-gradient}. Hence, passing to a
subsequence, we can distinguish two possibilities.
\begin{itemize}
\item If $l_u|T_u|\to 0$ then $\diam(\psi_u(T_u))\to 0$, and,
passing again to a subsequence, we may define $\TTT_{y}$ to be
the degenerate chain of gradient segments with constant value the
unique point $x\in X$ such that $\lim d(\psi_u(T_u),x)=0$.
 \item If $\lim l_u|T_u|\neq 0$ then for big enough $u$ we have $l_u\neq 0$, and by
 the theorem
$(l_u^{-1}\psi_u,l_uT_u)$ converges to a chain of gradient
segments, which we denote by $\TTT_{y}$.
\end{itemize}
Inequality (\ref{eq:estimate-phi0}) and formula
(\ref{eq:definicio-phi0}) in Theorem \ref{thm:small-energy} imply
that the renormalized sequence $(l_u^{-1}\phi_u,l_uM_u)$ converges
with gauge $\lambda$ to $\TTT_{y}$.
\end{enumerate}

Finally, to define the gluing data one only needs to pass to a
suitable subsequence, as is clear from the conditions {\it
Convergence of gluing angles} and {\it Convergence of gluing data}
in the definition of convergence of sequences of $c$-STHM (Section
\ref{ss:convergence-STHM}). Indeed, the set of possible gluing
data in a $c$-STHM is compact.

\subsection{Chains of gradient segments in generic marked points}
Let $\{x_u\}$ be a sequence of original marked points, where each
$x_u\in\bx_{\gen,u}$. Let $\lambda_u\in\imag\RR$ satisfy
$e^{2\pi\lambda}=\Hol(A_u,x_u)$, and let $\TTT_u$ be the chain of
gradient segments $\TTT_{x_u}$. Suppose that $\{x_u\}$ converges
to $x\in C'$. This means that there is a compact subset $K\subset
C\setminus(\bx\cup\bz)$ such that for big enough $u$ we have
$x_u\subset \iota_K(K)\subset C_u$. Passing to a subsequence we
can also assume that $\lambda_u\to\lambda\in\imag\RR$. In this
section we show that, passing to a subsequence, there is a well
defined limiting chain of gradient segments $\TTT_x$. The analysis
is very similar to the previous one, so we will be sketchy. We
distinguish two situations.
\begin{enumerate}
\item Suppose that $\lambda$ is not critical. In this case we define
$\TTT_x$ to be the unique degenerate chain which satisfies the
matching condition with $\phi_x$. (Note that for big enough $u$ the
residue $\lambda_u$ is also noncritical, hence $\TTT_u$ is degenerate.)
\item Suppose that $\lambda$ is critical. In this case, passing to
a subsequence, we can assume that one of the following possibilities holds:
\begin{itemize}
\item either $\lambda_u$ is critical for big enough $u$, in which this case
we define $\TTT_x$ to be the limit of the sequence of chains $\{\TTT_u\}$,
\item or $\lambda_u$ is not critical for big enough $u$, and
we follow the same strategy as in case (2) of the previous section
to define $\TTT_x$.
\end{itemize}
\end{enumerate}

\section{Local estimates and cylinders with noncritical residue}
\label{s:exp-decay}

\begin{theorem}
\label{thm:exp-decay} Let $P$ be a principal $S^1$ bundle on the
punctured disk $\DD^*$, and fix an element $\tau\in T(P,0)$.
Assume that $A$ is a connection on $P$ whose curvature is
uniformly bounded, $|F_A|_{L^{\infty}}<\infty$, and such that
$\Res(A,0,\tau)\notin\imag\ZZ$. Let $\phi$ be a section of the
trivial bundle $P\times_{S^1}X$ which satisfies
$\ov{\partial}_A\phi=0$ and $\|d_A\phi\|_{L^2(\DD^*)}<\infty$.
Then, denoting by $(r,\theta)$ the polar coordinates on $\DD^*$,
there are constants $K>0$ and $\nu>0$ such that
$$|d_A\phi(r,\theta)|\leq K r^{\nu}.$$
\end{theorem}

The proof of Theorem \ref{thm:exp-decay} will be given in Section
\ref{proof:thm:exp-decay}.

\begin{corollary}
\label{cor:existence-limit-orbit} Under the hypothesis of Theorem
\ref{thm:exp-decay}, the section $\phi$ extends at the origin to
give a section $\phi_0$ of $Y_0$. In particular, the composition
of $\phi$ with the moment map, $\mu(\phi)$, extends to a
continuous map from the disk $\DD$ to $\imag\RR$. Furthermore we
have:
\begin{enumerate}
\item The section $\phi_0$ is covariantly constant with respect to
the limiting connection $A_0$ on $P_0$, that is,
$d_{A_0}\phi_0=0$. \item Let $\lambda:=\Res(A,0,\tau)$. The
section $\phi_0$ takes values in $X^{\lambda}$. In particular, if
$\lambda$ is not critical, then $\phi_0$ takes values in the fixed
point set $F$, hence $\phi_0$ is constant and the following limit
exists
\begin{equation}
\label{eq:limit-point} \phi(0):=\lim_{z\to 0}\phi(z)\in F.
\end{equation}
\end{enumerate}
\end{corollary}
\begin{pf}
We prove (1). Take a trivialisation of $P$ around
$0$ for which $d_A=d+\alpha+\lambda d\theta$, where $\alpha$ is of
type $C^1$ on the whole disk $\DD$ and is in radial gauge, i.e.,
$\alpha=a d\theta$ for some function $a$ on $\DD$ vanishing at the
origin. The trivialisation of $P$ induces a trivialisation of
$P_0\to S_0$, with respect to which $d_{A_0}=d+\lambda d\theta$.
Using the trivialisation of $P$ we look at the section $\phi$ as a
map $\phi:\DD^*\to X$. Now define for every $0<r<1$ the map
$\phi_r:S^1\to X$ by $\phi_r(\theta):=\phi(r,\theta)$. It follows
from the estimate Theorem \ref{thm:exp-decay} that the limit
$\phi_0:=\lim_{r\to 0}\phi_r$ exists and is of type $C^1$. Define
now the connection $A_r$ on the trivial $S^1$ bundle over the
circle using the $1$-form $\alpha_r(r,\theta):=r
a(r,\theta)d\theta$. It follows also from Theorem
\ref{thm:exp-decay} that $d_{A_0}\phi_0=\lim_{r\to
0}d_{A_r}\phi_r=0$.
Finally, (2) follows from (1), observing that $\phi_0$ takes values in
a unique orbit of the action of $S^1$ on $X$.
\end{pf}

To study the local properties of equation (\ref{eq:holomorf}) we
restrict ourselves to considering trivial principal $S^1$ bundles
$P$ over a (nonnecessarily compact $C$). Then $Y=C\times X$, the
sections of $Y$ can be identified with maps $C\to X$, and the
connections on $Y$ are the same as forms
$\alpha\in\Omega^1(C,\imag\RR)$. Finally, there is a canonical
splitting $TY=TC\oplus TX$ (we omit the pullbacks). With respect
to a given form $\alpha$ the bundle of horizontal tangent vectors
is $T^{\hor}_{\alpha}=\{(u,\imag\alpha(u)\XXX)\mid u\in
TC\}\subset TY.$ Consequently, the covariant derivative of a
section $\phi:C\to X$ is
$d_{\alpha}\phi=d\phi-\imag\alpha\XXX(\phi).$ Finally, if $I$ is
an almost complex structure on $X$, we have
$$\ov{\partial}_{I,\alpha}\phi=\ov{\partial}_I\phi
-\frac{1}{2}(\imag\alpha\XXX(\phi)+\imag(\alpha\circ
I_C)(I\XXX)(\phi)).$$

A pair $(\alpha,\phi)$ consisting of a
$1$-form $\alpha$ on $C$ with values on $\imag\RR$ 
and a map $\phi:C\to X$ 
will be simply called a {\bf pair}; to specify both the curve $C$
and the target manifold $X$ we will write
$$(\alpha,\phi):C\to X.$$
If $I$ is any $S^1$-invariant almost complex structure on $X$, we
will say that the pair $(\alpha,\phi)$ is {\bf $I$-holomorphic} if
the equation $\ov{\partial}_{I,\alpha}\phi=0$ is satisfied.

For any natural number $N$, denote $C_N:=[-N,N]\times S^1$.
We can now state the second main result of this section.

\begin{theorem}
\label{thm:exp-decay-noncritical} For any noncritical
$\lambda\in\imag\RR\setminus\Lambda_{\crit}$ there exist some
$K>0$, $\sigma>0$ and $\epsilon>0$, depending continuously on
$\lambda$, with the following property. Let $(\alpha,\phi):C_N\to
X$ be a $I$-holomorphic pair satisfying
$\|\alpha-\lambda
d\theta\|_{L^{\infty}}<\epsilon$
and
$\|d_{\alpha}\phi\|_{L^{\infty}(C_N)}<\epsilon.$
Then the following inequality holds for any $t,\theta$:
\begin{equation}
\label{eq:energia-L2-puntual} |d_{\alpha}\phi(t,\theta)|\leq K
e^{-\sigma(N-|t|)}\|d_{\alpha}\phi\|_{L^{\infty}(C_N)}.
\end{equation}
In particular this implies
\begin{equation}
\label{eq:diametre-energia} \diam_{S^1}
(\phi(C_N))<K\|d_{\alpha}\phi\|_{L^{\infty}(C_N)}.
\end{equation}
\end{theorem}

The proof of Theorem
\ref{thm:exp-decay-noncritical} will be
given in Section \ref{ss:proof-thm:exp-decay-noncritical}

\subsection{Equivariant charts}

Let $I_0$ (resp. $g_0$) denote the standard complex structure on
$\CC^n$, viewed as a differentiable manifold.

\begin{lemma}
\label{lemma:equi-charts}
There exists a finite set $W=W(X)\subset\ZZ^n$ and
a real number $r=r(X,\omega,I)>0$ with the following property.
For any $\epsilon>0$ and any $x\in X$ there exists
an equivariant open neighbourhood $U$ of $x$,
an action of $S^1$ on $\CC^n$ whose collection of weights belongs to $W(X)$,
an equivariant almost complex structure $I_x$ and metric $g_x$ on $\CC^n$
satisfying $\|I_x-I_0\|_{L^{\infty}}<\epsilon$ and
$\|g_x-g_0\|_{L^{\infty}}<\epsilon$ (both norms taken with respect to
$g_0$), and
an equivariant map $\xi:U\to\CC^n$ which is a complex isometry.
Furthermore, we can assume that $\xi(U)$ is contained
in the ball $B(0,2r)\subset\CC^n$.
\end{lemma}

\begin{pf}
Take $r$ to be the length of the longest orbit in $X$ divided by
$2\pi$. If $x$ is a fixed point the chart is easily constructed
using the exponential map with respect to the invariant metric
$g$. Suppose that $x$ is not a fixed point, and that the length of
the orbit through $x$ is $2\pi\rho$. Let $\Gamma$ be the
stabiliser of $x$, which acts linearly on $T_xX$. Let $L\subset
T_xX$ be the complex subspace generated by $\XXX(x)$ and
$I\XXX(x)$ (recall that $\XXX\in\Gamma(TX)$ is the vector field
generated by the infinitesimal action of $\Lie S^1$). Both
$\XXX(x)$ and $I\XXX(x)$ are fixed by the action of $\Gamma$, so
$L$ is $\Gamma$-invariant. Let $N\subset T_xX$ be the Hermitian
ortogonal of $L$, which is a $\Gamma$-invariant and complex vector
subspace. Take a complex isomorphism $\CC^{n-1}\simeq N$. Consider
the induced action of $\Gamma$ on $\CC^{n-1}$. Fix any linear
extension $\rho$ of this action to $S^1$, which we may assume
(chosing appropiatedly the identification $\CC^{n-1}\simeq N$) to
be diagonal and compatible with the standard metric in
$\CC^{n-1}$. Let $k\in\NN$ be the order of $\Gamma$, consider the
action of $S^1$ on $\CC^*\times\CC^{n-1}$ given by
$\theta\cdot(x,y):=(\theta^kx,\rho(\theta)(y))$. There exists an
$S^1$-invariant neighbourhood $R\subset\CC^*$ of the circle
$S_{\rho}$ of radius $\rho$ centered at the origin in $\CC$ and an
$S^1$-equivariant embedding
$$f:R\times\CC^{n-1}\to X$$
(here we view $R\times\CC^{n-1}$ as a subset of $\CC^*\times\CC^{n-1}$)
which satisfies, for any $t\in\RR$ and $n\in\CC^{n-1}$,
$f(e^{t},u):=\exp_x(ktI\XXX(x)+u)$.
The differential of $f$ at any point in $S^1\times\{0\}$ is complex.
Hence, if $R$ is small enough, there exists an equivariant almost
complex structure $I_x$ on $\CC^n$ satisfying
$\|I_x-I_0\|_{L^\infty}<\epsilon$ and which induces a complex structure
on $R\times\CC^{n-1}$ with respect to which the restriction
of $f$ to a neighbourhood $V$ of $R\times \CC^{n-1}$
is complex. Also, since the restriction of $f$ to
$S_{\rho}\times\{0\}$ is an isometry, we can assume that the pullback metric
$f^*g$ extends to $\CC^n$ satisfying $\|g-g_0\|_{L^\infty}<\epsilon$.
Then we set $U=f(V)$ and $\xi:=(f|_V)^{-1}$.

It is easy to see, using the compactness of $X$ and the rigidity
of representations of compact groups, that the set of
weights of the representations which we construct as $x$ moves along $X$
forms a finite set.
\end{pf}

\subsection{An inequality for the energy on cylinders}

\label{ss:inequality-cylinders}

Let $C:=[-2,3]\times S^1$ with the standard product metric.
Whenever we write any norm ($L^\infty$, $L^2$, etc.) of
either the connection or the section of a pair
defined over $C$, unless we specify some other domain,
we mean the norm over $C$.
Define also the following subsets of $C$:
\begin{equation}
\label{eq:definicio-CIII}
Z:=[-1,2]\times S^1,\quad
\ZI:=[-1,2]\times S^1,\quad
\ZII:=[0,1]\times S^1,\quad
\ZIII:=[1,2]\times S^1.
\end{equation}

\begin{theorem}
\label{thm:ineq-perturb}
For any noncritical residue $\lambda\in\imag\RR\setminus\Lambda_{\crit}$
there exist real numbers $\epsilon=\epsilon(\lambda,X,I)>0$
$\gamma=\gamma(\lambda,X,I)\in (0,1/2)$,
and $K=K(\lambda,X,I)>0$,
depending continuously on $\lambda$, such that if
$(\alpha,\phi):C\to X$ is an $I$-holomorphic pair
satisfying the conditions
$\|\alpha-\lambda d\theta\|_{L^{\infty}(C)}\leq\epsilon$
and
$\|d_{\alpha}\phi\|_{L^2(C)}\leq\epsilon,$
then
\begin{enumerate}
\item the following inequality holds:
$$\|d_{\alpha}\phi\|^2_{L^2(\ZII)}
\leq\gamma\left(\|d_{\alpha}\phi\|^2_{L^2(\ZI)}+
\|d_{\alpha}\phi\|^2_{L^2(\ZIII)}\right);$$
\item if $I$, $\alpha$ and $\phi$ are of class $C^1$, then
$\sup_Z|d_{\alpha}\phi|\leq
K\|d_{\alpha}\phi\|_{L^2(C)}.$
\end{enumerate}

\end{theorem}

To prove Theorem \ref{thm:ineq-perturb} we need an analogous
result for the case $X=\CC^n$ which we now state. Fix a nontrivial
diagonal action of $S^1$ on $\CC^n$, and denote its weights by
$w=(w_1,\dots,w_n)\in\ZZ^n$. Consider on $\CC^n$ the standard
Riemannian metric $g_0$, and denote by $I_0$ the standard complex
structure. Finally, define, for every $\eta\in\RR$,
\begin{equation}
\label{eq:def-gamma}
\gamma(\eta):=\frac{1}{e^{\eta/2}+e^{-\eta/2}}.
\end{equation}

\begin{lemma}
\label{lemma:ineq-energia} There exist real numbers
$\epsilon=\epsilon(w)>0$ and $K=K(w,\lambda)$ with the following
property. Suppose that $I$ is a smooth equivariant almost complex
structure in $\CC^n$ such that $\|I-I_0\|_{L^{\infty}}<\epsilon$,
and that $(\alpha,\phi):C\to\CC^n$ is a $I$-holomorphic pair
satisfying $\|\alpha-\lambda
d\theta\|_{L^{\infty}(C)}\leq\epsilon$ and
$\|d_{\alpha}\phi\|_{L^2(C)}\leq\epsilon.$ Then
\begin{enumerate}
\item the following inequality holds:
\begin{equation}
\label{eq:average-bound-energy-loc} 
\|d_{\alpha}\phi\|^2_{L^2(\ZII)}
\leq\gamma(2l_{\min})\left( \|d_{\alpha}\phi\|^2_{L^2(\ZI)}+
\|d_{\alpha}\phi\|^2_{L^2(\ZIII)}\right),
\end{equation}
where
$l_{\min}:=\min\{|-k+\imag w_j\lambda|\mid
k\in\ZZ,\quad 1\leq j\leq n\};$
\item we can bound:
\begin{equation}
\label{eq:pointwise-bound-energy-loc}
\sup_Z|d_{\alpha}\phi|\leq K \|d_{\alpha}\phi\|_{L^2(C)}.
\end{equation}
\end{enumerate}
\end{lemma}

Lemma \ref{lemma:ineq-energia} will be proved in Section
\ref{proof:lemma:ineq-energia} below. We now prove Theorem
\ref{thm:ineq-perturb}. Take $\epsilon$ smaller than the infimum
of $\epsilon(w)$ as $w$ moves along $W(X)$. Take, for this choice
of $\epsilon$ and for each $x\in X$, an equivariant chart
$\xi_x:U_x\to\CC^n$ as in Lemma \ref{lemma:equi-charts}. Then take
a finite subcovering of $\{U_x\}_{x\in X}$ and denote it by
$\{U_1,\dots,U_r\}$. Denote also by $\xi_j$ and $w_j$ the
corresponding embeddings in $\CC^n$ and weights. Since the sets
$U_j$ are $S^1$-invariant and open, it follows that there is a
constant $\epsilon'>0$ such that for any $K\subset X$ satisfying
$\diam_{S^1}K<\epsilon'$ there is at least one $j$ such that
$K\subset U_j$. Now, if we take $\delta>0$ small enough, it
follows from Lemma \ref{cor:pocaenergia-diampetit} that
$\diam_{S^1}\phi(C)<\epsilon'$. Then the statement of Theorem
\ref{thm:ineq-perturb} follows from Lemma \ref{lemma:ineq-energia}
applied to $\xi_j\circ\phi:C\to\CC^n$. Indeed, the condition that
$\lambda$ is not critical implies that, applying Lemma
\ref{lemma:ineq-energia} to any chart, we have that $l_{\min}>0$,
so $\gamma(2l_{\min})<1/2$. Then $\gamma$ can be defined as the
maximum of the numbers $\gamma(2l_{\min})$ computed for each of
the charts $U_1,\dots,U_r$.

\subsection{From mean inequalities to exponential decay}

The following lemma will allow to obtain exponential decay from
inequalities of the type given by (1) in Lemma
\ref{lemma:ineq-energia}. We will use it several times in this
paper.

\begin{lemma}
\label{lemma:des-recur} Let $\{x_0,\dots,x_N\}$ be a sequence of
nonnegative real numbers satisfying, for some $\gamma\in(0,1/2)$ and
each $k$ between $1$ and $N-1$,
\begin{equation}
\label{eq:condrecur} x_{k}\leq\gamma(x_{k+1}+x_{k-1}).
\end{equation}
Let $\xi:=(1+\sqrt{1-4\gamma^2})/(2\gamma)$. Then
$\xi>1$ and for any $k$ between
$0$ and $N$ we have
\begin{equation}
\label{eq:cotes-exponencials} x_k\leq x_0\xi^{-k}+x_N\xi^{-(N-k)}.
\end{equation}
\end{lemma}
\begin{pf}
That $\xi>1$ is an easy computation.
As for (\ref{eq:cotes-exponencials}), observe that
$\xi^{-1}=\gamma(1+\xi^{-2})$. Consequently, if we define
$y_k:=x_k-(x_0\xi^{-k}+x_N\xi^{-(N-k)})$ then
$y_{k}\leq\gamma(y_{k+1}+y_{k-1})$. Since $\gamma\in (0,1/2)$,
it follows that the sequence $\{y_k\}$ attains its maximum
at $y_0$ or $y_N$. Both of these numbers are $\leq 0$, so
all $y_j$ are nonpositive. This proves
(\ref{eq:cotes-exponencials}).
\end{pf}

\begin{corollary}
\label{cor:des-recur} Suppose that $\{x_0,x_1,x_2,\dots\}$ is an
infinite sequence of nonegative real numbers satisfying
$x_{k}\leq\gamma(x_{k+1}+x_{k-1})$
for some $\gamma\in (0,1/2)$ and any $k\geq 1$.
Suppose that $x_j\to 0$ as $j$ goes to infinity. Define $\xi>1$ as
in Lemma \ref{lemma:des-recur}. Then for any $k$ we have $x_k\leq
x_0\xi^{-k}$.
\end{corollary}
\begin{pf}
Apply Lemma \ref{lemma:des-recur} to the first
$N$ terms and make $N$ go to infinity.
\end{pf}

\begin{lemma}
\label{lemma:recur-perturb} Suppose that
$\{x_{-N},x_{-N+1},\dots,x_N\}$ and
$\{z_{-N},x_{-N+1},\dots,z_N\}$ are sequences of nonnegative real
numbers. Assume that there are positive constants $\gamma\in
(0,1/2)$, $\chi$, $K$ and $\epsilon$ satisfying:
\begin{itemize}
\item for every $j$, $z_j\leq Ke^{-\chi(N-|j|)}$;
\item for every $-N+2\leq j\leq N-2$, if
$z_{j-2}+\dots+z_{j+2}\leq \epsilon(x_{j-2}+\dots+x_{j+2})$
then $x_j\leq\gamma(x_{j-1}+x_{j+1})$.
\end{itemize}
Let $\xi=\xi(\gamma)$ be as in Lemma \ref{lemma:des-recur},
and let $\sigma:=\min\{\chi,\ln\xi\}$. Then we have,
for any $-N+1\leq j\leq N-1$,
\begin{equation}
\label{eq:cota-x-jota}
x_j\leq (\epsilon^{-1}10 e^{2\chi}K+x_{-N+1}+x_{N-1})
e^{-\sigma(N-|j|)}.
\end{equation}
\end{lemma}
\begin{pf}
For any $-N+2\leq j\leq N-2$ we say that $x_j$ is big if $x_j\geq
\epsilon^{-1}(z_{j-2}+\dots+z_{j+2})$, and otherwise we say that
$x_j$ is small. We have $z_{j-2}+\dots+z_{j+2}\leq
5e^{2\chi}Ke^{-\chi(N-|j|)}.$ In particular, if $x_j$ is small
then $x_j\leq\epsilon^{-1}5e^{2\chi}Ke^{-\chi(N-|j|)}$. Take now
any $x_j$. If $x_j$ is small then the previous inequality implies
(\ref{eq:cota-x-jota}). So we assume that $x_j$ is big and we take
the longest sequence of consecutive big elements
$x_{p+1},\dots,x_j,\dots,x_{q-1}$ containing $x_j$. Suppose that
$x_p$ and $x_q$ are small (this need not be the case, since we
could have $p=-N+1$ or $q=N-1$). By Lemma \ref{lemma:des-recur} we
have $x_j\leq x_pe^{-\sigma(j-p)}+x_qe^{-\sigma(q-j)}.$ Now, since
both $x_p$ and $x_q$ are small we can bound (using $\sigma\leq
\chi$)
$$x_j \leq\epsilon^{-1}5e^{2\chi}K(e^{-\sigma(N-|p|+j-p)}
+e^{-\sigma(N-|q|+q-j)}) \leq
\epsilon^{-1}10e^{2\chi}K(e^{-\sigma(-N+|j|)}),$$ because both
$N-|p|+j-p$ and $N-|q|+q-j$ are greater than $N-|j|$. When
$p=-N+1$ or $q=N-1$ we proceed similarly.
\end{pf}

\subsection{Proof of Theorem \ref{thm:exp-decay}}

\label{proof:thm:exp-decay}
By Lemma \ref{lemma:bones-coordenades}
(see Remark \ref{rmk:curvatura-acotada}) we can trivialize $P$ in
such a way that $d_A=d+\lambda d\theta+\alpha'$, where
$\lambda=\Res(A,0,\tau)$ and $\alpha'$ is a continuous $1$-form on
$\DD$. Let $\alpha=\lambda d\theta+\alpha'$. We distinguish two
situations, according to whether the residue is critical or not.

Suppose first that the residue is critical.
In particular, $\imag\lambda$ is rational, so we can write
$\lambda=\imag p/q$, where $p,q$ are relatively prime integers
and $q\neq 1$. Let
$\pi:\DD^*\to\DD^*$ be the map given by $\pi(z):=z^q$. Then
$d_{\pi^*A}=d+q\lambda d\theta+\pi^*\alpha'
=d+pd\theta+\pi^*\alpha'$. Let $g:\DD^*\to S^1$ be defined by
$g(r,\theta):=e^{2\pi\imag p\theta}$. Then
$$d_B:=d_{g^*\pi^*A}=d+\pi^*\alpha'$$
is a continuous connection on $\DD$. Let $I(B)$ be the continuous
almost complex structure on $\DD\times X$ induced by $B$ (see
Section \ref{s:symplectic-fibrations}). By (\ref{eq:holomorf}),
the map $\Phi=(\iota,g\cdot (\phi\circ\pi)):\DD^*\to\DD\times X$
is $I(B)$-holomorphic (here $\iota:\DD^*\to\DD$ denotes the
inclusion), and by (\ref{eq:norma-dPhi}) $d\Phi$ has bounded $L^2$
norm. Hence we can apply the theorem on removal of singularities
for holomorphic curves (as proved for continuous almost complex
structures in Corollary 3.6 of \cite{IS}), deduce that $\Phi$
extends to a $I(B)$-holomorphic map $\Phi:\DD\to\DD\times X$, and
consequently obtain an extension $g\cdot (\phi\circ\pi):\DD\to X$.
Furthermore, this extension is of type $C^1$ (because the complex
structure is continuous). Let now $\zeta:=e^{2\pi\imag/q}$. Then
$(g\cdot (\phi\circ\pi))(\zeta z)=(g\cdot (\phi\circ\pi))(z)$, so
that $d(g\cdot (\phi\circ\pi))(0)=0$. Furthermore, since
$d\pi(0)=0$, $\pi^*\alpha'$ vanishes at $0$. It follows from this
that $d_{g^*\pi^* A}(g\cdot (\phi\circ\pi))=
(d+\pi^*\alpha')(g\cdot (\phi\circ\pi))=0$, from which we deduce
that $|d_{g^*\pi^* A}(g\cdot (\phi\circ\pi))(r,\theta)|= |d_{\pi^*
A}(\phi\circ\pi)(r,\theta)|<Kr$ for some constant $K$. This then
implies that $|d_A\phi(r,\theta)|<K'r^{1/q}$ for some other
constant $K'$, so the claim is proved.

Now consider the case of noncritical residue
$\lambda\notin\Lambda_{\crit}$.
Consider the cylinder $\RR^+\times S^1$ as a conformal model of
$\DD^*$, with coordinates $t\in\RR^+$ and $\theta\in S^1$ and the
standard flat metric $dt^2+d\theta^2$. We have:
\begin{equation}
\label{eq:est-tub} |\alpha'(t,\theta)|<Ke^{-t} \qquad\text{and}
\qquad \|d_{\alpha}\phi\|_{L^2}<\infty.
\end{equation}

We look at $\phi$ as a map from $\DD^*$ to $X$, so that
$(\alpha,\phi)$ is a pair. Define, for any $n\in\NN$,
$Z_n=[n,n+1]\times S^1$ and let
$f_n:=\|d_{\alpha}\phi\|^2_{L^2(Z_n)}.$ We claim that, for big
enough $n$, we have
\begin{equation}
\label{eq:recur} f_{n+1}\leq
\gamma(\lambda,X)(f_n+f_{n+2}),
\end{equation}
where $\gamma(\lambda,X)$ is as in Theorem \ref{thm:ineq-perturb}.
Indeed, thanks to (\ref{eq:est-tub}), if $n$ is big enough then
both the $L^{\infty}$ norm of the restriction of $\alpha'$ to
$[n,n+3]\times S^1$ and the $L^2$ norm of
$d_{\alpha}\phi$ restricted to
$[n,n+3]\times S^1$ are
less than the value of $\epsilon$ given by Theorem
\ref{thm:ineq-perturb}. Hence, we can apply (1) in Theorem
\ref{thm:ineq-perturb} (identifying
$\ZI=Z_n$, $\ZII=Z_{n+1}$ and $\ZIII=Z_{n+2}$)
and deduce (\ref{eq:recur}). Combining
Corollary \ref{cor:des-recur} with (\ref{eq:recur}) we deduce that
$$\|d_{\alpha}\phi\|^2_{L^2(Z_n)}=f_n\leq K \xi^{-n},$$
where $\xi>1$. Now, using (2) in
Theorem \ref{thm:ineq-perturb} we deduce the pointwise bound
$$|d_{\alpha}\phi(t,\theta)|\leq K'\xi^{-t}$$
for some other constant $K'$. This finishes the proof of the case
of noncritical residue.

\subsection{Proof of Theorem \ref{thm:exp-decay-noncritical}}
\label{ss:proof-thm:exp-decay-noncritical} The proof is exactly
like that of the case of noncritical residue of Theorem
\ref{thm:exp-decay}, in Section \ref{proof:thm:exp-decay} above.
Namely, we define $f_n$ to be $\|d_{\alpha}\phi\|_{L^2([n,n+1]\times
S^1)}$ and we use Theorem \ref{thm:ineq-perturb} (provided
$\epsilon$ has been chosen small enough) to deduce that
$f_{n+1}\leq\gamma (f_n+f_{n+2})$ for some $\gamma\in (0,1/2)$
depending on $\lambda$. Then Lemma \ref{lemma:des-recur} gives
$f_n\leq Ke^{-\sigma(N-|n|)}(f_{-N}+f_{N-1})$ for some $K$ and
$\sigma$. To deduce the pointwise bound
(\ref{eq:energia-L2-puntual}) we use
(2) in Theorem
\ref{thm:ineq-perturb}. Finally, (\ref{eq:diametre-energia})
follows easily from (\ref{eq:energia-L2-puntual}) combined with
Lemma \ref{lemma:cov-der-diam-S1}.

\subsection{Some convexity properties of holomorphic maps}
In this subsection we state a result which will be crucial for
proving most of the convexity results in this paper. Part of it
can be seen as a particular case of Lemma \ref{lemma:ineq-energia}
(concretely, the case $I=I_0$ and $\alpha'=0$).

For any $\delta\in (0,1/2)$ we will denote
\begin{equation}
\label{eq:def-Z-delta}
 Z^{\delta}:=[-1+\delta,2-\delta],
 \qquad
 \ZI^{\delta}:=\ZI\cap Z^{\delta},
 \qquad
\ZIII^{\delta}:=\ZIII\cap Z^{\delta}.
\end{equation}

\begin{lemma}
\label{lemma:convex-hol} Take on $\CC^n$ the standard almost
complex structure $I_0$ and metric $g_0$, and suppose that $S^1$
acts diagonally on $\CC^n$ with weights $w_1,\dots,w_n$. Let
$\lambda\in\imag\RR$ and let $\phi:C\to\CC^n$ be a smooth map
satisfying $\ov{\partial}_{\lambda d\theta}\phi=0$. Let
$$l_{\min}:=\min\{|-k+\imag w_j\lambda| \mid k\in\ZZ,\ 1\leq j\leq
n\}.$$ Let $\gamma$ be the function defined in
(\ref{eq:def-gamma}).
\begin{itemize}
\item[(1)] If $l_{\min}\neq 0$ then we have, for some $\delta>0$
depending on $l_{\min}$ but independent of $\phi$:
$$\|d_{\lambda d\theta}\phi\|^2_{L^2(\ZII)}<
\gamma(2l_{\min})( \|d_{\lambda
d\theta}\phi\|^2_{L^2(\ZI^{\delta})} +\|d_{\lambda
d\theta}\phi\|^2_{L^2(\ZIII^{\delta})}).$$ \item[(2)] Suppose that
$\lambda=0$ (so that $\phi$ is holomorphic) and define
$$\phi_{\av}(t,\theta):=\phi(t,\theta)-\frac{1}{2\pi}\int
\phi(t,\nu)d\nu.$$ (Here the subscript $\av$ stands for {\it
average}.) Then, if $\phi$ is not constant, the following holds:
$$\|\phi_{\av}\|^2_{L^2(\ZII)} \leq \frac{1}{e^2+e^{-2}}(
\|\phi_{\av}\|^2_{L^2(\ZI)}+\|\phi_{\av}\|^2_{L^2(\ZIII)}).$$
\end{itemize}
\end{lemma}
\begin{pf}
The condition $\ov{\partial}_{\lambda d\theta}\phi=0$ is
equivalent to
\begin{equation}
\label{eq:es-holomorfa}\frac{\partial \phi}{\partial t}
=I_0\left(\frac{\partial \phi}{\partial\theta}-\imag\lambda
\XXX(\phi)\right),
\end{equation}
where, denoting by $(\phi_1,\dots,\phi_n)$ the coordinates of
$\phi$, we have $\XXX(\phi)=(\imag w_1\phi_1,\dots,\imag
w_n\phi_n)$. Writing for each $j$ the Fourier expansion
$\phi_j(t,\theta)=\sum_{k\in\ZZ}a_{k,j}(t)e^{\imag k\theta}$
formula (\ref{eq:es-holomorfa}) is equivalent to the equation
$a_{k,j}'=(-k+\imag \lambda w_j)a_{k,j}$ for each Fourier
coefficient. Hence $a_{k,j}(t)=a_{k,j}(0)e^{(-k+\imag\lambda
w_j)t}$. Also, since $\ov{\partial}_{\lambda d\theta}\phi=0$ we
have $|d_{\lambda d\theta}\phi|^2=2|\partial\phi/\partial t|^2$.
Putting everything together we write
\begin{equation}
\label{eq:energia-fourier} \int_{[t_0,t_1]\times S^1}|d_{\lambda
d\theta}\phi|^2 =\sum_{j=1}^n\sum_{k\in\ZZ \atop \imag
k+\lambda\neq 0} \int_{t_0}^{t_1}2|a_{k,j}(0)|^2|\imag
k+\lambda|^2 e^{2(-k+\imag\lambda w_j)t} dt.
\end{equation}
The proof now follows from
Lemma \ref{lemma:integrals-dexponencials} below.
This proves (1). The formula in (2) is proved similarly, noting
that the integral which computes $\|\phi_{\av}\|^2$
is of the form $\int \sum_{j\in\ZZ\setminus\{0\}}|a_{k,j}|^2
e^{2jt}dt$.
\end{pf}

\begin{lemma}
\label{lemma:integrals-dexponencials}
For any $\eta_0>0$ there is some $\delta>0$ such that
for every $\eta\in\RR$ satisfying $|\eta|\geq\eta_0$ we have
\begin{equation}
\label{eq:integrals-dexponencials}
\int_0^1 e^{\eta x}dx<\gamma(\eta_0)
\left(\int_{-1+\delta}^0 e^{\eta x}dx
+\int_1^{2-\delta} e^{\eta x}dx\right).
\end{equation}
\end{lemma}
\begin{pf}
By symmetry it suffices to consider the case $\eta>0$.
A simple computation shows that the inequality
\begin{equation}
\label{eq:sense-delta}
\int_0^1 e^{\eta x}dx<\beta
\left(\int_{-1}^0 e^{\eta x}dx
+\int_1^{2} e^{\eta x}dx\right)
\end{equation}
is equivalent to $\beta>1/(e^{\gamma}+e^{-\gamma})$. Since the
function $f(x):=e^x+e^{-x}$ is increasing for positive $x$, if we
set $\beta:=\gamma(\eta_0/2)$ then (\ref{eq:sense-delta}) holds
for every $\eta\geq\eta_0$. On the other hand, if we set
$\delta:=1/2$ then (\ref{eq:integrals-dexponencials}) holds for
every $\eta$ bigger than some $\eta'>\eta_0>0$. Indeed, the left
hand side in (\ref{eq:integrals-dexponencials}) grows as a
function of $\eta$ as $e^{\eta}$, whereas the right hand side
grows as $e^{\eta(2-\delta)}=e^{3\eta/2}$. On the other hand,
since the inequality (\ref{eq:sense-delta}) is strict for
$\beta=\gamma(\eta_0)$ and every $\eta\geq \eta_0$, we can pick
some very small $0<\delta<1/2$ such that
(\ref{eq:integrals-dexponencials}) holds for every
$\eta\in[\eta_0,\eta]$. Then (\ref{eq:integrals-dexponencials})
will also hold for every $\eta\geq\eta_0$ (of course, if for a
given $\eta$ (\ref{eq:integrals-dexponencials}) is true for some
$\delta=\delta_0$, then it is also true for any $0\leq
\delta\leq\delta_0$).
\end{pf}

\subsection{Compactness of the set of pairs of small energy}

A sequence of pairs $\{(\alpha_n,\phi_n):C\to X\}$ will be said to
{\bf converge} to a pair $(\alpha,\phi):C\to X$ if
$\alpha_n\to\alpha_0$ in the $L^{\infty}$ norm and
$\phi_n\to\phi_0$ in $L^p_{1,\loc}$ for any $p<\infty$. This
implies that for any compact subset $K\subset C$ the energies
$\|d_{\alpha_n}\phi_n\|_{L^2(K)}$ converge to
$\|d_{\alpha_0}\phi_0\|_{L^2(K)}$, and that the maps $\phi_n$
converge to $\phi_0$ in the continuous topology.

The following lemma is, except for the last statement, a
combination of Lemma 3.1 and Corollary 3.3 in \cite{IS} (see
Definition 3.1 in \cite{IS} for the notion of uniformly continuous
almost complex structure on a manifold).

\begin{lemma}
\label{lemma:IS-2} Let $Y$ be a (nonnecessarily compact) manifold.
Let $h$ be some metric on $Y$, and let $I_0$ be a continuous
almost complex structure on $Y$. Let $Z\subset Y$ be a closed
$h$-complete subset, such that $I_0$ is uniformly continuous on
$Z$ w.r.t. $h$. There exists a real number
$\epsilon=\epsilon(I_0,Z,h)$ with the following property. Let
$\{I_n\}$ be a sequence of continuous almost complex structures on
$Y$ such that $I_n\to I_0$ in $C^0$-topology of $Y$. Let $u_n\in
C^0\cap L^2_{1,\loc}(\DD,Y)$ be a sequence of $I_n$-holomorphic
maps such that $u_n(\DD)\subset Z$,
$\|du_n\|_{L^2(\DD)}\leq\epsilon$ and $u_n(0)$ is bounded in $X$.
Then there exists a subsequence $\{u_{n_k}\}$ which $L^p_{1,\loc}$
converges to a $I_0$-holomorphic map $u_{\infty}$ for all
$p<\infty$. In particular, for any $K\Subset\DD$ the norms
$\|du_n\|_{L^2(K)}$ tend to $\|du_{\infty}\|_{L^2(K)}$.
Furthermore, if for some $p>2$ and $k>0$ the $L^p_k$ norms of
$I_n$ are uniformly bounded and $\|I_n-I_0\|_{L^p_k}\to 0$, then
$u_{\infty}$ is of class $L^p_k$ and the subsequence can be chosen
to converge to $u_{\infty}$ in $L^p_{k}$.
\end{lemma}
\begin{pf}
We prove the last statement: by Proposition B.4.7 in \cite{McDS},
each $\phi_n$ is of class $L^p_{k+1}$, and since $L^p_{k+1}\to
L^p_k$ is compact, passing to a subsequence there is a limit
$\phi_n\to\phi$ in $L^p_k$.
\end{pf}

Using the previous lemma we prove the following result on
convergence of pairs with small energy.

\begin{lemma}
\label{lemma:easy-compactness} Let $\{I_n\}$ be a sequence of
continuous almost complex structures on $X$ such that $I_n\to I$
in the $L^{\infty}$ norm, where $I$ is also a continuous almost
complex structure. Let $C$ be a (nonnecessarily compact) complex
curve with a Riemannian metric. Let $(\alpha_n,\phi_n):C\to X$ be
a sequence of $I_n$-holomorphic pairs. Assume that
$\|d\alpha_n\|_{L^{\infty}}$ is uniformly bounded, and that
$\|d_{\alpha_n}\phi_n\|_{L^2(C)}\leq \epsilon/2$ (where $\epsilon$
is as in Lemma \ref{lemma:IS-2}). Then, for any compact subset
$M\subset C$, there exists a subsequence of pairs
$\{(\alpha_{n_k},\phi_{n_k})\}$ whose restriction to $M$
converges, up to regauging, to a $I$-holomorphic pair
$(\alpha,\phi):M\to X$. If for some $p>2$ and $k>0$ the $L^p_k$
norms of $\alpha_n$, $I_n$, $\alpha$ and $I$ are uniformly
bounded, then $\phi$ is of class $L^p_k$ and the subsequence can
be taken so that $\phi_n\to\phi$ in $L^p_{k}$.
\end{lemma}

\begin{pf}
Since $d\alpha_n$ has uniformly bounded $L^{\infty}$ norm we can
assume, up to regauging, that the sequence $\alpha_n$ is uniformly
continuous and $L^{\infty}$ bounded. By Ascoli--Arzela it follows
that there is a subsequence which converges in $L^{\infty}$ norm
to a continuous $1$-form $\alpha\in\Omega^1(D,\imag\RR)$. Take the
corresponding sequence of pairs and denote it again by
$\{(\alpha_n,\phi_n)\}$, so that $\alpha_n$ converges to $\alpha$
in the $L^{\infty}$ norm as $n\to\infty$. By a standard patching
argument in gauge theory (see \S 4.4.2 in \cite{DK}) it suffices
to prove the result of the lemma for a finite collection of disks
covering $C$. Let $h_0$ be the product metric on $Y=C\times X$,
let $\Phi_n:C\to Y$ be the section corresponding to $\phi_n$, and
let $h_n:=g(\alpha_n)$ and $h:=g(\alpha)$. Then we have (see
(\ref{eq:norma-dPhi}))
$|d\Phi_n|^2_{h_n}=1+|d_{\alpha_n}\phi_n|^2.$ Since
$\alpha_n\to\alpha$ in $L^{\infty}$ we deduce that for a big
enough $n$ we can bound $h_n<\sqrt{2} h$. It follows that, if
$D\Subset\DD$ is a disk of area $\epsilon^2/(\sqrt{2}K)$, where
$K$ is a suitable constant depending on the metrics $h_0$ and $h$,
then $\|d\Phi_n\|_{h,L^2(D)}\leq\epsilon.$ Hence, if $f:\DD\to D$
is an affine biholomorphism, then $\|d(\Phi_n\circ
f)\|_{h,L^2(\DD)}\leq\epsilon$. Finally, it follows from
$\alpha_n\to\alpha$ that $I_n(\alpha_n)\to I(\alpha)$ in
$L^{\infty}$ norm. Taking into account (\ref{eq:holomorf}), we can
apply Lemma \ref{lemma:IS-2} (it is clear that $Z:=D\times X$
satisfies the hypothesis of Lemma \ref{lemma:IS-2}) to the
sequence of maps $\{\Phi_n\circ f\}$, and deduce that a
subsequence converges to a map defined on $D$. The last statement
follows from the last statement in Lemma \ref{lemma:IS-2} (since
$p>2$ and $k>1$, if $\alpha_n$ has bounded $L^p_k$ norm then
$I(\alpha_n)$ is also bounded in $L^p_k$).
\end{pf}

In the following two corollaries we use the same notation as in
Section \ref{ss:inequality-cylinders}, so $C$ denotes the cylinder
$[-1,2]\times S^1$ and $Z=[0,1]\times S^1$.

\begin{corollary}
\label{cor:pocaenergia-diampetit} Let $I$ be an almost complex
structure on $X$. For any $\epsilon>0$ there exists some
$\delta>0$ such that if $(\alpha,\phi):C\to X$ is a
$I$-holomorphic pair satisfying
$\|d\alpha\|_{L^{\infty}(C)}<\delta$ and
$\|d\phi\|_{L^2(C)}<\delta$ then
$$\diam_{S^1}\phi(Z)<\epsilon.$$
\end{corollary}
\begin{pf}
We prove the corollary by contradiction. Take any $\epsilon>0$ and
assume that there is a real number $\delta>0$ and a sequence of
$I$-holomorphic pairs $(\alpha_n,\phi_n):C\to X$ such that
$\|d\alpha_n\|_{L^{\infty}(C)}\to 0$,
$\|d_{\alpha_n}\phi_n\|_{L^2(C)}\to 0$ and
$\diam_{S^1}\phi_n(Z)\geq\delta$. By Lemma
\ref{lemma:easy-compactness} we can assume, up to restricting to a
subsequence and regauging, that there is a $I$-holomorphic pair
$(\alpha,\phi):Z\to X$ such that the restriction of
$\{(\alpha_n,\phi_n)\}$ to $Z$ converges to $(\alpha,\phi)$. This
implies in particular that $\|d_{\alpha}\phi\|_{L^2(Z)}=0$ and that
$\diam_{S^1}\phi(Z)\geq\delta$, which is impossible.
\end{pf}

\begin{corollary}
\label{cor:cota-inferior} For any $p>2$, $K$ and $\epsilon>0$
there is some $\delta>0$ with the following property. Suppose that
$(\alpha,\phi):Z\to X$ is a $I$-holomorphic pair, that
$\|d\alpha\|_{L^p_1}\leq K$, and that $\sup_{\{1/2\}\times
S^1}|d_{\alpha}\phi|\geq\epsilon$. Then
$\|d_{\alpha}\phi\|_{L^2(Z)}\geq\delta$.
\end{corollary}
\begin{pf}
Fix $p$, $K$ and $\epsilon$, and suppose there is no $\delta>0$
satisfying the hypothesis of the Corollary. Then there is a
sequence of pairs $(\alpha_u,\phi_u):Z\to X$ such that
$\|d\alpha_u\|_{L^p_1}\leq K$, $\sup_{\{1/2\}\times
S^1}|d_{\alpha_u}\phi_u|\geq\epsilon$ and
$\|d_{\alpha_u}\phi_u\|_{L^2(Z)}\to 0$. Using the last statement
in Lemma \ref{lemma:easy-compactness}, we deduce that there is a
subsequence (which we denote again by $(\alpha_u,\phi_u)$) which
converges to $(\alpha,\phi)$ in $L^p_2$. In particular,
$d_{\alpha}\phi=0$, but also $d_{\alpha_u}\phi_u$ converges
pointwise to $d_{\alpha}\phi$, which is a contradiction.
\end{pf}

\subsection{Proof of Lemma \ref{lemma:ineq-energia}}
\label{proof:lemma:ineq-energia}

The strategy will be to reduce the lemma to the case of the standard
complex structure $I_0$ on $\CC^n$ and $\alpha=\lambda d\theta$.
We will do this by means of a compactness argument. Before stating
the argument, we need some preliminaries.
Permuting the coordinates if necessary, we can assume that
for some $p$ we have
\begin{equation}
\label{eq:weights-sorted} \text{for any $1\leq j\leq p$, $w_j\neq
0$} \qquad \text{and} \qquad w_{p+1}=\dots=w_{n}=0.
\end{equation}
Let $q=n-p$, so that we have a splitting $\CC^n=\CC^p\times \CC^q$.
Let $\pi_p$, $\pi_q$ be the projections from $\CC^n$ to $\CC^p$ and
$\CC^q$.
The following lemma will be proved in Section \ref{ss:proof-C0-bounds}.

\begin{lemma}
\label{lemma:C0-bounds}
There exists some constant $K>0$ with the following property.
Let $(\alpha,\phi):Z\to\CC^n$ be a $I$-holomorphic pair such that
$\|\alpha-\lambda d\theta\|_{L^{\infty}}<|\lambda|/2$.
For any $0<\delta<1/2$ there exists some $y\in\CC^q$ such that
$$\sup_{z\in Z^{\delta}}
|\phi(z)-(0,y)|\leq l_{\min}^{-1}\delta^{-1/2}K
(\|d_{\alpha}\phi\|_{L^2}+1).$$
\end{lemma}

We first prove (\ref{eq:average-bound-energy-loc}).
Suppose that no matter how small $\epsilon$ is there are pairs
for which (\ref{eq:average-bound-energy-loc}) does not hold.
Then one can chose a sequence of equivariant
almost complex structures $I_u$ and
$I_u$-holomorphic pairs $(\alpha_u,\phi_u):C\to\CC^n$
satisfying:
$$\|I_u-I_0\|_{L^{\infty}}\to 0,\qquad
\|\alpha_u-\lambda d\theta\|_{L^{\infty}}\to 0,\qquad
\|d_{\alpha_u}\phi_u\|_{L^2(C)}\to 0,$$
and for which the inequality (\ref{eq:average-bound-energy-loc})
holds in the opposite direction:
\begin{equation}
\label{eq:desoposada}
\|d_{\alpha}\phi\|^2_{L^2(\ZII)}
>\gamma(2 l_{\min})\left(
\|d_{\alpha}\phi\|^2_{L^2(\ZI)}+
\|d_{\alpha}\phi\|^2_{L^2(\ZIII)}\right).
\end{equation}
In particular, the energy $E_u$ of $(\alpha_u,\phi_u)$ is nonzero.
Let $\epsilon$ be less that the epsilons appearing in Lemmae
\ref{lemma:C0-bounds} and \ref{lemma:easy-compactness}, and define
for every $n$ a new pair $(\alpha_u',\phi_u')$ by setting
$\alpha_u':=\alpha_u$ and $\phi_u'(z):=\phi_u(z)\epsilon/2E_u$.
Let also $I_u'$ be the pullback of $I_u$ by the homotecy of ratio
$2E_u/\epsilon$. Then $(\alpha_u',\phi_u')$ is a
$I_u'$-holomorphic pair and its energy is $\epsilon'/2$.
Furthermore, since $I_0$ is invariant under homotecies, it follows
that $\|I_u'-I_0\|_{L^{\infty}}\to 0$.

Let $\delta>0$ be the number given by (1) in Lemma
\ref{lemma:convex-hol} for our value of $l_{\min}$. Let
$y_u\in\CC^q$ a point, as given by Lemma \ref{lemma:C0-bounds},
such that $\phi_n'(Z^\delta)$ is contained in the ball centered at
$(0,y)$ and of radius $K l_{\min}^{-1}\delta^{-1}$. Finally, let
$\phi_u'':=\phi_u'-(0,y_u)$ and let $I_u''$ be the pullback of
$I_u'$ by the translation in the direction $(0,y_u)$. It follows
that $(\alpha_u',\phi_u'')$ is $I_u''$-holomorphic (note that
$I_u''$ is $S^1$-equivariant) and $\phi_u''(Z_\delta)$ is
contained in the ball centered at $0$ of radius $(\epsilon+1)
K l_{\min}^{-1}\delta^{-1/2}$. At this point we use Lemma
\ref{lemma:easy-compactness} to deduce that, up to regauging,
there is a subsequence of $\{(\alpha_u',\phi_u'')\}$ which
converges to a $I_0$-holomorphic pair $(\lambda
d\theta,\phi):Z^\delta\to\CC^n$ which has energy $\epsilon/2$ and
such that
\begin{equation}
\label{eq:cosa-rara} \|d_{\lambda d\theta}\phi\|_{L^2(\ZII)}\geq
\gamma(2l_{\min})( \|d_{\lambda
d\theta}\phi\|_{L^2(\ZI^{\delta})} +\|d_{\lambda
d\theta}\phi\|_{L^2(\ZIII^{\delta})}).
\end{equation}
Now, comparing this with (1) in Lemma \ref{lemma:convex-hol} we
obtain a contradiction. This finishes the proof of
(\ref{eq:average-bound-energy-loc}).

The inequality (\ref{eq:pointwise-bound-energy-loc}) can be reduced,
following the same strategy, to the case of the standard complex structure
$I_0$ on $\CC^n$ and the connection $\alpha=\lambda d\theta$.
Suppose for simplicity that $n=1$ and that $w_1=1$ (the general case
offers no extra difficulty). Let $\phi:C\to\CC$ be a map which satisfies
\begin{equation}
\label{eq:cond-on-phi}
\frac{\partial \phi}{\partial t}=\imag\left(
\frac{\partial \phi}{\partial\theta}
+\lambda\phi\right).
\end{equation}
We want to prove that
$$\sup_Z \left|\frac{\partial \phi}{\partial t}\right|
\leq K \left\|\frac{\partial \phi}{\partial t}\right\|_{L^2(C)}$$
for some constant $K$. Let $\psi:=\partial \phi/\partial t$. Then
$\psi$ satisfies the same equation (\ref{eq:cond-on-phi}) as
$\phi$. It follows from that that the function
$\zeta=e^{-\imag\lambda t}\psi$ is holomorphic on $C$. Define the
constants $K_0:=\sup_{t\in[0,1]} e^{\imag\lambda t}$ and
$K_1:=\sup_{t\in[0,1]} e^{-\imag\lambda t}.$ Take any $z\in Z$ and
denote by $\DD\subset C$ the disk of radius $1$ centered at $z$.
Since $\zeta$ is holomorphic we have $|\zeta(z)|\leq
\|\zeta\|_{L^2(\DD)}$. Now we bound:
$$|\psi(z)|\leq K_0 |\zeta(z)|\leq K_0\|\zeta\|_{L^2(\DD)}
\leq K_0 K_1 \|\psi\|_{L^2(\DD)} \leq K_0 K_1 \|\psi\|_{L^2(C)},$$
which is what we wanted to prove.

\subsection{Proof of Lemma \ref{lemma:C0-bounds}}
\label{ss:proof-C0-bounds} First of all we state and prove two
lemmae which will be necessary for proving Lemma
\ref{lemma:C0-bounds}. For the first one we introduce the
following notation. Let $M$ be the manifold $\RR\times S^1\times
\CC^n$. Take a $1$-form $\alpha$ on $\RR\times S^1$ with values in
$\imag\RR$, and let $J:=I(\alpha)$ be the complex structure on $M$
induced by $\alpha$ and the action of $S^1$ on $\CC^n$ with
weights $w_1,\dots,w_n$ as in (\ref{eq:weights-sorted}). Consider
also the metric $h:=g(\alpha)$ on $M$ induced by $\alpha$ and the
standard metric $g_0$ on $\CC^n$.

\begin{lemma}
\label{lemma:radi-densitat} For any $K'>0$ there exists some
$K''>0$ with the following meaning. Let
$\alpha\in\Omega^1(\RR\times S^1,\imag\RR)$ be a connection
$1$-form. Suppose that $\|\alpha\|_{L^{\infty}}\leq K'$ and that
$\Phi=(i,\phi):Z\to M$ is a $J$-holomorphic map, where
$i:Z\to\RR\times S^1$ is the inclusion; suppose also that for some
$z\in Z$ and $R>0$ we have $d(\phi(z),\phi(\partial Z))>R.$ Then
$$\|d\Phi\|^2_{L^2}=\area(Z)+\|d_{\alpha}\phi\|_{L^2}^2>K''R^2.$$
\end{lemma}

\begin{pf}
This is similar to a standard result in the theory of
pseudo-holomorphic curves. First, since $\Phi$ is $J$-holomorphic,
$\|d\Phi\|^2$ is equal to the area of $\Phi(Z)$, so we have to
prove that $\area(\Phi(Z))>K''R^2$. Define, for any $R\geq r>0$,
the following sets:
$$M_r=\{(t,\theta,x)\in M \mid |x-\phi(z)|\leq r\},
\qquad B_r=\partial M_r,
\qquad Z_r=\Phi^{-1}(M_r).$$
Define also $f(r):=\area(\Phi(Z_r))$.
A real number $r$ will be called regular if it is a regular value
of the function $Z\ni y\mapsto |\phi(y)-\phi(z)|$. For any such $r$,
$\Phi(Z_r)\subset M$ is a smooth subsurface with boundary. Also,
the set of nonregular values has measure zero.

We first prove an upper bound for $f(r)$. There exists some
constant $K_1>0$ such that for any regular $r$ there is a smooth
subsurface $\Sigma\subset M$ such that
$\partial\Sigma=\partial\Phi(Z_r)$ and such that
$\area_h(\Sigma)\leq K_1 \length(\partial\Phi(Z_r))^2$. Let
$\omega(\alpha)$ be the minimal coupling symplectic on $M$ (see
Section \ref{s:symplectic-fibrations}). Since $H^2(M,\RR)=0$ and
since $\omega(\alpha)$ has bounded $L^{\infty}$ norm (because
$\alpha$ is bounded), it follows that, for some constant $K_2$ and
any regular $r$,
$$\int_{\Phi(Z_r)}\omega(\alpha)=-\int_{\Sigma}\omega(\alpha)
\leq K_2\area(\Sigma)\leq K_2K_1 \length(\partial \Phi(Z_r))^2.$$
On the other hand, since $\Phi$ is $I(\alpha)$-holomorphic we can identify
the first integral with the area of $\Phi(Z_r)$. Hence we have proved
\begin{equation}
\label{eq:upper-bound-area}
f(r)\leq K_2K_1 \length(\partial \Phi(Z_r))^2.
\end{equation}
To finish the proof of the lemma, observe that there exists some
constant $K_3$ such that for any regular $r$ we have $f'(r)\geq
K_3\length(\partial \Phi(Z_r)).$ Combining this with
(\ref{eq:upper-bound-area}) we deduce $f'(r)^2\geq K_3
f(r)/(K_2K_1).$ Integrating over $r$ we deduce the result.
\end{pf}

Let $\XXX:\CC^n\to\CC^n$ be the vector field generated by the
infinitesimal action of $S^1$ on $\CC^n$.

\begin{lemma}
\label{lemma:norma-L11} There exists some constant $K_1$,
independent of $\lambda$, such that if $f:S^1\to\CC^p\subset\CC^n$
is any smooth map, and $\beta:S^1\to\imag\RR$ satisfies
$|\beta-\lambda|<\lambda/2$ then $$\sup|f|\leq
l_{\min}^{-1}K_1\|f'+\beta\XXX(f)\|_{L^2(S^1)}.$$
\end{lemma}
\begin{pf}
Let $\beta_0:=\frac{1}{2\pi}\int\beta(\theta)d\theta$. Define
$g:=e^sf$, where $s:S^1\to\imag\RR$ satisfies $ds=\beta-\beta_0$.
We have $|f|=|g|$, and one checks easily that
$|f'+\beta\XXX(f)|=|g'+\beta_0\XXX(g)|$ pointwise. Hence it
suffices to prove that for some universal constant $K_1$
(depending on $\lambda$) and any constant $\beta_0$ satisfying
$|\beta_0-\lambda|<\lambda$ one has the inequality
$$\sup|g|\leq l_{\min}^{-1}K_1\|g'+\beta_0\XXX(g)\|_{L^2(S^1)}.$$
Now, this can be proved by considering the Fourier series of $g$.
\end{pf}

After these preliminaries, we now prove Lemma
\ref{lemma:C0-bounds}. Pick $K_1$ satisfying the requirement of
the previous theorem and in such a way that for any
$h:S^1\to\CC^q\subset\CC^n$ and any $\theta,\eta\in S^1$ we have
\begin{equation}
\label{eq:norma-L2} |h(\theta)-h(\eta)|\leq K_1\|h'\|_{L^2(S^1)}
\end{equation}
(this is possible in view of Cauchy--Schwartz). Let $K''$ be the
constant given by Lemma \ref{lemma:radi-densitat} for the value
$K'=3|\lambda|/2$. Take $K>0$ big enough so that for any positive
$E>0$ we have
\begin{equation}
\label{eq:choice-K} l_{\min}^{-1}\delta^{-1/2}K(E+1)
>6 \max\left\{l_{\min}^{-1}\delta^{-1/2}K_1E,
\sqrt{(E^2+2\pi)/K''}\right\}.
\end{equation}
There exists some $a\in [0,\delta]$ such that
$$\|d_{\alpha}\phi\|_{L^2(\{a\}\times S^1)}^2
\leq\delta^{-1}\|d_{\alpha}\phi\|^2.$$ Denote the components of
$\phi(a,\theta)$ by $(f(\theta),h(\theta))$. Let
$\beta:S^1\to\imag\RR$ be the map such that for any $\theta$ we
have $\alpha(a,\theta)=\beta(\theta)d\theta+\gamma(\theta)dt$.
Applying Lemma \ref{lemma:norma-L11} to the component $f$ and
inequality (\ref{eq:norma-L2}) to the component $h$ (together with
$l_{\min}<1$) we deduce that
$$\sup_{\theta\in S^1}|\phi(a,\theta)-(0,y_a)|\leq
\rho:=l_{\min}^{-1}\delta^{-1/2}K_1\|d_{\alpha}\phi\|_{L^2},$$
where $y_a=h(1)$ (indeed, $f'+\beta\XXX(f)$ can be identified with
$d_{\alpha}\pi_p(\phi)$ to the circle $\{a\}\times S^1$, and we
have $\|d_{\alpha}\pi_p(\phi)\|\leq \|d_{\alpha}\phi\|$). Hence,
$\phi(a,S^1)\subset\CC^n$ is contained in the ball centered at
$z_a:=(0,y_a)$ of radius $\rho$. Similarly, there is some
$b\in[1-\delta,1]$ and $y_b\in\CC^q$ such that
$\phi(b,S^1)\subset\CC^n$ is contained in the ball centered at
$z_b:=(0,y_b)$ of the same radius $\rho$.

Let $Z':=[a,b]\times S^1$ and let
$$R:=l_{\min}^{-1}\delta^{-1/2}K(\|d_\alpha\phi\|_{L^2}+1).$$
We claim that $\phi(Z')$ is contained at least in one of the balls
$B(z_a,R)$ or $B(z_b,R)$. Suppose this is not the case. Then
$\phi(Z')$ is not contained in the union $B(z_a,R/3)\cup
B(z_b,R/3)$. Indeed, if this were true, then, since $\phi(Z')$ is
connected, we should have $B(z_a,R/3)\cap B(z_b,R/3)\neq 0$
(otherwise $\phi(Z')$ would be included either in $B(z_a,R/2)$ or
in $B(z_b,R/2)$, in contradiction with our assumption). But in
this case $d(z_a,z_b)<2R/3$, whence $\phi(Z')\subset
B(z_a,R/3)\cup B(z_b,R/3)\subset B(z_a,R)$, which we assume not to
be true. It follows from all this that there is some
$z\in\phi(Z')$ satisfying simultaneously
$$d(z,z_a)>R/3\qquad\text{ and }\qquad d(z,z_b)>R/3.$$
In other words, the ball $B(z,R/6)$ is disjoint both with
$B(z_a,R/6)$ and $B(z_b,R/6)$. By our choice of $R$ we have
$R\geq 6\rho$, so the discussion above proves that
$\phi(a,S^1)\subset B(z_a,R/6)$ and $\phi(b,S^1)\subset B(z_b,R/6)$.
We deduce that $B(z,R/6)\cap \phi(\partial Z')=\emptyset$.
Hence we can apply Lemma \ref{lemma:radi-densitat} to
$(\phi,\alpha):Z'\to \CC^n$ and deduce that
$$2\pi+\|d_{\alpha}\phi\|^2_{L^2}\geq\area(Z')+\|d_{\alpha}\phi\|^2_{L^2(Z')}
\geq K''R^2/36.$$ However, it follows from (\ref{eq:choice-K})
that we must have $K''R^2/36>2\pi+\|d_{\alpha}\phi\|^2_{L^2},$
which is a contradiction. This finishes the proof.

\section{Long cylinders with small energy and nearly critical residue}
\label{s:long-cylinders}

In this section we use the same notations as in the previous one.

For any natural number $N$, let $C_N=[-N,N]\times S^1$ with the
standard product metric and the induced conformal structure. We
denote as always by $(t,\theta)$ the usual coordinates in $C_N$.
Let $v=f dt\wedge d\theta$ be a volume form, and let $\eta>0$ be
any number. We say that $v$ is {\bf exponentially
$\eta$-bounded} if
\begin{equation}
\label{eq:exp-bounded} |f|<\eta e^{|t|-N}.
\end{equation}
Volume forms satisfying this property arise when we consider $C_N$
as a conformal model for open sets $U$ in stable curves $(C,\bx)$
of the form $U\simeq \{xy=\delta\}\subset\CC^2$ for small
$\delta$, and we take on $C_N$ the volume form corresponding to
the metric $\nu_{[C,\bx]}$.

\begin{theorem}
\label{thm:small-energy} Fix some number $c\in\imag\RR$. There
exist numbers $\epsilon>0$, $\eta_0>0$, $\sigma>0$ and $K>0$
(depending only on $X$, the action of $S^1$ on $X$, and $I$) with
the following property. Let $(\alpha,\phi):C_N\times S^1\to X$ be
a pair. Suppose that $\ov{\partial}_{I,\alpha}\phi = 0$ and that
there exists some critical residue
$\lambda_{\crit}\in\Lambda_{\crit}$ such that
\begin{equation}
\label{eq:condicions-pair}
\|\alpha-\lambda_{\crit} d\theta\|_{L^{\infty}}<\epsilon
\qquad\text{and}\qquad
\|d_{\alpha}\phi\|_{L^{\infty}}<\epsilon.
\end{equation}
Then there exist maps $\psi:[-N,N]\to X^{\lambda_{\crit}}$ and
$\phi_0:C_N\to TX$ satisfying $\phi_0(t,\theta)\in T_{\psi(t)}X$,
\begin{equation}
\label{eq:definicio-phi0} \phi(t,\theta)= e^{-\lambda_{\crit}
\theta}
\exp^g_{\psi(t)}(e^{\lambda_{\crit}\theta}\phi_0(t,\theta))
\qquad\text{and}\qquad \int
e^{\lambda_{\crit}\theta}\phi_0(t,\theta)d\theta=0,
\end{equation}
where $\exp^g_x:T_xX\to X$ denotes the exponential map on $X$ with
respect to the metric $g$ (see also Remark \ref{rmk:ambiguitats}
for the precise meaning of the formulae). We distinguish two
cases.
\begin{enumerate}
\item Let $v$ be an exponentially $\eta$-bounded volume form on
$C_N$, where $\eta<\eta_0$. Suppose that
\begin{equation}
\label{eq:equacions-a-phi-vor} \iota_{v}d\alpha+\mu(\phi) = c;
\end{equation}
then the following inequality holds:
\begin{equation}
\label{eq:estimate-phi0} |\phi_0(t,\theta)|< K e^{-\sigma(N-|t|)}
\left( \|d_{\alpha}\phi\|_{L^{\infty}}+\|\alpha\|_{L^{\infty}}
+\eta(|c|+\sup|\mu|)\right)^{1/4}.
\end{equation}
\item Suppose now that $d\alpha=0$;
then the previous estimate can be improved to:
\begin{equation}
\label{eq:estimate-phi0-flat} |\phi_0(t,\theta)|< K e^{-\sigma(N-|t|)}
(\|d_{\alpha}\phi\|_{L^{\infty}}+\|\alpha\|_{L^{\infty}})^{1/4}.
\end{equation}
\end{enumerate}
\end{theorem}

\begin{remark}
\label{rmk:ambiguitats}
Note that the first formula in (\ref{eq:definicio-phi0}) is well defined.
To see this,
denote $\lambda_{\crit}=\imag p/q$, where $p$ and $q$ are relatively
primer integers. Then the condition $\psi\subset X^{\lambda_{\crit}}$
implies that
for any $t$ the point $\psi(t)$ is fixed by the action of
$e^{2\pi\imag/q}$. On the other hand, $e^{\pm\lambda_{\crit}\theta}$ is well
defined up to multiplication by powers of $e^{2\pi\imag/q}$. Finally,
since the exponential map is equivariant we have
$e^{-2\pi\imag/q}\exp_{\psi(t)}(e^{2\pi\imag/q}v)=\exp_{\psi(t)}(v)$
for any vector $v\in T_{\psi(t)}X$.
As for the second formula in (\ref{eq:definicio-phi0}), it should be
understood as
$$\int e^{\lambda_{\crit}\theta}\phi(t,\theta)d\theta=\int_0^{2\pi q}
e^{\lambda_{\crit}\theta}\phi(t,\theta)d\theta.$$
\end{remark}

The proof of Theorem \ref{thm:small-energy} will be given in
Section \ref{ss:proof-thm:small-energy} below. Next theorem states
that the map $\psi$ constructed in Theorem \ref{thm:small-energy}
in approximately a gradient flow line for the moment map $\mu$.
Before stating it we introduce some notation. Suppose that
$(\alpha,\phi)$ is a pair satisfying the hypothesis of Theorem
\ref{thm:small-energy}. Applying a gauge transformation if
necessary, we can assume that that $\alpha$ is in temporal gauge,
$\alpha=\alpha_{\theta}d\theta$, and that $\alpha_{\theta}$
restricted to $\{0\}\times S^1$ takes a constant value
$\lambda+\lambda_{\crit}\in\imag\RR$. Then we can write
\begin{equation}
\label{eq:def-lambda-beta}
\alpha=(\lambda+\lambda_{\crit}+\beta)d\theta
\end{equation}
for some function $\beta:C_N\to\imag\RR$.

\begin{theorem}
\label{thm:psi-gradient} Following the notation of Theorem
\ref{thm:small-energy}, we distinguish again two cases.
\begin{enumerate}
\item Suppose that $\iota_{v}d\alpha+\mu(\phi)=c$, where $v$
is $\eta$-bounded and $\eta<\eta_0$; then
\begin{equation}
\label{eq:der-psi} |\psi'(t)+\imag\lambda
I(\psi(t))\XXX(\psi(t))|<Ke^{-\sigma (N-|t|)} \left(
\|d_{\alpha}\phi\|_{L^{\infty}}+\|\alpha\|_{L^{\infty}}
+\eta(|c|+\sup|\mu|)\right)^{1/4},
\end{equation}
holds for any $t$, where $K$ is a constant depending only on $X$,
the action of $S^1$, and the almost complex structure $I$; \item
suppose that $d\alpha=0$; then, for any $t$, we have
\begin{equation}
\label{eq:der-psi-flat} |\psi'(t)+\imag\lambda
I(\psi(t))\XXX(\psi(t))|<Ke^{-\sigma (N-|t|)}
(\|d_{\alpha}\phi\|_{L^{\infty}}+\|\alpha\|_{L^{\infty}})^{1/4},
\end{equation}
\end{enumerate}
\end{theorem}

The proof of Theorem \ref{thm:psi-gradient} will be given in
Section \ref{ss:proof-thm:psi-gradient}

\subsection{Proof of Theorem \ref{thm:small-energy}}
\label{ss:proof-thm:small-energy}

We first prove the theorem under the assumption that the critical
residue $\lambda_{\crit}$ is equal to $0$.

The existence of $\psi$ and $\phi_0$ satisfying
(\ref{eq:definicio-phi0}) follows easily from the implicit function
theorem (see for example Section 14 in \cite{FO}) provided $\epsilon$
is small enough (so that $\|\alpha\|_{L^{\infty}}<\epsilon$ together
with $\|d_{\alpha}\phi\|_{L^2}<\epsilon$ imply that for any $t$
the image of $\phi(t,\cdot):S^1\to X$ is contained in a small ball
of radius less than the injectivity radius of $X$ and $g$).

To prove that the estimate (\ref{eq:estimate-phi0}) follows
from (\ref{eq:equacions-a-phi-vor})
we use the same strategy as in the proof of Theorem
\ref{proof:thm:exp-decay}, namely, we make use of
an inequality involving
the $L^2$ norm of $\phi_0$ restricted to three consecutive pieces
of the cylinder and then using the fact that $(\alpha,\phi)$ is
$I$-holomorphic to deduce pointwise bounds on $\phi_0$ from $L^2$
bounds. To state the inequality we define as in the previous
section $C=[-2,3]\times S^1$. Recall also that when we write any
norm ($L^\infty$, $L^2$, etc.) of either the connection or the
section of a pair defined over $C$, or of a complex structure
depending on points of $C$, we mean the norm over $C$ (unless we
specify some other domain). Let also $\ZI$, $\ZII$ and $\ZIII$ be
the subsets of $C$ defined
in (\ref{eq:definicio-CIII}).

\begin{theorem}
\label{thm:small-energy-average} There exist some constants
$\epsilon>0$ and $K>0$, depending only on $X$ and $I$, with the
following property. Let $(\alpha,\phi):C\times S^1\to X$ be a
$I$-holomorphic pair satisfying $\|\alpha\|_{L^{\infty}}<\epsilon$
and $\|d_{\alpha}\phi\|_{L^{\infty}}<\epsilon.$ Let
$\psi:[-2,3]\to X$ and $\phi_0:C\to TX$ be defined as in Theorem
\ref{thm:small-energy}.
\begin{itemize}
\item[(1)] Suppose that
$\|d\alpha\|_{L^2(C)}<\epsilon\|\phi_0\|_{L^2(C)}$.
Then
\begin{equation}
\label{eq:average-norma-phi0}
\|\phi_0\|^2_{L^2(\ZII)}\leq
\frac{1}{e+e^{-1}}
\left(
\|\phi_0\|^2_{L^2(\ZI)}+
\|\phi_0\|^2_{L^2(\ZIII)}\right)
\end{equation}
\item[(2)] Suppose now that $\|d\alpha\|_{L^2_1}<K$ and
$\|d\alpha\|_{L^{\infty}}<K.$ Then $\sup_Z|\phi_0|\leq
K\|\phi_0\|_{L^2}^{1/4}.$
\end{itemize}
\end{theorem}

Theorem \ref{thm:small-energy-average} will be proved in Section
\ref{proof:thm:small-energy-average}, and we now prove Theorem
\ref{thm:small-energy}. Take $\epsilon$ as in Theorem
\ref{thm:small-energy-average}. Since
$d_\alpha\phi=d\phi-\imag\alpha\XXX(\phi)$, we can bound
\begin{equation}
\label{eq:dphi-daphi} \|d\phi\|_{L^{\infty}}\leq
K(\|d_{\alpha}\phi\|_{L^{\infty}}+\|\alpha\|_{L^{\infty}}).
\end{equation}
As a consequence, and taking (\ref{eq:condicions-pair}) into
account, we deduce that $\|d\phi\|_{L^{\infty}}$ is uniformly
bounded. Let $v=f dt\wedge d\theta$. Equation
(\ref{eq:equacions-a-phi-vor}) can be written as $d\alpha=f(c-\mu(\phi))
dt\wedge d\theta.$ Since $\mu(\phi)$ is uniformly bounded (because
$X$ is compact), once $c$ has been chosen we can take $\eta_0$ in
such a way that the requirement of $v$ being exponentially
$\eta$-bounded (and hence $\eta_0$-bounded as well) implies that
$\|d\alpha\|_{L^2_1(C_N)}$ is necessarily less
than the $K$ in Theorem \ref{thm:small-energy-average}
(here we use that $\|d\phi\|_{L^{\infty}}$ is uniformly bounded).
On the other hand, we have
\begin{equation}
\label{eq:dalpha-exp} |d\alpha|(t,\theta)<\eta (|c|+\sup|\mu|)e^{|t|-N}.
\end{equation}
In the sequel we will denote for convenience
$R:=\eta(|c|+\sup|\mu|)$.

Define, for every $-N\leq n\leq N$, $Z_n=[n,n+1]$ and
$x_n:=\|\phi_0\|_{L^2(Z_n)}$.
Formula (\ref{eq:dphi-daphi}) implies that for any $n$ we have
\begin{equation}
\label{eq:fj-puntual} x_n \leq
K(\|d_{\alpha}\phi\|_{L^{\infty}}+\|\alpha\|_{L^{\infty}}).
\end{equation}
Let also $z_n:=\|d\alpha\|_{L^2(Z_n)}$.
Formula (\ref{eq:dalpha-exp}) implies that
$z_n\leq K_0Re^{-(N-|n|)}$ for every $n$, where $K_0$ is a universal
constant.
On the other hand, Theorem \ref{thm:small-energy-average}
implies that if
$$z_{n-2}+\dots+z_{n+2}\leq\epsilon(x_{n-2}+\dots+x_{n+2})$$
then $x_n\leq (e+e^{-1})^{-1}(x_{n-1}+x_{n+1}).$ Hence the
sequences $\{x_n\}$ and $\{z_n\}$ satisfy the requirements of
Lemma \ref{lemma:recur-perturb}. Applying Lemma
\ref{lemma:recur-perturb} together with (\ref{eq:fj-puntual}) we
deduce, for some $\sigma_0>0$ independent of $(\alpha,\phi)$, an
estimate
$$\|\phi_0\|_{L^2(Z_j)}=x_j\leq \epsilon^{-1}K_0
(R+\|d_{\alpha}\phi\|_{L^{\infty}} +\|\alpha\|_{L^{\infty}})
e^{-\sigma_0(N-|j|)}.$$ To finish the proof of the theorem and
obtain (\ref{eq:estimate-phi0}), combine the previous inequality
with (2) of Theorem \ref{thm:small-energy-average} (the two
necessary conditions for section (2) of Theorem
\ref{thm:small-energy-average} are satisfied in our situation: the
$L^2_1$ norm of $d\alpha$ on the whole $C_N$ is less than $K$ and
$d\alpha$ has bounded $L^{\infty}$ norm by (\ref{eq:dalpha-exp})).

The proof that $d\alpha=0$ implies (\ref{eq:estimate-phi0-flat})
follows exactly the same scheme (and is even easier).

Now we consider the case of general critical residue
$\lambda_{\crit}\in\Lambda_{\crit}$. Let us write
$\lambda_{\crit}=\imag p/q$ for some
relatively prime integers $p,q$ satisfying $q\geq 1$.
Consider the covering map
$$\pi:C_{N/q}\to C_N$$ defined as
$\pi(t,\theta):=(qt,q\theta)$, and define a new pair
$(\alpha',\phi'):=\pi^*(\alpha,\phi)$. The first inequality in
(\ref{eq:condicions-pair}) is equivalent to
$\|\alpha'-pd\theta\|_{L^{\infty}}<q\epsilon$. Applying the gauge
transformation $g(t,\theta):=e^{\imag p\theta}$ we obtain another
pair $(\alpha'',\phi''):=g^*(\alpha',\phi')
=(\alpha'-pd\theta,e^{\imag\theta p}\phi').$ The new connection
$\alpha''$ now satisfies $\|\alpha''\|_{L^{\infty}}<q\epsilon$,
$\|d_{\alpha''}\phi''\|_{L^{\infty}}<q\epsilon$ and
$\ov{\partial}_{I,\alpha''}\phi''=0$. Also, equation
(\ref{eq:equacions-a-phi-vor}) is satisfied with $v$ replaced by
$\pi^*v$ (which now is exponentially $q\eta$-bounded). Hence,
provided $\epsilon$ and $\eta$ are small enough, we can apply the
case of zero residue to the pair $(\alpha'',\phi'')$. We thus
arrive at a pair of maps $\psi'':[-N/q,N/q]\to X$ and
$\phi_0'':C_{N/q}\to TX$ satisfying $\phi''=\exp_{\psi''}\phi_0''$
and $\int\phi_0''(t,\nu)d\nu=0$. We now prove that $\psi''\subset
X^{\lambda_{\crit}}$. Let $\alpha:=2\pi/q$. For any $t,\theta$ we
compute
$$\phi''(t,\alpha+\theta)=e^{\imag (\alpha+\theta)p}\phi'(t,\alpha+\theta)
=e^{\imag \alpha p}e^{\imag \theta p}\phi'(t,\theta)
=e^{\imag \alpha p}\phi''(t,\theta),$$
which, combined with the fact that $g$ is $S^1$-invariant,
implies $\psi''\subset X^{\lambda_{\crit}}$.
The computation above implies similarly that
$\phi''_0(t,\alpha+\theta)=e^{\imag \alpha p}\phi''_0(t,\theta),$
so that the map
$$\tphi''_0(t,\theta):=e^{-\imag\theta p}\phi''_0(t,\theta)$$
satisfies $\tphi''_0(t,\theta)=\tphi''_0(t,\alpha+\theta)$ and hence
descends to give a map $\phi_0:C_N\to X$ such that
$$\tphi''_0(t,\theta)=\phi_0(qt,q\theta).$$
Define also $\psi:[-N,N]\to X^{\lambda_{\crit}}$ by the condition
$\psi''(t)=\psi(qt)$. Then the equality $\phi''=\exp_{\psi''}\phi_0''$
translates into
$$\phi(qt,q\theta)=e^{-\imag\theta p}\exp_{\psi(qt)}(e^{\imag\theta p}
\phi_0(qt,q\theta)),$$
which is equivalent to the first formula in
(\ref{eq:definicio-phi0}). On the other hand, the balancing condition
$\int \phi_0''(t,\nu)d\nu=0$ is clearly equivalent to
$\int_0^{2\pi q} e^{\lambda_{\crit}\nu}\phi_0(t,\nu)d\nu=0$.

Finally, since the set of $q$ which
appear as denominators of numbers in $\imag\Lambda_{\crit}$
is finite, $\epsilon$ and $\eta_0$ can be chosen in such a way that for any
critical residue we can reduce to the case of zero residue applying
the procedure which we just described.

\subsection{Proof of Theorem \ref{thm:psi-gradient}}
\label{ss:proof-thm:psi-gradient}

As in the proof of Theorem \ref{thm:small-energy}, we prove the
case $\lambda_{\crit}=0$. Using the same covering argument
as in Section \ref{ss:proof-thm:small-energy} the general case
is reduced to this case.
Furthermore, we only prove that $\iota_{v}d\alpha+\mu(\phi)=c$
implies (\ref{eq:der-psi}). That $d\alpha=0$ implies
(\ref{eq:der-psi-flat}) is even simpler and follows from the same ideas.

We begin by stating a local version of the theorem.

\begin{lemma}
Let $K_0>0$ be a real number,
let $I$ (resp. $g$, $\XXX$) be an almost complex structure (resp.
Riemannian metric, vector field) on $\CC^n$, and let
$V\subset\CC^n$ be a compact subset.
Let $Z:=[-1,2]\times S^1$, and let $\phi:Z\to V$ be a map satisfying
\begin{equation}
\label{eq:l'equacio-de-sempre}
\frac{\partial\phi}{\partial t}=
I(\phi)\left(\frac{\partial \phi}{\partial \theta}
-\imag(\lambda+\beta)\XXX(\phi)\right).
\end{equation}
Suppose that the diameter of $\phi(Z)$ is small enough so that
$\psi:[0,1]\to\CC^n$ and $\phi_0:Z\to\CC^n$ can be defined
as in Theorem \ref{thm:small-energy} (here we identify
$T_{\psi(t)}\CC^n\simeq\CC^n$) and that
\begin{equation}
\label{eq:phi-balla-poc}
|\phi_0|<K_0
\qquad\text{and}\qquad
\left|\frac{\partial\phi_0}{\partial \theta}\right|,
\left|\frac{\partial\phi_0}{\partial t}\right|<K_0.
\end{equation}
Then for any $t\in(0,1)$ we have
\begin{equation}
\label{eq:psi-quasi-gradient-local}
|\psi'(t)+\imag\lambda I(\psi(t))\XXX(\psi(t))|<K\sup_{Z_t}(|\phi_0|+|\beta|)
\end{equation}
for some constant $K$ independent of $\phi$, where $Z_t:=\{t\}\times S^1$.
\end{lemma}
\begin{pf}
Let $\tV:=\{\exp_x^gv\mid x\in V, |v|\leq K_0\}$.
By hypothesis we have $\psi([0,1])\subset\tV$.
Let $B\subset\CC^n$ be the closed ball of radius $K_0$,
let $E:\tV\times B\to\CC^n$ be the exponential map
$E(x,v):=\exp_x^gv$. Consider the first derivatives $F:=D_xE$
and $G:=D_vE$. Since the domain of $E$ is compact, the second
derivatives of $E$ are uniformly bounded; combining this observation
with the fact that $F(x,0)=\Id$ and $G(x,0)=\Id$ for any $x$, we obtain
bounds
\begin{equation}
\label{eq:bounds-F-G}
|F(x,v)-\Id|<K|v|
\qquad\text{and}\qquad
|G(x,v)-\Id|<K|v|.
\end{equation}
Taking partial derivatives with respect to $t$
in the equality $\exp_{\psi}\phi_0=\phi$ we obtain
$$F(\psi,\phi_0)\psi'+G(\psi,\phi_0)\frac{\partial \phi_0}{\partial t}
=\frac{\partial\phi}{\partial t}.$$
Writing $F=\Id+(F-\Id)$ and $G=\Id+(G-\Id)$ and integrating over $S^1$
the above equality, combined with
(\ref{eq:phi-balla-poc}), (\ref{eq:bounds-F-G}), and the fact
that $\int\phi_0(t,\theta)d\theta=0$ for every $t$, yields
\begin{equation}
\label{eq:derivada-psi}
\left|\psi'(t)-\frac{1}{2\pi}\int\frac{\partial \phi}{\partial t}\right|
<K\sup_{Z_t}|\phi_0|.
\end{equation}
Using (\ref{eq:l'equacio-de-sempre}) we compute
\begin{align*}
\frac{\partial \phi}{\partial t} =
& -\imag\lambda I(\psi)\XXX(\psi)
+I(\psi)\frac{\partial\phi}{\partial\theta}
+(I(\phi)-I(\psi))\left(\frac{\partial\phi}{\partial\theta}-\imag\lambda
\XXX(\phi)\right) \\
&-\imag\lambda I(\psi)(\XXX(\phi)-\XXX(\psi))
-\imag\beta I(\phi)\XXX(\phi).
\end{align*}
Integrating over $S^1$ and dividing by $2\pi$ the first
term in the right hand side does not change; the second term vanishes,
because $I(\psi)$ is independent of $\theta$ and the integral of
$\partial\phi/\partial\theta$ over $S^1$ is obviously $0$;
the third and fourth term
can be bounded by $K\sup_{Z_t}|\phi_0|$, since both
$|I(\psi)-I(\phi)|$
and $|\XXX(\psi)-\XXX(\phi)|$ are less than $K\sup_{Z_t}|\phi_0|$
(of course here $K$ depends on the derivatives of $I$ and $\XXX$)
and (\ref{eq:phi-balla-poc}) holds; the fifth term can be bounded
by $K\sup_{Z_t}|\beta|$.
Putting together all this observations and using
(\ref{eq:derivada-psi}) we obtain inequality
(\ref{eq:psi-quasi-gradient-local}).
\end{pf}

Covering $X$ with a finite number of charts we deduce from the lemma
that for any $t\in[-N,N]$:
\begin{equation}
\label{eq:psi-quasi-gradient-X}
|\psi'(t)+\imag\lambda I(\psi(t))\XXX(\psi(t))|
<K\sup_{\{t\}\times S^1}(|\phi_0|+|\beta|).
\end{equation}
Note in fact that the constant $K_0$ in the Lemma
can be taken independently of $(\alpha,\phi)$, and depending only
on the $\epsilon$ in Theorem \ref{thm:small-energy},
because (\ref{eq:condicions-pair}) implies the bounds
(\ref{eq:phi-balla-poc}). Combining this with the fact that the
number of charts which we use is finite, we deduce that the constant
$K$ in the above inequality is independent also of $(\alpha,\phi)$.

To finish the proof of Theorem \ref{thm:psi-gradient} we need
to bound the right hand side of (\ref{eq:psi-quasi-gradient-X}).
On the one hand, formula (\ref{eq:estimate-phi0})
in Theorem \ref{thm:small-energy} tells us that
\begin{equation}
\label{eq:cota-psi-phi0}
|\phi_0(t,\theta)|< K e^{\sigma(|t|-N)}
\left( \|d_{\alpha}\phi\|_{L^{\infty}}+\|\alpha\|_{L^{\infty}}
+\eta(|c|+\sup|\mu|)\right)^{1/4}.
\end{equation}
On the other hand, we can bound $\beta$ writing
$$\beta(t,\theta)=\beta(0,\theta)+\int_0^t\frac{\partial \beta}{\partial t}
(\tau,\theta)d\tau=
\int_0^t\frac{\partial \alpha_{\theta}}{\partial t}(\tau,\theta)d\tau,$$
then using
$d\alpha=\frac{\partial\alpha_{\theta}}{\partial t}dt\wedge d\theta$
and (\ref{eq:dalpha-exp}) to obtain
\begin{align}
|\beta(t,\theta)| &\leq \int_0^t |d\alpha(\tau,\theta)|d\tau
\leq \int_0^t e^{|\tau|-N} \eta (|c|+\sup|\mu|)d\tau \notag \\
&\leq e^{|\tau|-N} \eta (|c|+\sup|\mu|)
\leq e^{\sigma(|t|-N)} \eta (|c|+\sup|\mu|) \notag \\
&\leq K e^{\sigma(|t|-N)} (\eta (|c|+\sup|\mu|))^{1/4}
\label{eq:cota-psi-beta}
\end{align}
where we have used that $\sigma<1$
and the fact that $\eta (|c|+\sup|\mu|)$ is bounded above.
Combining the estimates (\ref{eq:cota-psi-phi0}) with
(\ref{eq:cota-psi-beta}) with
(\ref{eq:psi-quasi-gradient-X}) we obtain the desired bound.

\subsection{Proof of Theorem \ref{thm:small-energy-average}}
\label{proof:thm:small-energy-average}

Using Lemma \ref{lemma:equi-charts} we reduce the proof of
Theorem \ref{thm:small-energy-average} to a local version
of it, very much as we did in the proof of Theorem
\ref{thm:ineq-perturb}.
We use the same notation as in Lemma \ref{lemma:ineq-energia},
so we fix a diagonal action of $S^1$ on $\CC^n$
with weights $w=(w_1,\dots,w_n)\in\ZZ^n$. Let $I_0$ (resp. $g_0$)
denote the standard complex structure (resp. Riemannian metric)
on $\CC^n$.

\begin{lemma}
\label{lemma:small-energy-average}
For any $r>0$ and $K>0$
there exists some $\epsilon=\epsilon(K,w,r)>0$ with the following property.
Suppose that $I$ (resp. $g$)
is a smooth equivariant almost complex structure (resp. Riemannian
metric) on $\CC^n$ such that $\|I-I_0\|_{L^\infty}<\epsilon$,
$\|D I\|_{L^{\infty}}<K$, and
$\|D^2I\|_{L^{\infty}}<K$ (here $DI$ and $D^2I$ denote the first and
second derivatives of $I$ and the norms are taken with respect to
the standard metric $g_0$ on $\CC^n$), and
$\|g-g_0\|_{L^{\infty}}<\epsilon$
and $\|Dg\|_{L^{\infty}}<K$.
Let $(\alpha,\phi):C\times S^1\to \CC^n$ be a $I$-holomorphic pair
satisfying $\|\alpha\|_{L^{\infty}}<\epsilon$,
$\|d_{\alpha}\phi\|_{L^{\infty}}<\epsilon$ and
$$\phi(C)\subset B(0,2r)\subset\CC^n.$$
Let $\psi:[-2,3]\to \CC^n$ and $\phi_0:C\to T\CC^n$ be defined as in
Theorem \ref{thm:small-energy}.
\begin{itemize}
\item[(1)] Suppose that
$\|d\alpha\|_{L^2}<\epsilon\|\phi_0\|_{L^2}$. Then
\begin{equation}
\label{eq:average-norma-phi0-S1} \|\phi_0\|^2_{L^2(\ZII)}\leq
\frac{1}{e+e^{-1}} \left( \|\phi_0\|^2_{L^2(\ZI)}+
\|\phi_0\|^2_{L^2(\ZIII)}\right)
\end{equation}
\item[(2)] Suppose now that $\|d\alpha\|_{L^2_1}<K$ and
$\|d\alpha\|_{L^{\infty}}<K.$ Then $\sup_Z|\phi_0|\leq
K\|\phi_0\|_{L^2}^{1/4}.$
\end{itemize}
\end{lemma}

The proof of the lemma 
will be given in Section \ref{ss:proof-lemma:small-energy-average}
below.

\begin{remark}
We specify the dependence of $\epsilon$ on $DI$, $D^2I$ and $Dg$
because we want to be sure that the hypothesis of the lemma are
preserved when we zoom in (this will become clear in the course of
the proof). The same comment applies for Lemma
(\ref{lemma:exponencial-average}).
\end{remark}

\subsection{Comparison between $\phi_0$ and $\phi_{\av}$}

When $(\alpha,\phi)$ is a pair taking values in $\CC^n$, the map
$\phi_0$ defined in Theorem \ref{thm:small-energy-average} can be
roughly speaking
approximated by the map $\phi_{\av}$ defined in (2) of Lemma
\ref{lemma:convex-hol}. The following lemma makes this statement
precise, specifying to what extent $\phi_{\av}$ is a good
approximation of $\phi_0$, both in $L^2$ and $C^0$ norms. This
will be crucial in proving Lemma \ref{lemma:small-energy-average},
since it is clearly easier to deal with $\phi_{\av}$ than with
$\phi_0$.

\begin{lemma}
\label{lemma:exponencial-average} For any $K>0$ there exist
constants $K_0>0$ and $\epsilon>0$ with the following property.
Denote by $g_0$ the standard Riemannian flat metric in $\CC^n$.
Let $g$ be another metric on $\CC^n$ such that $K^{-1}g_0\leq
g\leq K g_0$ and such that $\|Dg\|_{L^{\infty}}<K$, where $Dg$
denotes the first derivatives of $g$ and the norm is taken with
respect to $g_0$. Let $x\in M$ and let $\gamma_0:S^1\to\CC^n$ be a
smooth map satisfying $\int\gamma_0=0$ and $\sup
|\gamma_0|<\epsilon$. Let
$\gamma(\theta):=\exp^g_x\gamma_0(\theta)$ (here we are
identifying $T_x\CC^n\simeq\CC^n$) and define
$$\gamma_{\av}(\theta):=\gamma(\theta)-
\frac{1}{2\pi}\int \gamma(\nu)d\nu.$$ Then we have
\begin{equation}
\label{eq:exp-ave-L2} \|\gamma_0\|_{L^2}(1-K_0\|\gamma_0\|_{L^2})
\leq \|\gamma_{\av}\|_{L^2} \leq
\|\gamma_0\|_{L^2}(1+K_0\|\gamma_0\|_{L^2})
\end{equation}
and similarly
\begin{equation}
\label{eq:exp-ave-sup} \sup |\gamma_0|(1-K_0\sup|\gamma_0|) \leq
\sup|\gamma_{\av}| \leq \sup |\gamma_0| (1+K_0\sup |\gamma_0|).
\end{equation}
\end{lemma}
\begin{pf}
We claim that given $K>0$, $\epsilon>0$, and a metric $g$ on
$\CC^n$ satisfying $K^{-1}g_0\leq g\leq K g_0$ and such that
$\|Dg\|_{L^{\infty}}<K$, there is a constant $K'$ depending only
on $K$ such that for any $x,v\in\CC^n$ satisfying $|v|<\epsilon$
we have
\begin{equation}
\label{eq:derivada-exp} |\exp^g_xv-x-v|<K'|v|^2.
\end{equation}
To see this, define $\gamma(t):=\exp^g_x tv$. Since $\gamma(t)$ is
a geodesic, we have $\frac{d \gamma'_k}{d
t}=-\Gamma_{ij}^k\gamma'_i\gamma'_j$ (here
$\gamma_1,\dots,\gamma_{2n}$ are the components of $\gamma$).
Integrating and using $\gamma'(0)=v$ we deduce that (provided
$\epsilon>0$ has been chosen small enough and $|v|<\epsilon$) for
any $t\in[0,1]$, $|\gamma'(t)-v|<K'|v|^2$, where $K'$ is
proportional to the sup norm of $\Gamma_{ij}^k$ which, on its
turn, can be estimated in terms of $(g^{ij})=(g_{ij})^{-1}$ and
the derivatives of $(g_{ij})$. Integrating this inequality we
obtain (\ref{eq:derivada-exp}).

It follows from  that (\ref{eq:derivada-exp}) that
$$\left|x-\frac{1}{2\pi}\int\gamma\right|=
\left|x-\frac{1}{2\pi}\int\exp_x^g\gamma_0\right|
<K'\|\gamma_0\|_{L^2}^2.$$ Consequently we have
\begin{align*}
|\gamma_{\av}(\theta)-\gamma_0(\theta)| &=\left|
\gamma(\theta)-\left(\frac{1}{2\pi}\int\gamma\right) -
\gamma_0(\theta)\right| \\
 &<\left| x + \gamma_0(\theta)-
\left(\frac{1}{2\pi}\int\gamma\right) - \gamma_0(\theta)\right|+K'|\gamma_0(\theta)|^2 \\
&<K'(\|\gamma_0\|_{L^2}^2+|\gamma_0(\theta)|^2).
\end{align*}
Integrating over $\theta$ we obtain (\ref{eq:exp-ave-L2}).
Inequality (\ref{eq:exp-ave-sup}) is proved similarly.
\end{pf}

\subsection{Proof of Lemma \ref{lemma:small-energy-average}}
\label{ss:proof-lemma:small-energy-average}

We follow an idea similar to the proof of Lemma
\ref{lemma:ineq-energia} in Section
\ref{proof:lemma:ineq-energia}. Assume that there exists sequences
of positive real numbers $\epsilon_u\to 0$, invariant almost
complex structures $I_u$ and metrics $g_u$, and $I_u$-holomorphic
pairs $(\alpha_u,\phi_u)$. Suppose that $I_u$, $g_u$ and
$(\alpha_u,\phi_u)$ satisfy the hypothesis of the lemma for
$\epsilon=\epsilon_u$. In particular, we have
\begin{equation}
\label{eq:convergencia-a-phi} \|\alpha_u\|_{L^{\infty}}\to 0
\qquad\text{and}\qquad \|d_{\alpha_u}\phi_u\|_{L^{\infty}}\to 0.
\end{equation}
Finally, suppose that for each $u$ either (1) or (2) of the lemma
fails to be true. We will see that this is impossible.

As we did in Lemma \ref{lemma:ineq-energia}, we assume that $w_j$
is nonzero for any $j$ between $1$ and $p$, and that
$w_{p+1}=\dots=w_n=0$. Let $q=n-p$, so that we have a splitting
$\CC^n=\CC^p\times\CC^q$. Let also $\pi_p$ and $\pi_q$ denote the
projections from $\CC^n$ to $\CC^p$ and $\CC^q$.

Define $\rho_u:=\|d_{\alpha_u}\phi_u\|_{L^{\infty}}$ and
$s_u:=\sup |\pi_p \phi_u(C)|$. Note that since $\phi_u(C)\subset
B(0,2r)$ the numbers $s_u$ are uniformly bounded above. Let
$y_u\in C$ be a point where $|\pi_p\phi_u(C)|$ attains the value
$s_u$, and let $x_u=\phi_u(y_u)$. Passing to a subsequence if
necessary, we can assume that the sequence $\rho_u^{-1}s_u$
converges somewhere in $\RR_{\geq 0}\cup\{\infty\}$. We
distinguish two possibilities.

Suppose first that $\rho_u^{-1}s_u\to e<\infty$. Let
$\phi_u':=\rho_u^{-1}(\phi_u-x_u).$ Define also $I_u'$ (resp.
$g_u$) to be the pullback of $I_u$ (resp. $g_u$) under the
composition of the translation along $x_u$ with the homotecy of
ratio $\rho$. Then $I_u'$ and $g_u'$ are $S^1$-invariant and still
satisfy the bounds in the hypothesis of the lemma. Furthermore,
$(\alpha_u,\phi_u')$ is a $I_u'$-holomorphic pair. Finally, since
$\|d_{\alpha_u}\phi_u'\|_{L^{\infty}}=1$,
$\|\alpha\|_{L^{\infty}}\to 0$ and $|\phi_u'(y)|\leq e$, the image
$\phi_u'(C)$ is contained in a compact set independent of $u$.
Hence we may apply Lemma \ref{lemma:easy-compactness} and deduce
that, up to regauging, the pairs $(\alpha_u,\phi_u')$ converge in
$C^1$ norm to a $I_0$-holomorphic pair $(0,\phi)$. So
$\|d\phi\|_{L^{\infty}}=1$ (hence $\phi$ is not constant) and
$\ov{\partial}_{I_0}\phi=0$.

Let $\phi_{\av}$ be as defined in (2) of Lemma
\ref{lemma:convex-hol}. Since $\phi$ is holomorphic and not
constant, by Lemma \ref{lemma:convex-hol}  we have
$$\|\phi_{\av}\|_{L^2(\ZII)} \leq \frac{1}{e^2+e^{-2}}
\|\phi_{\av}\|_{L^2(\ZI)}+\|\phi_{\av}\|_{L^2(\ZIII)}).$$ On the
other hand, since the convergence $\phi_u'\to\phi$ is in $C^0$, it
follows that $\phi_{u,\av}'$ converges pointwise to $\phi_{\av}$.
This implies that for big enough $u$ we have
$$\|\phi'_{u,\av}\|_{L^2(\ZII)} <
\frac{1}{e^{1.5}+e^{-1.5}}
(\|\phi'_{u,\av}\|_{L^2(\ZI)}+\|\phi'_{u,\av}\|_{L^2(\ZIII)}).$$
Since $\phi'_u$ is related to $\phi_u$ by a translation and a
homotecy, it follows that the same inequality is satisfied by
$\phi_u$. Finally, (\ref{eq:exp-ave-L2}) in Lemma
\ref{lemma:exponencial-average} implies (using the fact that
$\|\phi_{u,0}\|_{L^2}$ converges to $0$) that for big enough $u$
$$\|\phi_{u,0}\|_{L^2(\ZII)} < \frac{1}{e+e^{-1}}
(\|\phi_{u,0}\|_{L^2(\ZI)}+\|\phi_{u,0}\|_{L^2(\ZIII)}).$$ Hence
(1) in Lemma \ref{lemma:small-energy-average} has to hold for big
enough $u$.

On the other hand, since $\phi_{\av}$ is holomorphic and satisfies
$\int\phi_{\av}(t,\nu)d\nu=0$, standard elliptic estimates imply
that, for some $K$ independent of $\phi$, we have
$$\sup_Z|\phi_{\av}|\leq K\|\phi_{\av}\|_{L^2}.$$ Again using the
fact that $\phi_{u,\av}'$ converges to $\phi_{\av}$ we deduce that
the same inequality holds for $\phi_{u,\av}$ provided $u$ is big
enough and maybe after increasing slightly $K$. Finally,
(\ref{eq:exp-ave-sup}) in Lemma \ref{lemma:exponencial-average}
implies a similar inequality $\sup_Z|\phi_{u,0}|\leq
K\|\phi_{u,0}\|_{L^2}.$ But since $\|\phi_{u,0}\|_{L^2}$ is
smaller than $1$ for big enough $u$, this implies
$$\sup_Z|\phi_{u,0}|\leq K\|\phi_{u,0}\|_{L^2}^{1/4}.$$
Consequently, (2) in Lemma \ref{lemma:small-energy-average} has to
hold for big enough $u$. And this is in contradiction with our
assumptions.

Now suppose that $\rho_u^{-1}s_u\to\infty$ and define
$\phi_u':=s_u^{-1}(\phi_u-x_u).$ Then
$\|d_{\alpha_u}\phi_u'\|_{L^{\infty}}\to 0$, so the diameter of
$\phi_u'$ converges to $0$. On the other hand each $\phi_u'$
intersects the set
$$S=\{(x,0)\in\CC^p\times\CC^q\mid |x|=1\}\subset\CC^n.$$
Then (1) in Lemma \ref{lemma:a-la-vora-de-la-vora} below implies
that for big enough $u$ the pairs have to satisfy (1) in Lemma
\ref{lemma:small-energy-average}. On the other hand, (2) in Lemma
\ref{lemma:a-la-vora-de-la-vora} implies that for big enough $u$
we have $\sup_Z|\phi'_{u,0}|\leq K\|\phi'_{u,0}\|_{L^2}^{1/4}$.
Since $\phi_{u,0}=s_u\phi'_{u,0}$ and $s_u$ is uniformly bounded
above, maybe after increasing $K$ we also have the following
inequality for big enough $u$:
$$\sup_Z|\phi_{u,0}|\leq
K\|\phi_{u,0}\|_{L^2}^{1/4}.$$ This implies that (2) in Lemma
\ref{lemma:small-energy-average} has to hold for big enough $u$,
leading to a contradiction and finishing the proof of the lemma.

\begin{lemma}
\label{lemma:a-la-vora-de-la-vora} For any $K>0$ there exist
numbers $\epsilon=\epsilon(K,w)>0$ and $\delta=\delta(K,w)>0$ with
the following property. Suppose that $I$ (resp. $g$) is a smooth
equivariant almost complex structure (resp. Riemannian metric) on
$\CC^n$ such that $\|I-I_0\|_{L^{\infty}}<\epsilon$, $\|D
I\|_{L^{\infty}}<K$, and $\|D^2I\|_{L^{\infty}}<K$ and
$\|g-g_0\|_{L^{\infty}}<\epsilon$ and $\|Dg\|_{L^{\infty}}<K.$ Let
$(\alpha,\phi):C\times S^1\to \CC^n$ be a $I$-holomorphic pair
satisfying $\|\alpha\|_{L^{\infty}}<\epsilon$,
$\|d_{\alpha}\phi\|_{L^{\infty}}<\epsilon$, and
$$\phi(C)\subset B(y,\delta)\subset\CC^n,$$
where $y$ is a point in $S$. Let $\psi:[-2,3]\to \CC^n$ and
$\phi_0:C\to T\CC^n$ be defined as in Theorem
\ref{thm:small-energy}.
\begin{itemize}
\item[(1)] Suppose that
$\|d\alpha\|_{L^2}<\epsilon\|\phi_0\|_{L^2}$. Then
$$\|\phi_0\|^2_{L^2(\ZII)}\leq
\frac{1}{e+e^{-1}} \left( \|\phi_0\|^2_{L^2(\ZI)}+
\|\phi_0\|^2_{L^2(\ZIII)}\right).$$ \item[(2)] There is a constant
$K'$ depending only on $K$ with the following property. Suppose
that $\|d\alpha\|_{L^2_1}<K$ and $\|d\alpha\|_{L^{\infty}}<K$.
Then $\sup_Z |\phi_0|\leq K'\|\phi_0\|_{L^2}^{1/4}.$
\end{itemize}
\end{lemma}

\subsection{Proof of Lemma \ref{lemma:a-la-vora-de-la-vora}}

We deduce the lemma from a local statement similar
to Lemma \ref{lemma:small-energy-average} but accounting
for the degenerate situation
in which instead of a linear action of $S^1$ on a vector space we have
an action of $\RR$ by translations.
So now $\XXX$ denotes a constant vector in $\CC^n$.

If $(\alpha,\phi):C\to\CC^n$ is a pair in which
$\alpha=\alpha_{\theta} d\theta$, its energy density is defined by
$$|d_{\alpha}\phi|^2=\left|\frac{\partial\phi}{\partial t}\right|^2
+\left|\frac{\partial\phi}{\partial \theta}
-\imag\alpha_{\theta}\XXX\right|^2.$$ Furthermore, the condition
on $(\alpha,\phi)$ of being $I$-holomorphic reads
\begin{equation}
\label{eq:parell-holomorf-dege} \frac{\partial \phi}{\partial t}
=I(\phi)\left(\frac{\partial \phi}{\partial \theta}
-\imag\alpha_{\theta}\XXX\right).
\end{equation}

\begin{lemma}
\label{lemma:small-energy-average-boundary} For any $K>0$ there
exists a constant $\epsilon>0$ with the following property. Let
$I$ be a translation invariant almost complex structure on $\CC^n$
satisfying
$$\|I-I_0\|_{L^{\infty}}<\epsilon,\qquad
\|D I\|_{L^{\infty}}<K,\qquad\text{and}\qquad \|D^2
I\|_{L^{\infty}}<K$$ and let $g$ be a Riemannian metric on $\CC^n$
such that $\|g-g_0\|_{L^{\infty}}<\epsilon$ and
$\|Dg\|_{L^{\infty}}<K.$ Take a $I$-holomorphic pair
$(\alpha,\phi):C\to\CC^n$ satisfying
\begin{equation}
\label{eq:cotes-alpha-phi}
 \|\alpha\|_{L^{\infty}}<\epsilon,
\qquad\text{and}\qquad \|d_{\alpha}\phi\|_{L^{\infty}}<\epsilon,
\end{equation}
and define $\psi:[-2,3]\to\CC^n$ and $\phi_0:C\to T\CC^n$ as in
Theorem \ref{thm:small-energy}.
\begin{itemize}
\item[(1)] Suppose that
$\|d\alpha\|_{L^2}<\epsilon\|\phi_0\|_{L^2}.$ Then
\begin{equation}
\label{eq:average-norma-phi0-trans}
\|\phi_0\|^2_{L^2(\ZII)}\leq
\frac{1}{e+e^{-1}}
\left(
\|\phi_0\|^2_{L^2(\ZI)}+
\|\phi_0\|^2_{L^2(\ZIII)}\right).
\end{equation}
\item[(2)] There is a constant $K'$ depending only on $K$ with the
following property. Suppose that $\|d\alpha\|_{L^2_1}<K$ and
$\|d\alpha\|_{L^{\infty}}<K.$ Then we have
\begin{equation}
\label{eq:cota-diametre}
\sup_Z |\phi_0|\leq K'\|\phi_0\|_{L^2}^{1/4}.
\end{equation}
\end{itemize}
\end{lemma}

We now resume the proof of Lemma \ref{lemma:a-la-vora-de-la-vora}.
To avoid confusion, let $E$ denote $\CC^n$ with the linear action
of $S^1$ and let $F$ denote $\CC^n$ with an action of $\RR$
given by translations along a vector $\XXX\in F$.
We have the following analogue of the charts constructed
in Lemma \ref{lemma:equi-charts}.
Suppose that $I$ (resp. $g$) is any $S^1$-invariant
complex structure (resp. Riemannian metric) on $E$.
Let $\epsilon>0$ be as in Lemma \ref{lemma:small-energy-average-boundary}.
For any point $y\in S\subset E$
there exists a neighbourhood $U_x$ of $y$,
a translation invariant
complex structure $I_x$ and Riemannian metric $g_x$ on $F$
satisfying
$\|I_x-I_0\|_{L^{\infty}}<\epsilon$ and
$\|g_x-g_0\|_{L^{\infty}}<\epsilon$
(where $I_0$ and $g_0$ are the standard structures on $F$)
and such that $\|DI\|_{L^{\infty}}$, $\|D^2I\|_{L^{\infty}}$ and
$\|Dg\|_{L^{\infty}}$ are bounded,
and a complex isometry $\xi_x:U_x\to F$.
Then $\{U_x\}_{x\in S}$ cover $S$.
Consider a finite subcovering $\{U_1,\dots,U_l\}$.
For each open set $U_j$, let
$$K_j:=\max\{\|DI_j\|_{L^{\infty}},
\ \|D^2I_j\|_{L^{\infty}},
\ \|Dg_j\|_{L^{\infty}}\}.$$
Let $\epsilon_j$ be the value of $\epsilon$ given by Lemma
\ref{lemma:small-energy-average-boundary} for the choice $K=K_j$,
and set $\epsilon:=\min\{\epsilon_1,\dots,\epsilon_l\}$.
There exists some constant $\delta>0$ such that any ball
of radius $\delta$ centered somewhere in $S$ is contained
in one of the $U_j$'s. Taking this value of $\delta$,
the statement of Lemma \ref{lemma:a-la-vora-de-la-vora}
follows from considering one of the open sets
(say $U_j$) which contains the image of $\phi$,
taking the composition $\phi_j:=\xi_j\circ\phi$
and applying Lemma \ref{lemma:small-energy-average-boundary}
to $(\alpha,\phi_j)$.

\subsection{Proof of Lemma \ref{lemma:small-energy-average-boundary}}

\subsubsection{Proof of (1)}
We proceed by contradiction.
Assume that there exist sequences of
real numbers $\epsilon_u>0$, almost complex structures $I_u$,
Riemannian metrics $g_u$ on $\CC^n$ and $I_u$-holomorphic pairs
$(\alpha_u,\phi_u)$ with $\epsilon_u\to 0$,
\begin{equation}
\label{eq:estimates-I} \|I_u-I_0\|_{L^{\infty}}<\epsilon_u,\qquad
\|D I_u\|_{L^{\infty}}<K,\qquad \|D^2 I_u\|_{L^{\infty}}<K,
\end{equation}
\begin{equation}
\label{eq:estimates-g} \|g_u-g_0\|_{L^{\infty}}<\epsilon,\qquad
\|Dg_u\|_{L^{\infty}}<K,
\end{equation}
\begin{equation}
\label{eq:estimates-alpha-phi}
\|\alpha_u\|_{L^{\infty}}<\epsilon_u, \qquad
\|d_{\alpha_u}\phi_u\|_{L^{\infty}}<\epsilon_u,
\end{equation}
and such that defining $\psi_u$ (resp. $\phi_{u,0}$) as
$\psi$ (resp. $\phi_0$) in Theorem
\ref{thm:small-energy} we have the estimate
\begin{equation}
\label{eq:estimates-dalpha}
\|d\alpha_u\|_{L^2}<\epsilon_u\|\phi_{u,0}\|_{L^2},
\end{equation}
and the inequality oposite to (\ref{eq:average-norma-phi0-trans})
is satisfied
\begin{equation}
\label{eq:average-norma-phi0-oposite}
\|\phi_{u,0}\|^2_{L^2(\ZII)}>
\frac{1}{e+e^{-1}}
\left(
\|\phi_{u,0}\|^2_{L^2(\ZI)}+
\|\phi_{u,0}\|^2_{L^2(\ZIII)}\right).
\end{equation}
Define
$$\phi_{u,\av}(t,\theta):=\phi_u(t,\theta)-\frac{1}{2\pi}\int
\phi_u(t,\nu)d\nu.$$ It follows from
(\ref{eq:estimates-alpha-phi}) that $\|\phi_{u,0}\|_{L^\infty}\to
0$ (unless we specify something different, all limits in this
proof will implicitly mean as $u$ goes to $\infty$). Combining
this fact with Lemma \ref{lemma:exponencial-average} and
(\ref{eq:average-norma-phi0-oposite}) we deduce that for big
enough $u$ we also have
\begin{equation}
\label{eq:average-norma-phi0-oposite-app}
\|\phi_{u,\av}\|^2_{L^2(\ZII)}>
\frac{1}{e^{1.5}+e^{-1.5}} \left(
\|\phi_{u,\av}\|^2_{L^2(\ZI)}+
\|\phi_{u,\av}\|^2_{L^2(\ZIII)}\right).
\end{equation}
Lemma \ref{lemma:exponencial-average} combined with
(\ref{eq:estimates-dalpha}) also implies that
\begin{equation}
\label{eq:estimates-dalpha-app}
\|d\alpha_u\|_{L^2}/\|\phi_{u,\av}\|_{L^2}\to 0.
\end{equation}

Let
$$x_u:=\frac{1}{2\pi}\int \phi_u(0,\nu) d\nu.$$
Using a real linear transformation $\CC^n\to\CC^n$ we may assume
that $I_u(x_u)=I_0$ and $g_u(x_u)=g_0$. Furthermore, since
$\|I_u-I_0\|_{L^\infty}\to 0$, these linear transformations are
uniformly bounded and (maybe after increasing $K$ slightly) they
preserve the bounds on $DI$, $D^2I$ and $Dg$.

Using a gauge transformation we can assume that $\alpha_u$ in
balanced temporal gauge, so that $\alpha_u=a_u d\theta$ for some
function $a_u:C\to\imag\RR$ and the restriction of $\alpha_u$ to
$\{0\}\times S^1$ is equal to a constant $\lambda_u\in\imag\RR$.
Then we have $\alpha_u=\lambda_u d\theta+\beta_u d\theta$, where
\begin{equation}
\label{eq:alpha-temporal-gauge}
|\lambda_u|<\|\alpha_u\|_{L^{\infty}}\qquad\text{and}\qquad
\|\beta_u\|_{L^\infty}\leq 2\|d\alpha_u\|_{L^\infty}.
\end{equation}

Let us define
$\xi_u(t,\theta):=\phi_u(t,\theta)+\imag\lambda_u t I_0\XXX-x_u$
and
\begin{equation}
\label{eq:def-xi-app}
\xi_{u,\av}(t,\theta):=\xi_u(t,\theta)-\frac{1}{2\pi}\int
\xi_u(t,\nu)d\nu.
\end{equation}
Then we have $\xi_{u,\av}=\phi_{u,\av}$. It follows that
$\|d\xi_u\|_{L^2}>0$, for otherwise we would have
$\xi_{u,\av}=\phi_{u,\av}=0$, in contradiction with
(\ref{eq:average-norma-phi0-oposite-app}).

\begin{lemma}
\label{lemma:xi-quasi-hol}
 The following holds: $\|\ov{\partial}_{I_0}\xi_u\|_{L^2}/\|d\xi_u\|_{L^2}\to
 0.$
\end{lemma}
\begin{pf}
Applying (\ref{eq:parell-holomorf-dege}) to the pair
$(\alpha_u,\phi_u)$ we compute:
\begin{align}
\ov{\partial}_{I_0}\xi_u &= \frac{\partial \xi_u}{\partial
t}-I_0\frac{\partial \xi_u}{\partial \theta}
=(I_u(\phi_u)-I_0)\left(\frac{\partial\phi_u}{\partial\theta}
-\imag\lambda_u\XXX\right) -\imag\beta_u I_u(\phi_u)\XXX \notag \\
&=(I_u(\xi_u+x_u)-I_u(x_u))\left(\frac{\partial\phi_u}{\partial\theta}
-\imag\lambda_u\XXX\right) -\imag\beta_u I_u(\xi_u+x_u)\XXX,
\label{eq:d-bar-xi}
\end{align}
where in the second equality we have used the fact that $I_u$ is
invariant under translations along $\XXX$. We claim that for some
constant $K'$ depending on $K$
\begin{equation}
\label{eq:cota-I}
\|I_u(\xi_u+x_u)-I_u(x_u)\|_{L^2}<K'\|d\xi_u\|_{L^2}.
\end{equation}
Indeed, using (\ref{eq:estimates-I}) we deduce that for any
$x\in\CC^n$ at distance less than $1$ from $x_u$ we have, for some
$K'$ depending on $K$,
$\|I_u(x)-I_u(x_u)-DI(x-x_u)\|<K'\|x-x_u\|^2.$ Applying this to
$x=\xi_u+x_u$ and using the fact that $\|d\xi_u\|_{L^{\infty}}$ is
uniformly bounded (which follows from
(\ref{eq:estimates-alpha-phi}) and implies that $\xi_u(C)$ stays
uniformly not too far from $x_u$) we get (\ref{eq:cota-I}). The
estimates (\ref{eq:estimates-alpha-phi}) imply that
$\|d\phi_u\|_{L^{\infty}}\to 0$, which on its turn implies
\begin{equation}
\label{eq:cota-dphi} \left\|\frac{\partial\phi_u}{\partial
\theta}\right\|_{L^{\infty}}\to 0.
\end{equation}
Also, combining (\ref{eq:estimates-alpha-phi}) with
(\ref{eq:alpha-temporal-gauge}) we deduce that
\begin{equation}
\label{eq:cota-lambda} \lambda_u\to 0.
\end{equation}
Finally, combining (\ref{eq:estimates-dalpha-app}) with
$\|\phi_{u,\av}\|_{L^2}=\|\xi_{u,\av}\|_{L^2}<K'\|d\xi_u\|_{L^2}$
we deduce
\begin{equation}
\label{eq:cota-beta} \beta_u/\|d\xi_u\|_{L^2}\to 0.
\end{equation}
Combining (\ref{eq:d-bar-xi}) with
(\ref{eq:cota-I})---(\ref{eq:cota-beta}), we obtain the desired
limit.
\end{pf}

Now define $\xi_u':=\xi_u/\|d\xi_u\|_{L^2}$.
Since $\int \xi'_u(0,\nu)d\nu=0$
and $\|d\xi_u'\|_{L^2}=1$, the $L^2$ norm of $\xi_u'$ is uniformly
bounded above, so there is a constant $K_0$ such that
$$1\leq \|\xi'_u\|_{L^2_1}\leq K_0.$$
Since the inclusion $L^2_1\subset L^2$ is compact, it follows
that, up passing to a subsequence, there is some nonzero $\xi\in
L^2(C,\CC^n)$ such that $\xi'_u\to\xi$ in $L^2$. We claim that
$\xi$ is in fact holomorphic in the interior of $C$. Indeed, for
any test function $g$ supported in the interior of $C$ we have
$\la \xi,\ov{\partial}^*g\ra_{L^2} =\lim \la
\xi'_u,\ov{\partial}^*g\ra_{L^2} =\lim \la
\ov{\partial}\xi'_u,g\ra_{L^2}=0.$ Hence $\xi$ is a weak solution
of $\ov{\partial}=0$, and from standard regularity results we
deduce that $\xi$ is smooth and a strong solution:
$\ov{\partial}\xi=0$. Combining this with G\r{a}rding's inequality
(in the interior of $C$): $\|\xi'_u-\xi\|_{L^2_1}\leq
K'(\|\ov{\partial}(\xi'_u-\xi)\|_{L^2} + \|\xi'_u-\xi\|_{L^2})$ we
deduce that $\|\xi'_u-\xi\|_{L^2_1}\to 0$, so $\xi'_u$ converges
to $\xi$ in $L^2_1$ and hence
\begin{equation}
\label{eq:xi-no-constant} \|d\xi\|_{L^2}=1.
\end{equation}
Define $\xi'_{u,\av}$ and $\xi_{\av}$ exactly as
we defined $\xi_{u,\av}$ in (\ref{eq:def-xi-app}). Since
$\xi'_u\to\xi$ in $L^2$, it follows that $\xi'_{u,\av}\to
\xi_{\av}$ in $L^2$.  On the other hand, since
$\xi'_{u,\av}$ is a rescaling of $\xi_{u,\av}$ and we
have $\xi_{u,\av}=\phi_{u,\av}$, the formula
(\ref{eq:average-norma-phi0-oposite-app})
implies, passing to the limit, that
\begin{equation}
\label{eq:average-norma-xi-oposite-app}
\|\xi_{\av}\|^2_{L^2(\ZII)}\geq \frac{1}{e^{1.5}+e^{-1.5}}
\left( \|\xi_{\av}\|^2_{L^2(\ZI)}+
\|\xi_{\av}\|^2_{L^2(\ZIII)}\right).
\end{equation}
This contradicts (2) in Lemma \ref{lemma:convex-hol} (applied to
$\phi=\xi$) unless $\xi$ is constant, which is impossible by
(\ref{eq:xi-no-constant}). This finishes the proof of (1).

\subsubsection{Proof of (2)}
We follow the same scheme of proof as
before, assuming the existence of sequences $I_u$, $g_u$,
$I_u$-holomorphic pairs $(\alpha_u,\phi_u)$, and $\epsilon_u\to 0$
satisfying (\ref{eq:estimates-I})---(\ref{eq:estimates-alpha-phi}) and
\begin{equation}
\label{eq:estimate-dalpha-L21}
\|d\alpha\|_{L^2_1}<K,\qquad
\|d\alpha\|_{L^{\infty}}<K,
\end{equation}
but contradicting (\ref{eq:cota-diametre}). We will see that this
is not possible. Define $\xi_u$ and $\xi_{u,\av}$ as before. We
state two lemmae whose proof will be given below so as not to
break the argument. The proof of the Lemma \ref{lemma:xi-xi-app}
will be given in Section \ref{ss:xi-xi-app} and that of Lemma
\ref{lemma:norma-L-2-2} in Section \ref{ss:norma-L-2-2}.

\begin{lemma}
\label{lemma:xi-xi-app} For big enough $u$ we have
$\|\xi_u\|_{L^2} <2\|\xi_{u,\av}\|_{L^2}
=2\|\phi_{u,\av}\|_{L^2}.$
\end{lemma}

\begin{lemma}
\label{lemma:norma-L-2-2} There is some constant $K_0$ such that,
for big enough $u$, $\|\xi_u\|_{L^2_2}\leq K_0.$
\end{lemma}

Let $\|\xi_u\|_{1/2}$ denote the H\"older
$C^{0+\frac{1}{2}}$ norm of $\xi_u$. Since in real dimension
$2$ the Sobolev space $L^2_2$ is included in $C^{0+\frac{1}{2}}$,
Lemma \ref{lemma:norma-L-2-2} gives a uniform bound
(here $K_0$ may increase from one line to the other, but
will always be independent of $u$)
$$\|\xi_u\|_{1/2}<K_0.$$
Taking this into account and
applying Lemma \ref{lemma:Holder-cont} to $\xi_u$ we obtain
for big enough $u$
\begin{equation}
\label{eq:sup-L2}
\sup_Z |\xi_u| \leq K_0(\|\xi_u\|_{L^2}^{1/4}+\|\xi_u\|_{L^2})
\leq K_0\|\xi_u\|_{L^2}^{1/4},
\end{equation}
(recall that $\|\xi_u\|_{L^2}$ goes to $0$). It is easy to check
that $\sup_Z|\xi_{u,\av}|\leq K_0\sup_Z|\xi_u|.$ Combining this
observation with the equality $\xi_{u,\av}=\phi_{u,\av}$ and
(\ref{eq:exp-ave-sup}) in Lemma (\ref{lemma:exponencial-average})
we deduce
\begin{equation}
\label{eq:sup-L2-2}
\sup_Z|\phi_{u,0}|\leq K_0\sup_Z|\xi_u|.
\end{equation}
On the other hand, combining Lemma \ref{lemma:xi-xi-app} with
(\ref{eq:exp-ave-L2}) in Lemma \ref{lemma:exponencial-average} we
conclude
\begin{equation}
\label{eq:sup-L2-3}
\|\xi_u\|_{L^2}\leq K_0\|\phi_{u,0}\|_{L^2}.
\end{equation}
Putting together (\ref{eq:sup-L2})---(\ref{eq:sup-L2-3}) we deduce
$\sup_Z|\phi_{u,0}|\leq K_0\|\phi_{u,0}\|_{L^2}^{1/4},$ so
(\ref{eq:cota-diametre}) must hold, in contradiction with our
assumption. This finishes the proof of (2).

\begin{lemma}
\label{lemma:Holder-cont}
There is a constant
$K_0>0$ with the following property. Let $f:C\to\CC^n$ be a smooth
map and let $\|f\|_{1/2}$ denote its H\"older
$C^{0+\frac{1}{2}}$ norm. Then
\begin{equation}
\label{eq:Holder-cont} \sup_Z |f|\leq K_0
(\|f\|_{L^2}^{1/4}\|f\|_{1/2}^{4/5}+\|f\|_{L^2}).
\end{equation}
\end{lemma}
\begin{pf} Denote for convenience $H:=\|f\|_{1/2}$, and
let $x\in Z$ be such that $|f(x)|=\sup |f|$. Let
$$\rho=\frac{|f(x)|^2}{4H^2}.$$
Then for any $y\in C$ such that
$|x-y|\leq\rho$ we have $|f(y)|\geq |f(x)|/2$, because
$$\frac{|f(x)-f(y)|}{\rho^{1/2}}
\leq\frac{|f(x)-f(y)|}{|x-y|^{1/2}}\leq H
\ \Longrightarrow
\ |f(y)|\geq |f(x)|-\rho^{1/2}H \geq |f(x)|-\frac{|f(x)|}{2}.$$
Now suppose that $\rho\leq 1$, so that the ball $B$ of radius
$\rho$ centered at $x$ is contained in $C$. Then $\|f\|_{L^2}\geq
\|f\|_{L^2(B)}\geq K_0\rho^2|f(x)|,$ which implies, rearranging,
that
\begin{equation}
\label{eq:des-1} |f(x)|\leq K_0 \|f\|_{L^2}^{1/5}
\|f\|_{1/2}^{4/5}.
\end{equation}
On the other hand, if $\rho>1$ we estimate $|f(x)|\leq
K_0\|f\|_{L^2}.$ Summing the two inequalities we obtain
(\ref{eq:Holder-cont}).
\end{pf}

\subsection{Proof of Lemma \ref{lemma:xi-xi-app}}
\label{ss:xi-xi-app} Pick some big $u$ and
let $f:=\xi_u$ and $f_{\av}:=\xi_{u,\av}$. By definition we have
\begin{equation}
\label{eq:diferencia-f-f0}
 |f(t,\theta)-f_{\av}(t,\theta)|=
\frac{1}{2\pi}\int f(t,\nu)d\nu= \frac{1}{2\pi}\int_0^t\int
\frac{\partial f}{\partial t}(\tau,\nu)d\nu d\tau,
\end{equation}
since $\int f(0,\nu)d\nu=0$.
On the other hand, denoting $\ov{\partial}_{I_0}$ by
$\ov{\partial}$ we have
$$\frac{\partial f}{\partial t}(t,\nu)=I_0\frac{\partial
f}{\partial\theta}(t,\nu) +\ov{\partial}f(t,\nu).$$
Integrating for $\theta\in S^1$ the first term in the right hand
side vanishes, so we obtain
\begin{equation}
\label{eq:norma-L2-df} \left|\int_0^t\int \frac{\partial
f}{\partial t}(\tau,\nu)d\nu d\tau\right|\leq
K\|\ov{\partial}f\|_{L^2}
\end{equation}
for some constant $K$. Let $K'$ be the constant in G\r{a}rding's inequality
$$\|df\|_{L^2}\leq \|f\|_{L^2_1}\leq
K'(\|\ov{\partial}f\|_{L^2}+\|f\|_{L^2}).$$ Using Lemma
\ref{lemma:xi-quasi-hol} we know that, if $u$ is big enough,
$\|\ov{\partial}f\|_{L^2}<1/2K'\|df\|_{L^2}$. Rearranging the
terms in the inequality, this implies that
$\|df\|_{L^2}<2K'\|f\|_{L^2}$. Using again Lemma
\ref{lemma:xi-quasi-hol} we conclude that for big enough $u$ we
have $\|\ov{\partial}f\|_{L^2}<(2K\Vol(C))^{-1} \|f\|_{L^2}.$
Taking $L^2(C)$ norms in (\ref{eq:diferencia-f-f0}) and combining
the previous inequality with (\ref{eq:norma-L2-df}) we obtain
$$\|f\|_{L^2}\leq \|f_{\av}\|_{L^2}+\frac{1}{2}\|f\|_{L^2}.$$
Hence, $\|f\|_{L^2}\leq 2\|f_{\av}\|_{L^2}$, which is what we
wanted to prove.

\subsection{Proof of Lemma \ref{lemma:norma-L-2-2}}
\label{ss:norma-L-2-2}
As before, throughout this proof $K_0$ will denote a positive number which
may increase from line to line but which will always be independent of $u$.
We estimate the $L^2_1$ norm of
$\ov{\partial}_{I_0}\xi_u$ using formula (\ref{eq:d-bar-xi}).
First, since $DI$ is uniformly bounded and the diameter of
$\xi_u(C)$ tends to $0$ (because of
(\ref{eq:estimates-alpha-phi})), we have
$\|I_u(\xi_u+x_u)-I_u(x_u))\|_{L^{\infty}}\to 0$. Hence, for any
$\delta$ and big enough $u$ we have
\begin{align*}
\left\|(I_u(\xi_u+x_u)-I_u(x_u))\frac{\partial
\xi_u}{\partial\theta}\right\|_{L^2_1} &\leq
\|I_u(\xi_u+x_u)-I_u(x_u))\|_{L^{\infty}}\|\xi_u\|_{L^2_2} \\
&+ \|I_u(\xi_u+x_u)-I_u(x_u))\|_{L^2_1}\|\xi_u\|_{L^\infty} \\
&\leq \delta \|\xi_u\|_{L^2_2} +
\|I_u(\xi_u+x_u)-I_u(x_u))\|_{L^2_1}\|\xi_u\|_{L^\infty}.
\end{align*}
(Note by the way that
$\partial\xi_u/\partial\theta=\partial\phi_u/\partial\theta$.) Now
we can estimate
$$\|I_u(\xi_u+x_u)-I_u(x_u))\|_{L^2_1}\leq
K_0\|\xi_u\|_{L^2_1}$$ We also know that if $u$ is big enough
$\|\xi_u\|_{L^2_1}\leq 2\|\xi_u\|_{L^2}$ (this follows from Lemma
\ref{lemma:xi-quasi-hol} and G\r{a}rding's inequality, see for
example the argument in Section \ref{ss:xi-xi-app}). Furthermore,
$\|\xi_u\|_{L^2}$ goes to $0$. Also, since
$\int\xi_u(0,\nu)d\nu=0$ and the diameter of $\xi_u(C)$ goes to
$0$, $\|\xi_u\|_{L^\infty}<K_0$ holds if $u$ is big enough.
Putting these observations together we obtain
\begin{equation}
\label{eq:primer-term-L21}
\left\|(I_u(\xi_u+x_u)-I_u(x_u))\frac{\partial
\xi_u}{\partial\theta}\right\|_{L^2_1}\leq \delta
\|\xi_u\|_{L^2_2} + K_0.
\end{equation}
For the next term in (\ref{eq:d-bar-xi}) the following is easy to
prove (if $u$ is big enough), taking into account the previous arguments:
\begin{equation}
\label{eq:segon-term-L21}
\|(I_u(\xi_u+x_u)-I_u(x_u))\imag\lambda_u\XXX\|_{L^2_1}\leq
K_0\|\xi_u\|_{L^2_1}\leq K_0.
\end{equation}
The remaining term satisfies, provided $u$ is big enough,
\begin{align}
\|\beta_u I_u(\xi_u)\XXX\|_{L^2_1} &\leq
K_0(\|\beta_u\|_{L^2_1}+\|\beta_u\|_{L^{\infty}}\|\xi_u\|_{L^2_1}) \notag \\
&\leq
K_0(\|d\alpha_u\|_{L^2_1}+K\|d\alpha_u\|_{L^{\infty}}\|\xi_u\|_{L^2_1})
\notag \\
&\leq K_0(\|d\alpha_u\|_{L^2_1}+\|\xi_u\|_{L^2_1})\leq K_0,
\label{eq:tercer-term-L21}
\end{align}
where here we use (\ref{eq:estimate-dalpha-L21}) and the
inequality $\|\beta_u\|_{L^2_1}\leq K_0\|d\alpha_u\|_{L^2_1}$,
which is easy to prove. Taking this into account and combining
(\ref{eq:primer-term-L21}), (\ref{eq:segon-term-L21}) and
(\ref{eq:tercer-term-L21}), formula (\ref{eq:d-bar-xi}) implies
that $\|\ov{\partial}_{I_0}\xi_u\|_{L^2_1}\leq
\delta\|\xi_u\|_{L^2_2}+K_0.$ Finally, using the standard
inequality
$$\|\xi_u\|_{L^2_2}\leq
K'(\|\ov{\partial}_{I_0}\xi_u\|_{L^2_1}+\|\xi_u\|_{L^2_1})$$ and
taking $\delta$ smaller than $1/(2K')$ we obtain the desired bound.

\section{Limits of approximate gradient lines}
\label{s:limits-gradient-lines}

Recall that $H=-\imag\mu$. We denote for convenience $V:=I\XXX$,
so that $V$ is the negative gradient of $H$.

\begin{theorem}
\label{thm:convergeix-cap-a-gradient} Let $\sigma>0$ be a real
number. Suppose that $\{\psi_u:T_u\to X, l_u, G_u\}$ is a sequence
of triples in which each $\psi_u$ is a smooth map with domain a
finite closed interval $T_u\subset\RR$, each $l_u$ is a nonzero
real number and each $G_u>0$ is a real number. Suppose that for
each $u$ and $t\in T_u$ we have
\begin{equation}
\label{eq:psi-aprox-grad} |\psi_u'(t)-l_u V(\psi_u(t))|\leq
G_ue^{-\sigma d(t,\partial T_u)}.
\end{equation}
Suppose also that $G_u\to 0$ and that $l_u\to 0$. Passing to a
subsequence, we can assume that $l_u |T_u|$ converges somewhere in
$\RR\cup\{\pm\infty\}$ (here $|T_u|$ denotes the length of $T_u$).
Then we have the following.
\begin{enumerate}
\item If $\lim l_u|T_u|=0$ then $\lim \diam \psi_u(T_u)=0$. \item
If $\lim l_u|T_u|\neq 0$,  define for big enough $u$ and for every
$t\in S_u$ the rescaled objects $S_u:=l_u T_u$ and
$f_u(t):=\psi_u(t/l_u)$. There is a subsequence of $\{f_u,S_u\}$
which converges to a chain of gradient segments $\TTT$ in $X$.
\end{enumerate}
\end{theorem}

\begin{pf}
The case $l_u|T_u|\to 0$ is obvious, so we consider the case $\lim
l_u|T_u|>0$. We begin modifying slightly the definition of $S_u$.
Taking $u$ big enough we can assume that $l_u<1$. Suppose that
$T_u=[a,b]$ and let $T_u':=[a,a-\ln l_u/\sigma]$ and
$T_u'':=[b+\ln l_u/\sigma,b]$. Using (\ref{eq:psi-aprox-grad}) we
can bound
\begin{align*}
\diam (\psi_u(T'_u)) & \leq\int_{T_u'} l_u|V| dt
+G_u\int_0^{-\ln l_u/\sigma}e^{-\sigma x}dx \\
&\leq l_u (\sup|V|) |\ln l_u|/\sigma +G_u(1-l_u)/\sigma\to 0,
\end{align*}
and similarly $\diam (\psi_u(T_u''))\to 0$. Consequently, if we define
\begin{equation}
\label{eq:bona-def-tu}
S_u:=[(a-\ln l_u/\sigma)/l_u,(b+\ln l_u/\sigma)/l_u]
\qquad \text{and}\qquad
f_u(t):=\psi(t/l_u),
\end{equation}
then the statement of the theorem is equivalent to saying that
the sequence $(f_u,S_u)$ defined by (\ref{eq:bona-def-tu}) has
a subsequence converging to a chain of gradient segments $\TTT$.
Furthermore, equation (\ref{eq:psi-aprox-grad}) implies that for
any $t\in S_u$ we have
\begin{equation}
\label{eq:psi-aprox-grad-t} |f_u'(t)-V(f_u(t))|\leq G_ue^{-\sigma
d(t,\partial S_u)/l_u} \leq G_u e^{-\sigma d(t,\partial S_u)}.
\end{equation}

\begin{lemma}
\label{lemma:H-monotone} For any connected component $F_0\subset
F$ of the fixed point set and any small enough $\delta>0$ there
exist numbers $0<d<\delta$, $\eta$ and $c>0$ with the following
property: let $f:[a,b]\to X$ be a map satisfying
$|f'(t)-V(f(t))|\leq \eta$ for every $t$ and, for some $a<t<t'<b$,
$$f(t)\in F_0^{d}\qquad\text{and}\qquad
f(t')\notin F_0^{2\delta}$$ (recall that $F_0^d$ and
$F_0^{2\delta}$ denote the $d$ and $2\delta$-neighbourhoods of
$F_0$ respectively). Then, for every $t''>t'+\delta$ we have
$$H(f(t''))\leq H(F_0)-c.$$
In particular, if $d$ is small enough the for every
$t''>t'+\delta$ we have $f(t'')\notin F_0^d$ (i.e., $f$ never
comes back to $F_0^d$ after time $t'+\delta$).
\end{lemma}

The proof of Lemma \ref{lemma:H-monotone} will be given in Section
\ref{ss:proof-lemma:H-monotone}. Now we continue with the proof of
Theorem \ref{thm:convergeix-cap-a-gradient}. Passing to a
subsequence, we can assume that there is some $K>0$ such that for
any connected component $F'\subset F$ either
$\lim_{u\to\infty}d(\psi_u(S_u),F')=0$ or $d(\psi_u(S_u),F')\geq
K$ for every $u$. Let $F_1,\dots,F_l\subset F$ be the connected
components which fall in the first case. We claim that the values
of $H$ in each of the components $F_1,\dots,F_l$ are all
different.

To prove this, suppose on the contrary that for some $i$ and $j$
we have $H=H(F_i)=H(F_j)$. Take some small $\delta$, and let
$0<d<\delta$, $\eta$ and $c$ be the numbers obtained by taking
$F_0:=F_i$ in Lemma \ref{lemma:H-monotone}. Let also $\alpha>0$ be
so small so that $H(F_j^\alpha)\subset [H-c/2,H-c+2]$. Then, if
$G_u$ is smaller than $\eta$, and for some $t$ we have $f_u(t)\in
F_i^d$, then for every $t'\geq t$ we have $f_u(t)\notin
F_j^\alpha$. Interchanging the roles of $F_i$ and $F_j$ we prove
in the same way that there are some $d',\alpha'$ such that if
$f_u(t)\in F_j^{d'}$ and $t'\geq t$ then $f_u(t')\notin
F_i^{\alpha'}$. So we cannot have simultaneously
$\lim_{u\to\infty}d(\psi_u(S_u),F_i)=0$ and
$\lim_{u\to\infty}d(\psi_u(S_u),F_j)=0$.

Hence we can suppose that $H(F_1)>H(F_2)>\dots> H(F_l)$. Now we
pick some very small $\delta$. Suppose that $u$ is so big that we
can apply Lemma \ref{lemma:H-monotone} to each of the components
$F_1,\dots,F_l$. Let $d_1,\dots,d_l$ be the numbers given by Lemma
\ref{lemma:H-monotone} for our choice of $\delta$. Then we define,
for every $1\leq j\leq l$,
$$E_{u,j}^{\delta}:=[\inf\{t\in S_u\mid f(t)\in F_j^{d_j}\},
\inf\{ t\in S_u\mid f(t)\notin F_k^{2\delta}\}+\delta].$$ Let
$E_u^{\delta}$ be the union $E_{u,1}^{\delta}\cup\dots\cup
E_{u_l}^{\delta}$ and let $T_u^{\delta}$ denote the closure of the
complementary $S_u\setminus E_u^{\delta}$.

Passing to a subsequence we can assume that for all $\delta$ and
$u$ all sets $T_u^{\delta}$ have the same number of connected
components: $T_u^{\delta}=T_{u,1}^{\delta}\cup\dots\cup
T_{u,p}^{\delta},$ where $p$ lies between $l-1$ and $l+1$, and
that the inequality $T_{u,1}^{\delta}\leq E_{u,1}^{\delta}$ holds
either for all $u$ or for none of them; if the inequality holds
and furthermore
$$\limsup_{\delta\to 0}\limsup_{u\to\infty}\diam
(f_u(T_{u,0}^{\delta}))=0,$$ then we remove $T_{u,1}^{\delta}$
from $T_u^{\delta}$ and attach it to $E_{u,1}^{\delta}$.
Similarly, we remove $T_{u,p}^{\delta}$ from $T_u^{\delta}$ if the
diameter of its images converges to $0$. After this operations we
end up with a new set $T_u^{\delta}=T_{u,1}^{\delta}\cup\dots\cup
T_{u,k}^{\delta}$ for each $u$ and $\delta$.

It follows from the definition of $E_{u,j}^{\delta}$ and from
Lemma \ref{lemma:H-monotone}, is that for any $t\in T_u^{\delta}$
the point $f(t)$ lies away from the set $X'=F_1^{d_1}\cup\dots\cup
F_l^{d_l}$. Since $\sup_{X\setminus X'} |V|>0$, it follows that if
$G_u$ is small enough then the size of $f'(t)$ is comparable to
that of $V$. As a consequence, for fixed $\delta$ the intervals
$T_{u,j}^{\delta}$ have bounded length. From this it is rather
straightforward to prove, passing to a subsequence, that
(\ref{eq:f_u-quasi-gradient}) holds and that there is a limiting
gradient segment $(x_j^{\delta},T_j^{\delta})$ to which the images
of $T_{u,j}^{\delta}$ converge.

Hence we only need to prove that the images of the intervals
$E_{u,j}^{\delta}$ accumulate near the fixed point component $F_j$
and that their diameters tend to $0$. The former is almost obvious
from the definition, whereas the latter is a bit more subtle and
follows from next lemma, using the fact that
$f_u(E_{u,j}^{\delta})\subset F_j^{2\delta}$.
\end{pf}

\begin{lemma}
\label{lemma:diametre-petit} There is a constant $K>0$ with the
following property. For any small enough $\delta>0$, any big
enough $u$ (depending on $\delta$), and any interval $E\subset
S_u$ such that $f_u(E)\subset F^{\delta}$, we have $\diam
(f_u(E))\leq \delta K.$
\end{lemma}
\begin{pf}
If $\delta$ is small enough we can assume that $f_u(E)$ lies
in the $\delta$-neighbourhood $F_0^{\delta}$ of a unique connected
component $F_0\subset F$. Let $N$ be the normal bundle of $F_0$.
Using the exponential map with respect to the $S^1$-equivariant
metric $g$, we can identify $F_0^{\delta}$ with a neighbourhood
$N^{\delta}$ of the zero section of $N$. Then we can pullback the
vector field $V$ to a vector field on $N^{\delta}$, which we
denote by the same symbol $V$. We also think of $f_u|_E$ as taking
values in $N^{\delta}$. Consider on $N$ the restriction of the
metric $g$. This gives an equivariant Euclidean metric on $N$.

The bundle $N$ carries a linear action of $S^1$, and we can split
$N=N_1\oplus\dots\oplus N_k$ in such a way that $S^1$ acts on
$N_j$ with weight $w_j\neq 0$. Let $V_0\in\Gamma(TN^{\vert})$ be
the vertican tangent field whose value at a vector
$x=(x_1,\dots,x_k)\in N$ is $V_0(x):=(w_1x_1,\dots,w_kx_k)$ (here
we are identifying $TN^{\vert}\simeq N$). It is well known that
$V_0$ approximates at first order $V$ near $F_0$ (which we view as
the zero section of $N$). More precisely, there is a constant $K$
such that for any $x\in N^{\delta}$ we have
\begin{equation}
\label{eq:V0-aproxima-V} |V_0(x)-V(x)|\leq K|x|^2.
\end{equation}

Let $d:E\to\RR$ be the function defined as $d(t):=|f_u(t)|^2$. We
want to prove that $d(t)$ decays exponentially as $t$ goes away
from the extremes of $E$. We will follow the same idea as in the
proof of formula (\ref{eq:estimate-phi0}) in Theorem
\ref{thm:small-energy}; however, in this situation the analysis
will be much simpler.

Shifting $E$ (and modifying accordingly $f_u$) and removing if
necessary a small interval of length $<2$ at the end of $E$, we
can assume that $E=[-L,L]$, where $L$ is a natural number (it is
clear that the truncation does not affect the estimate). Define
for every natural number $-L\leq n<L$ the energy
$d_n:=\int_n^{n+1}|f_u(t)| dt$. Define also
$g_n:=\int_n^{n+1}|f_u'(t)-V(f_u(t))|dt$.

\begin{lemma}
\label{lemma:recurrencia-dels-d} Let
$\gamma:=1/(e^{1/2}+e^{-1/2})$. If $\delta>0$ is small enough,
then there is some $\epsilon>0$ such that, for every $-L<n<L-1$
satisfying $g_{n-1}+g_n+g_{n+1}\leq\epsilon d_{n}$, we have
$$d_{n}\leq\gamma(d_{n-1}+d_{n+1}).$$
\end{lemma}

The proof of Lemma \ref{lemma:recurrencia-dels-d} will be given in
Section \ref{ss:proof-recurrencia-dels-d}. We now finish the proof
of Lemma \ref{lemma:diametre-petit}. Inequality
(\ref{eq:psi-aprox-grad-t}) implies that $g_n\leq
G_ue^{-\sigma(L-n)}.$ Arguing exactly as in the proof of Lemma
\ref{lemma:recur-perturb}, we deduce that for some $\beta>0$
independent of $f_u$ there is a bound $d_n\leq (\epsilon^{-1}K
G_u+d_{-L}+d_L)e^{-\beta(L-|n|)}.$ Then we have (using the fact
that $|V(x)|\leq K|x|$)
\begin{align*}
\diam(f_u(E)) &\leq K\int_E|f'_u|
\leq K\int_E(|V(f_u)|+|f'_u-V(f_u)|) \\
&\leq K\int |f_u|+K \int_E G_ue^{-\sigma(t,\partial S_u)}dt
\leq K\sum_{n=-L}^L d_n+K G_u \\
&\leq K(G_u+d_{-L}+d_L)\leq K(G_u+\delta).
\end{align*}
Taking $u$ big enough so that $G_u<\delta$, the result follows.
\end{pf}

\subsection{Proof of Lemma \ref{lemma:H-monotone}}
\label{ss:proof-lemma:H-monotone}

We can assume without loss of generality that $\eta\leq 1$. Define
$M:=\sup |V|+1$. Then for every $t$ we have
\begin{equation}
\label{eq:cota-derivada-f} |f'(t)|\leq M.
\end{equation}
Let $h:=H(f)$ and let $h_0:=H(F_0)$. We have $h'=-\la V(f),f'\ra$,
so that
\begin{equation}
\label{eq:h-decreix} |h'(t)+V(f(t))|^2\leq \eta|V(f(t)|.
\end{equation}
This implies the following: for every $\delta>0$ there exists some
$K>0$ such that, if $|\eta|\leq K/2$ and $f(t)\in X^{\delta}$ then
$h'(t)\leq K/2$.

\begin{lemma}
\label{lemma:mono-1} Given any small $\epsilon>0$, if $\eta$ is
small enough (depending on $\epsilon$) we have, for any $t$,
\begin{itemize}
 \item[(1)] if $f(t)\in X^{\epsilon M+\delta}$ and
 $h(t)\leq h_0+\epsilon K/4$ then, for every $t'\geq t+\epsilon$,
 we have $h(t')\leq h_0-\epsilon K/4$.
\item[(2)] if $h(t)\leq h_0-\epsilon K/4$, then for any $t'\geq t$
we have $h(t')\leq h_0-\epsilon K/4$;
 \item[(3)] if $h(t)\leq h_0+\epsilon K/4$, then for any $t'\geq t$
we have $h(t')\leq h_0+\epsilon K/4$.
\end{itemize}
\end{lemma}
\begin{pf}
Suppose that $\eta$ is so small that whenever $f(t)\in X^{\delta}$
we have $h'(t)\leq K/2$. It follows from
(\ref{eq:cota-derivada-f}) that if $f(t)\in X^{\epsilon
M+\delta}$, then for every $\tau\in[t,t+\epsilon]$ we have
$f(\tau)\in X^{\delta}$, so $h'(\tau)\leq K/2$. Integrating we
obtain (1). We now prove (2). If $\epsilon$ is small enough then
$V$ does not vanish in $\{H=h_0-\epsilon K/4\}$. Taking $\eta$
smaller than one half of the supremum of $|V|$ on the level set
$\{H=h_0-\epsilon K/4\}$ we deduce from (\ref{eq:h-decreix}) that
for any $\tau$ such that $f(\tau)=h_0-\epsilon K/4$ we have
$h'(\tau)<0$. This clearly implies (2). The same argument proves
(3).
\end{pf}

Take any $\epsilon$ satisfying the hypothesis of the previous
lemma and also $\epsilon\leq M^{-1}\delta$. Take $d<\delta$ small
enough so that we have an inclusion
$$F_0^d\subset \{|H-h_0|\leq \epsilon K/4\}.$$ Define also
$c:=\epsilon K/4$. We claim that this choice of $d$ and $c$
satisfies the requirements of the lemma. Indeed, suppose that for
some $t<t'$ we have $f(t)\in F_0^d$ and $f(t')\notin
F_0^{2\delta}$. Then $h(t)\leq h_0+\epsilon K/4$, so by (3) in
Lemma \ref{lemma:mono-1} we also have $h(t')\leq h_0+\epsilon
K/4$. If $\delta$ is small enough so that the
$2\delta$-neighbourhoods of each connected component of $F$ are
all disjoint, then we can assume that $f(t')\in X^{2\delta}\subset
X^{\epsilon M+\delta}$.  Combining this fact with the bound on
$h(t')$ and (1) in Lemma \ref{lemma:mono-1} we deduce that
$h(t'+\epsilon)\leq h_0-K/4$. Finally, (2) in Lemma
\ref{lemma:mono-1} implies that for every $t''\geq t'+\epsilon$ we
have $h(t'')\leq h_0-K/4$, which is what we wanted to prove
(recall that $M\geq 1$, so $\epsilon\leq M^{-1}\delta\leq
\delta$).

\subsection{Proof of Lemma \ref{lemma:recurrencia-dels-d}}
\label{ss:proof-recurrencia-dels-d}
Suppose that $f_u|_E$ takes values on the restriction $N'$ of $N$
to the connected component $F'\subset F$ of the fixed point set.
Let $N'=N_1\oplus\dots\oplus N_k$ be the decomposition in weights
of the $S^1$ action, and denote the corresponding weights
by $w_1,\dots,w_k$.
Let also $A$ denote the endomorphism of $N$ acting on the
subbundle $N_j$ as multiplication by $w_j$.

Define, for every $t\in[-1,2]$, $\phi(t):=f_u(n+t)$,
$G(t):=\phi'(t)-V(\phi(t))$ and $H(t):=V(\phi(t))-V_0(\phi(t))$.
Let $(x_1,\dots,x_k)$ be the coordinates of $\phi(0)$. The
relation between $\phi$ and the integral curves of the linear
vector field $V_0$ is given by Duhamel's formula:
\begin{align}
\phi(t) &= e^{tA}\phi(0)+\int_0^te^{(s-t)A}(G(s)+H(s))ds \notag \\
&=(e^{w_1t}x_1,\dots,e^{w_kt}x_k)+\int_0^te^{(s-t)A}(G(s)+H(s))ds.
\label{eq:duhamel}
\end{align}
Let $\phi_0(t):=(e^{w_1t}x_1,\dots,e^{w_kt}x_k)$ and let
$\chi(t):=\int_0^te^{(s-t)A}(G(s)+H(s))ds$. Define also
$t_j:=\int_{j-n}^{j-n+1}|\phi_0(t)| dt$ for every
$j\in\{n-1,n,n+1\}$. One checks easily (see for example the proof
of Lemma \ref{lemma:integrals-dexponencials}) that
\begin{equation}
\label{eq:recur-t} t_n\leq \frac{1}{e+e^{-1}}(t_{n-1}+t_{n+1}).
\end{equation}
By assumption we have $$\|G(t)\|_{L^1(I)}\leq
(g_{n-1}+g_n+g_{n+1})\leq \epsilon d_n.$$ On the other hand, using
(\ref{eq:V0-aproxima-V}) and $|\phi|<K\delta$ we obtain
$$\|H(t)\|_{L^1(I)}\leq K\left(\int_I |\phi|^2\right)\leq
K\delta \|\phi\|_{L^1(I)}.$$ Since $e^{(s-t)A}$ is bounded
for $s,t\in[-1,2]$, we can bound
$$\|\chi\|_{L^1(I)} \leq K(\|G(t)\|_{L^1(I)}+\|H(t)\|_{L^1(I)}) \\
\leq K(\epsilon d_n+\delta\|\phi\|_{L^1(I)}).$$ Taking this
estimate for $I$ equal to $[n-1,n]$, $[n,n+1]$ or $[n+1,n+2]$ it
follows from (\ref{eq:duhamel}) that, for any $j\in\{n_1,n,n+1\}$,
the inequality $|d_j-t_j| \leq K\delta (\epsilon d_n+\delta d_j)$
holds. Combining this estimate with (\ref{eq:recur-t}), and taking
into account that $\gamma>1/(e+e^{-1})$, we deduce that if
$\delta$ and $\epsilon$ are small enough then we must have
$$d_n\leq\gamma(d_{n-1}+d_{n+1}),$$
which is what we wanted to prove.

\end{document}